\documentclass[10pt]{article}

\usepackage[cp1252]{inputenc} 
\usepackage[T1]{fontenc}
\usepackage[a4paper]{geometry}
\usepackage[french,english]{babel}

\usepackage{amsmath, amssymb,mathrsfs,empheq}
\usepackage{amsthm}
\usepackage{graphicx,enumerate}
\usepackage{fancybox}
\usepackage{morefloats}
\usepackage{color}
\usepackage{cancel}
\usepackage{url}

\allowdisplaybreaks[4]

\theoremstyle{plain}
\newtheorem{prop}{Proposition}[section]
\newtheorem{thm}[prop]{Theorem}
\newtheorem{lem}[prop]{Lemma}

\theoremstyle{definition}

\theoremstyle{remark}
\newtheorem{rmq}{Remark}[section]

\newcommand{\R}{\mathbb{R}}
\newcommand{\Z}{\mathbb{Z}}
\newcommand{\N}{\mathbb{N}}
\newcommand{\dd}{\mathrm{d}}
\newcommand{\ep}{\varepsilon}
\newcommand{\alp}{\alpha}

\newcommand{\e}{{\mathrm{e}}}
\newcommand{\dt}{{\Delta t}}
\newcommand{\dx}{{\Delta x}}
\newcommand{\dy}{{\Delta y}}

\newcommand{\ds}{\displaystyle}
\newcommand{\tin}{{\mathrm{in}}}
\newcommand{\LL}{{L}}
\newcommand{\ccl}{[\![}
\newcommand{\ccr}{]\!]}

\renewcommand{\phi}{\varphi}
\renewcommand{\ge}{\geqslant}
\renewcommand{\le}{\leqslant}
\newcommand{\xstar}{{x^*}}
\newcommand{\xistar}{{\xi^*}}

\title{Concentration in Lotka-Volterra parabolic equations: an asymptotic-preserving scheme
}
\author{Vincent Calvez\footnote{Univ Lyon, CNRS, Université Claude Bernard Lyon 1, UMR5208, Institut Camille Jordan, F-69622 Villeurbanne, France. vincent.calvez\at math.cnrs.fr} , 
H\'el\`ene Hivert\footnote{Univ Lyon, \'Ecole centrale de Lyon, CNRS UMR 5208, Institut Camille Jordan, F-69134 \'Ecully, France. helene.hivert\at ec-lyon.fr} , 
Havva Yolda\c s\footnote{Faculty of Mathematics, University of Vienna, Oskar-Morgenstern-Platz 1, 1090 Vienna, Austria. havva.yoldas\at univie.ac.at}}
\date{}

\begin{document}

\selectlanguage{english}
\maketitle

\begin{abstract}
In this paper, we introduce and analyze an asymptotic-preserving scheme for  Lotka-Volterra parabolic equations. It is a class of nonlinear and nonlocal stiff equations, which describes the evolution of a population structured with phenotypic trait. In a regime of long time and small mutations, the population concentrates at a set of dominant traits. The dynamics of this concentration is described by a constrained Hamilton-Jacobi equation, which is a system coupling a Hamilton-Jacobi equation with a Lagrange multiplier determined by a constraint. This coupling makes the equation nonlocal. Moreover, the constraint does not enjoy much regularity, since it can have jumps.

The scheme we propose is convergent in all the regimes, and enjoys stability in the long time and small mutations limit.  Moreover, we prove that the limiting scheme converges towards the viscosity solution of the constrained Hamilton-Jacobi equation, despite the lack of regularity of the constraint. The theoretical analysis of the schemes is illustrated and complemented with numerical simulations.
\end{abstract}


\section{Introduction}

We are interested in the numerical analysis of a Lotka-Volterra parabolic equation
\begin{equation} 
 \label{eq:n_epsilon}
 \left\{ 
 \begin{array}{l l}
 \ds \partial_t n_\ep(t,x) - \ep \Delta_x n_\ep(t,x) = \frac{n_\ep(t,x)}{\ep} R(x,I_\ep(t)), & \ds x\in\R^d, t\ge 0 \vspace{4pt} \\
 I_\ep(t)=\int_{\R^d} \psi(x) n_\ep(t,x)\dd x, & \ds t\ge 0,
 \end{array}
 \right.
\end{equation}
supplemented with the initial condition $n_\ep(t=0,x)=n_\ep^{\tin}(x) \in\LL^1(\R^d)$, such that $n_\ep^\tin> 0$. 
It is a particular case of models arising in the theory of adaptative evolution \cite{MetzGeritzMesznaAl1996, GeritzMetzKisdiAl1997, GeritzKisdiMeszenaAl1998, MeszenaGyllenbergAl2005,Diekmann2004}.
It describes the evolution of a population structured with phenotypic trait, where $n_\ep(t,x)$ denotes the amount of individuals with trait $x\in\R^d$ at time $t\ge 0$. The evolution of the population is driven by births and deaths, synthesized in the net growth rate $R$. 
Note that the birth and death rates depend on the phenotypic trait, meaning that some individuals may be advantaged, because they are better adapted. The function $R$ also depends on $I_\ep$, defined in \eqref{eq:n_epsilon}, accounting for the total population burden on each individual growth rate. 
Mutations in the model are represented by unbiased random changes of phenotypes, with  the Laplacian term in the left-hand side of \eqref{eq:n_epsilon}. 
The parameter $\ep\in(0,1]$ in \eqref{eq:n_epsilon} is a scaling parameter,  so that considering the limit $\ep\to 0$ stands for the study of the population in an asymptotic regime of long time and small mutations. This is usually referred to as the separation of ecological and evolutionary time scales.

The asymptotic analysis of phenotype-structured population models has been carried out for various situations, we refer for instance to \cite{DiekmannJabinMischlerPerthame2005, PerthameTranspEqBio, CarrilloCuadradoPerthame2007, DesvillettesJabinMischlerRaoul2008, NordmannPerthameTaing2018, LorzLorenziClairambaultAl2015, LorenziPouchol2020}. The particular case of \eqref{eq:n_epsilon} was studied in \cite{BarlesPerthame, BarlesMirrahimiPerthame2009, LorzMirrahimiPerthame2011}, and more general mutation operators than the one in \eqref{eq:n_epsilon} were considered in \cite{BarlesPerthame2007, BarlesMirrahimiPerthame2009}. Generally speaking, because of the selection and the dynamics of adaptation, the population density is expected to concentrate on  a set of dominant traits,  meaning that it degenerates to a Dirac mass, or a sum of Dirac masses, located at the dominant trait(s). In particular, in the asymptotic regime $\ep\to 0$, the solution is expected to enjoy no better than measure regularity, requiring dedicated analytical methods. The Hopf-Cole transform, a logarithmic transformation of the unknown, is introduced to circumvent the regularity issues and study the dynamics of the concentration points. Coming back to \eqref{eq:n_epsilon}, the Hopf-Cole transform $u_\ep$ of $n_\ep$ is introduced 
\begin{equation}
 \label{eq:HopfCole}
 \forall\; t\ge 0,\;\forall\;x\in\R^d,\; n_\ep(t,x)=\e^{-u_\ep(t,x)/\ep},
\end{equation}
such that $u_\ep$ satisfies the following problem
\begin{equation}
\tag{$\mathrm{P}_\ep$}
 \label{eq:u_epsilon}
 \left\{ 
\begin{array}{l l}
\ds \partial_t u_\ep(t,x) + \left| \nabla_x u_\ep(t,x) \right|^2 = \ep \Delta_x u_\ep(t,x) - R(x,I_\ep(t)), & \ds x\in\R^d, t\ge 0, \vspace{4pt} \\ 
\ds I_\ep(t) = \int_{\R^d} \psi(x) \e^{-u_\ep(t,x)/\ep} \dd x, & \ds t\ge 0,
\end{array}
 \right.
\end{equation}
with the initial data $u_\ep(t=0,x)=u_\ep^\tin(x)= -\ep \ln n_\ep^\tin(x)$.

The asymptotic behavior of $u_\ep$ when $\ep\to 0$ is  studied in  \cite{BarlesMirrahimiPerthame2009}, under suitable assumptions on the parameters. Following \cite{BarlesMirrahimiPerthame2009}, we will suppose that  there are two constants $\psi_m$, $\psi_M$ such that 
\begin{equation}
 \label{assumption:psi} 
 \tag{A1}
\forall\; x\in\R^d,\; 0< \psi_m \le \psi(x) \le \psi_M < +\infty, \;\;\mathrm{and}\;\; \psi\in W^{2,\infty}(\R^d).
\end{equation}
It is also assumed that there exist two constants $0<I_m\le I_M <+\infty$ satisfying 
\begin{equation}
 \label{assumption:R_extrema}
 \tag{A2}
 \min\limits_{x\in\R^d} R(x,I_m) = 0, \;\;\;\; \max\limits_{x\in\R^d} R(x,I_M)=0,
\end{equation}
that $R$ is decreasing with respect to its second variable, and that there exists a constant $K>0$ such that 
\begin{equation}
\tag{A3}
 \label{assumption:R_bounded_decreasing} 
\forall x\in\R^d,\; \forall I\in\R,\; -K \le \partial_I R(x,I) \le -1/K <0, \;\;\mathrm{and}\;\; \sup\limits_{I_m/2\le I\le 2I_M} \left\| R(\cdot,I)\right\|_{W^{2,\infty}(\R^d)}\le K.
\end{equation}
In section \ref{sec:scheme_ep}, we will also use slightly stronger assumptions for $R$,
namely that  $I\mapsto \|R(\cdot,I)\|_{W^{2,\infty}(\R)}$ is bounded on all compact sets of $\R_+$.
The initial data $u_\ep^\tin$ in \eqref{eq:u_epsilon} is chosen such that 
\begin{equation}
\label{assumption:initial_data}
\tag{A4}
 \e^{-u_\ep^\tin/\ep}\in\LL^\infty(\R^d), \;\;\;\; I_m\le \int_{\R^d} \psi(x) \e^{-u_\ep^\tin(x)/\ep}\dd x\le I_M. 
\end{equation}
Because of assumption \eqref{assumption:initial_data},  $u_\ep^\tin$ has to be large at infinity. However, in what follows, a quantitative estimate of this behavior will be needed.  Still following \cite{BarlesMirrahimiPerthame2009}, we will then rather suppose that  
\begin{equation}
\label{assumption:u_increasing_infty}
\tag{A5}
 \exists\; \underline{a}, \overline{a}>0,\; \exists\; \underline{b},\overline{b}\in\R, \;\forall\; \ep >0,\; \forall x\in\R^d,\;  \underline{a}|x-x_0|+\underline{b}\le u^\tin_\ep(x)\le  \overline{a}|x-x_0|+\overline{b},
\end{equation}
where the upper bound is introduced for technical reasons, see Section \ref{sec:scheme_0}. 
Moreover, we will supose that $u^\tin_\ep$ enjoys Lipschitz regularity, uniformly with respect to $\ep>0$. Its Lipschitz constant is denoted by $L_0$,
\begin{equation}
 \label{assumption:u0Lipschitz}
 \tag{A6}
 \forall x,y\in\R^d,\; |u^\tin_\ep(x) - u^\tin_\ep(y)|\le L_0 |x-y|.
\end{equation}
Eventually, a refined assumption is made on the minimum of $u^\tin_\ep$, as we suppose that there exist two constants $c^\tin_m$ and $c^\tin_M$ such that 
\begin{equation}
 \label{assumption:min_u0}
 \tag{A7}
 c^\tin_m \ep \le \min u^\tin_\ep \le c^\tin_M \ep.
\end{equation}
Under these assumptions, the following theorem holds 
\begin{thm}[\cite{BarlesMirrahimiPerthame2009, CalvezLam2020}]
\label{thm:continuous}
Suppose that assumptions \eqref{assumption:psi}-\eqref{assumption:R_extrema}-\eqref{assumption:R_bounded_decreasing}-\eqref{assumption:initial_data}-\eqref{assumption:u_increasing_infty}-\eqref{assumption:u0Lipschitz}-\eqref{assumption:min_u0}
are satisfied. Let $u_\ep$ be the solution of \eqref{eq:u_epsilon} and $I_\ep$ be defined in \eqref{eq:u_epsilon}. Suppose also that $(u^\tin_\ep)_\ep$ is a sequence of uniformly continuous functions which converges locally uniformly to $v^\tin$. Then, $(u_\ep)_\ep$ converges locally uniformly to a function $v\in\mathcal{C}([0,+\infty[\times\R^d)$, and $(I_\ep)_\ep$ converges almost everywhere to a function $J$, such that $J\in BV(0,T)$ for all $T>0$, and that $(v,J)$ is the unique viscosity solution of the following equation
\begin{equation}
\tag{$\mathrm{P}_0$}
 \label{eq:limit}
 \left\{ 
 \begin{array}{l l}
  \ds \partial_t v(t,x) + \left| \nabla_x v(t,x) \right|^2 = - R(x,J(t)), & \ds x\in\R^d, \; t>0 \vspace{4pt} \\
 \ds \min\limits_{x\in\R^d} v(t,x) = 0 & \ds t>0,
 \end{array}
 \right. 
\end{equation}
with initial data $v^\tin$. 
\end{thm}

Equation \eqref{eq:limit} is a constrained Hamilton-Jacobi equation, with quadratic Hamiltonian 
\begin{equation}
\label{eq:Hamiltonian_continuous}
 \forall p\in\R^d, \;\mathcal{H}(p)=|p|^2,
\end{equation}
where $|\cdot|$ stands for the Euclidean norm on $\R^d$,
 and the unknown $J$ behaves as a Lagrange multiplier regarding the constraint $\min v(t,\cdot) = 0$. 
In Theorem \ref{thm:continuous}, uniqueness of the pair $(v,J)$ holds true in the class of locally Lipschitz-continuous functions $v$, and locally BV functions $J$. On the one side, Lipschitz regularity is a natural setting for viscosity solutions of Hamilton-Jacobi equations \cite{BarlesBook, Evans, CrandallIshiiLions1992}. On the other side,  the limiting function  
$J$ may have jump discontinuities \cite{BarlesPerthame, PerthameTranspEqBio, BarlesMirrahimiPerthame2009}, so that BV is the appropriate functional space for well-posedness.  
The existence of a solution $(v,J)$ of \eqref{eq:limit} is a consequence of \cite{BarlesMirrahimiPerthame2009}, where it is obtained as the limit of the sequence $(u_\ep,I_\ep)_\ep$ of solutions of \eqref{eq:u_epsilon}, together with locally uniform Lipschitz and BV estimates, respectively. 
The uniqueness of the pair $(v,J)$ has been adressed in some particular cases in \cite{BarlesPerthame, MirrahimiRoquejoffre2016, Kim2021}, then in  \cite{CalvezLam2020} in a more general setting including the problem under study. 
It is in fact composed of two companion results. Considering \eqref{eq:limit}, the following holds
\begin{thm}[\cite{CalvezLam2020}]
\label{thm:uniqueness}
\begin{enumerate}
 \item \label{thm:uniqueness_constrained} Suppose that $(v_1,J_1)$ and $(v_2,J_2)$ are two solutions of \eqref{eq:limit} in $W_{\mathrm{loc}}^{1,\infty}\times BV_{\mathrm{loc}}$ with the same initial data $v^\tin$. Assume that $v^\tin$ is coercive, that $\min v^\tin = 0$, and that $R$ is uniformly decreasing with respect to its second argument \eqref{assumption:R_bounded_decreasing}. Then, $v_1=v_2$,  and $J_1=J_2$ almost everywhere.
 \item \label{thm:uniqueness_unconstrained} Let $J\in BV(0,T)$ be given. Then, the variational solution $v$ of 
 \begin{equation}
 \label{eq:HJ}
  \partial_t v(t,x) + |\nabla_x v(t,x)|^2 = -R(x,J(t)),\;\;t>0, \;x\in\R^d,
 \end{equation}
with initial data $v^\tin$, is the unique locally Lipschitz viscosity solution of \eqref{eq:HJ} over $(0,T]\times\R^d$. Moreover, $v$ is independent of the choice of a representative of $J$ in $BV$. Namely, if \eqref{eq:HJ} is considered with two source terms $J_1$ and $J_2$ in $BV(0,T)$ such that $J_1=J_2$ almost everywhere in \eqref{eq:HJ}, then $v_1=v_2$. 
\end{enumerate}
\end{thm}

Theorem \ref{thm:uniqueness} suggests the following argument, which will be a key strategy in the present work. It is possible to consider $J$ being given as a source term in \eqref{eq:limit}, and show separately that the solution satisfies the constraint $\min v(t,\cdot)=0$, in order to prove that $(v,J)$ is the unique viscosity solution of \eqref{eq:limit}. This enables decoupling the Hamilton-Jacobi equation from its constraint.

In this paper, we propose and investigate a numerical scheme for \eqref{eq:u_epsilon} which enjoys stability properties when the parameter $\ep$ goes to $0$. Indeed, because of the definition of $I_\ep$ in \eqref{eq:u_epsilon}, the problem becomes stiff in the small-$\ep$ regime. If no specific strategy was employed, the accuracy of the numerical approximation of \eqref{eq:u_epsilon} would hence be deteriorated in the asymptotic regime. Schemes specifically designed for such singular problems are called \emph{Asymptotic-Preserving} (AP). They were introduced for kinetic equations  \cite{Jin, Klar1, Klar2}, and their properties are usually summarized by the following diagram 
\[
 \begin{array}{r c l}
  \eqref{eq:u_epsilon} &\xrightarrow{\;\;\;\ds\ep\to 0\;\;\;} & \eqref{eq:limit}  \vspace{4pt} \\
  \left. \begin{array}{c} $ $ \\ h\to 0 \\ $ $ \end{array} \right\uparrow \hspace{7pt} &  & \hspace{7pt}
  \left\uparrow  \begin{array}{c} $ $ \\ h\to 0 \\ $ $ \end{array} \right. \vspace{4pt}
  \\
  \left(S_\ep^h\right) & \xrightarrow[\;\;\;\ds\ep\to 0\;\;\;]{} &  \left(S_0^h\right)
 \end{array},
 \]
that should be understood as follows: an equation \eqref{eq:u_epsilon} depending on a parameter $\ep>0$ is given, and its solution converges when $\ep\to 0$ to the solution of another equation \eqref{eq:limit}. The scheme $\left(S_\ep^h\right)$, where all the discretization parameters are included in the notation $h$,  enjoys the AP property if it converges to the solution of \eqref{eq:u_epsilon} when $\ep>0$ is fixed and $h\to 0$, and if its solution converges when $h>0$ is fixed and $\ep\to 0$, to the solution of another  scheme $\left(S_0^h\right)$. The latter scheme is required to be convergent to the solution of problem \eqref{eq:limit}, when $h\to 0$. Even if it is in general not true, an AP scheme can also enjoy the stronger property of being \emph{Uniformly Accurate} (UA), meaning that its precision is independent of $\ep$. There is a large literature about AP schemes for various asymptotics of kinetic equations \cite{JinReview2012, DimarcoPareschiReview2014}, but, to the best of our knowledge, there are few results in case the asymptotic problem belongs to the class of  Hamilton-Jacobi equations: a scheme for front propagation in a one-dimensional kinetic linear BGK equation is analyzed in \cite{Hivert_2}, a scheme for dynamics of concentration in a selection-mutation equation close to \eqref{eq:u_epsilon} but with an integral mutation kernel is proposed, tested but not analyzed in \cite{CEMRACS2018}, and a model structured with age but where mutations are not considered is treated in \cite{AlmeidaPerthameRuan2020}. In contrast with AP schemes designed for linear kinetic equations, the latter works share the following features: the nonlinear character of the continuous problem \eqref{eq:u_epsilon}, and the need of a specific numerical analysis for the approximation of Hamilton-Jacobi equations \eqref{eq:limit}. 

The discretizations of the two problems \eqref{eq:u_epsilon} and \eqref{eq:limit} raise  several challenges. 
Concerning \eqref{eq:u_epsilon}, the stiffest term $I_\ep$ is handled implicitly in the numerical approximation. It implies stability in the small $\ep$ limit, but it requires the resolution of a nonlinear scalar equation, whose cost is independent of $\ep$. The other terms are discretized according to the properties expected for the scheme in the limit $\ep\to 0$. 

The numerical analysis of the constrained Hamilton-Jacobi problem \eqref{eq:limit} is original, to the best of our knowledge. We identified two important difficulties: the unbounded character of the solution on the one hand, and  the lack of regularity of $J$ on the other hand. 
We propose a finite-difference scheme for \eqref{eq:limit}, which enjoys partial monotonicity properties. The classical Hamilton-Jacobi side of the problem is handled with a standard monotonic scheme compatible with the discrete maximum principle \cite{CrandallLions, Souganidis}.
The contribution issued from the constraint comes with a nonlinear scalar problem to solve. During this step, the monotonicity of $R$ with respect to its second argument is crucially used to handle the lack of regularity.
Thanks to this construction, the scheme enjoys strong stability properties even if it is nonlocal, nonlinear, and it is used for unbounded data. 

The paper is organized as follows: the scheme for \eqref{eq:u_epsilon} is constructed, in Section \ref{sec:results}, as well as the scheme for the limit problem \eqref{eq:limit}. 
The AP property of the scheme for \eqref{eq:u_epsilon} is proved in Section \ref{sec:AP}. The convergence of the scheme for \eqref{eq:limit} is proved in Section \ref{sec:scheme_0}, while the convergence of the scheme for \eqref{eq:u_epsilon} for a given positive $\ep>0$ is treated in Section \ref{sec:scheme_ep}. Finally, various properties of the schemes are illustrated and discussed via numerical tests in Section \ref{sec:tests}.

\medskip

\noindent\textbf{Acknowledgment.} The authors wish to thank Beno\^it Gaudeul for the proofreading of this paper. 
This project has received funding from the European Research Council (ERC) under the European Union's Horizon $2020$ research and innovation program (ERC consolidator grant WACONDY n\textsuperscript{o} $865711$).
HY was partially supported by the Vienna Science and Technology Fund (WWTF) with a Vienna Research Groups for Young Investigators project, grant VRG17-014 (since October 2021). 
The third author would like to thank the Isaac Newton Institute for Mathematical Sciences for support and hospitality during the programme ``Frontiers in kinetic theory: connecting microscopic to macroscopic scales - KineCon 2022'' when work on this paper was undertaken. This work was supported by EPSRC Grant Number EP/R014604/1.

\section{Construction of the scheme and main results}
\label{sec:results}

In this section, we present the construction of an AP scheme  for \eqref{eq:u_epsilon} in dimension $d=1$,
and we state its properties.  Presenting the results in dimension $1$ avoids useless technical complications in what follows. However, the scheme  can be  generalized to any finite dimension, and its properties can be proved as in dimension $1$. The generalization of the scheme in higher dimension is presented in  Section \ref{sec:dim2}.

Let $T>0$ be fixed, the number $N_t$ of time steps be given. The time step is defined as $\dt = T/N_t$, and let $t_n=n\dt$ for $n\in\ccl 0,N_t\ccr$. The trait step is denoted $\dx>0$, and the grid is defined with $x_i= x_0+i\dx$ for a given $x_0\in\R$ and for all $i\in\Z$.  
For $n\in\ccl 0,N_t-1\ccr$ and $i\in\Z$, the scheme for \eqref{eq:u_epsilon} is given by 
\begin{equation}
\tag{$\mathrm{S}_\ep$}
 \label{scheme:epsilon}
 \left\{ 
\begin{array}{l}
 \ds \frac{u^{n+1}_i - u^n_i }{\dt } + H\left( \frac{u^n_i-u^n_{i-1}}{\dx}, \frac{u^n_{i+1}-u^n_i}{\dx} \right) = \ep \frac{u^n_{i+1}-2u^n_i + u^n_{i-1}}{\dx^2} - R(x_i,I^{n+1}) \vspace{4pt} \\
 I^{n+1}= \dx \sum\limits_{i\in\Z} \psi(x_i)\e^{-u^{n+1}/\ep}.
\end{array}
 \right.
\end{equation}
Note that sequences $(u^n_i)_{n,i}$ and $(I^{n+1})_n$ depend on $\ep$, although it is ommitted to simplify the notation. 
The scheme is initialized with $u^0_i=u_\ep^\tin(x_i)$ for all $i\in\Z$. 
The function $H$ is given by 
\begin{equation}
 \label{eq:H}
 H(p,q)=\max \left \{H^+(p),H^-(q) \right \},
 \end{equation}
 with
 \begin{equation}
 \label{eq:H+H-}
 H^+(p)=\left\{ \begin{array}{l l}\ds p^2 & \;\ds\mathrm{if}\; p>0  \vspace{4pt}\\ \ds 0 & \ds \mathrm{otherwise}, \end{array}\right.  
 \;\;\mathrm{and}\;\;
 H^-(q)=\left\{ \begin{array}{l l}\ds q^2 & \;\ds\mathrm{if}\; q< 0  \vspace{4pt}\\ \ds 0 &\ds \mathrm{otherwise}. \end{array}\right.  
\end{equation}
Such a choice of discretization for the Hamiltonian $\mathcal{H}$ defined in \eqref{eq:Hamiltonian_continuous} makes the scheme \eqref{scheme:epsilon} enjoy monotonicity properties. It is a classical assumption in numerical schemes for Hamilton-Jacobi equations, see \cite{CrandallLions,Souganidis}, and discretizations like \eqref{eq:H} were for instance used in \cite{GuerandKoumaiha2019}. Here, together with the implicit definition of $I^{n+1}$ in \eqref{scheme:epsilon}, it provides stability properties in the small $\ep$ limit. Moreover, we will show that the monotonicity is conserved when $\ep\to 0$. It is a key ingredient of the convergence of the scheme in the asymptotic regime.

In what follows, we will denote, for a given $L>0$,   
\begin{align}
\label{eq:CH_L}
C_H(L)
&=\sup\limits_{|p|\le L}|(H^+)'(p)|+\sup\limits_{|q|\le L}|(H^-)'(q)|
 = 4L.
\end{align}
Then, the following results hold:

\begin{prop}[Convergence of the scheme \eqref{scheme:epsilon}]
\label{prop:scheme_ep}
 Suppose that assumptions \eqref{assumption:psi}-\eqref{assumption:R_extrema}-\eqref{assumption:R_bounded_decreasing}-\eqref{assumption:u_increasing_infty}-\eqref{assumption:u0Lipschitz} are satisfied, and that $\ep>0$ and $T>0$ are fixed. Let $\Lambda\in (0,1)$.
 There exists $I_{M'}>0$, and  $\dt_0>0$ such that for all $\dt<\dt_0$ and $\dx$ satisfying
 \begin{equation}
  \label{eq:epfixed_CFL} 
  \tag{$\text{CFL}_\ep$}
  2\ep \frac{\dt}{\dx^2}+ C_H(L_0+T\kappa) \frac{\dt}{\dx} =\Lambda,
 \end{equation}
with $L_0$ defined in \eqref{assumption:u0Lipschitz},  $C_H$  in \eqref{eq:CH_L}, and 
 \begin{equation}
  \label{eq:kappa}
   \kappa=\sup\limits_{0\le I\le I_{M'}} \| R(\cdot, I)\|_{W^{2,\infty}(\R)},
  \end{equation}
 scheme \eqref{scheme:epsilon} is well defined. Moreover, there exists a constant $C(\ep)$, depending on $T$, $\|\partial_t^2 u_\ep\|_{\infty, [0,T]\times\R}$, $\|\partial_x^k u_\ep \|_{\infty,[0,T]\times\R}$ for $k=1,2,3$, and $\left\|\partial_x^2 \left(\psi n_\ep\right)\right\|_{\infty,[0,T]\times\R}$, such that for all $n\in\ccl0,N_t-1\ccr$, 
 \begin{equation}
 \label{eq:epfixed_CVu}
  \sup\limits_{i\in\Z} |u^{n+1}_i - u_\ep(t_{n+1},x_i)|\le C(\ep)(|\ln(\dt)|\dt+\dx),
 \end{equation}
 and
 \begin{equation}
 \label{eq:epfixed_CVI}
  \left|I_\ep(t_{n+1}) - I^{n+1}\right| \le C(\ep) \left(|\ln(\dt)|\dt + \dx\right),
 \end{equation}
 where $u_\ep$ and $I_\ep$ are defined in \eqref{eq:u_epsilon}, $u^{n+1}=(u^{n+1}_i)_{i\in\Z}$ and $I^{n+1}$ in \eqref{scheme:epsilon}, $n_\ep$ in \eqref{eq:n_epsilon} and $\psi$ in \eqref{assumption:psi}.
\end{prop}

\begin{rmq} It is worth remarking that the $\LL^\infty$ norms of derivatives of $u_\ep$ and $n_\ep$ in Prop. \ref{prop:scheme_ep} are well defined, provided that $\psi$ is smooth enough. Indeed, the bound for $\|\partial_x u_\ep\|_{\infty,[0,T]\times\R}$ is a consequence of the Lipschitz property of $u_\ep$ in $x$, and comes from the  maximum principle applied to \eqref{eq:u_epsilon} derivated with respect to $x$. Bounds for higher order derivatives, as well as derivatives of $n_\ep$, are consequences of Duhamel's formula for \eqref{eq:n_epsilon} and \eqref{eq:u_epsilon}, and of regularizing effects of the Laplacian.  As it is not the purpose of this paper, we omit the details of these properties.  One can refer to \cite{Evans} for the necessary tools.  
\end{rmq}

\begin{rmq} The estimate in $|\ln(\dt)|\dt$ in \eqref{eq:epfixed_CVu}-\eqref{eq:epfixed_CVI} comes from the quadrature rule in the approximation of $I_\ep$. At first sight, this could be seen as a reduction of order of the scheme, compared to the order $1$ in $\dt$ that could be expected. However, because of \eqref{eq:epfixed_CFL}, one has $\dx=_{\dt\to 0}\mathcal{O}(\sqrt{\dt})$, so that the order reduction in time has no impact on the precision of scheme \eqref{scheme:epsilon}. 
\end{rmq}

\begin{rmq}
\label{rmq:UA}
  The behavior of $C(\ep)$ when $\ep$ goes to $0$ brings serious difficulties. Indeed, it goes to $+\infty$ when $\ep\to 0$, meaning that the time step $\dt$ should be refined according to $\ep$ to make \eqref{scheme:epsilon} approximate \eqref{eq:u_epsilon} properly. The asymptotic behavior of $C(\ep)$ for small $\ep$ does not only come from regularity issues of $u_\ep$ when $\ep\to 0$. Indeed, it is strongly related to the fact that Prop. \ref{prop:scheme_ep} holds for fixed $\ep>0$ only. In particular, the constant $\kappa$ in \eqref{eq:epfixed_CFL} depends on $\ep$ and may go to $+\infty$ when $\ep$ goes to $0$. 
 To overpass this difficulty,  Prop. \ref{prop:scheme_ep} is supplemented by the two forthcoming propositions, that give the behavior of \eqref{scheme:epsilon} when $\ep$ is small. \end{rmq}

\begin{rmq} \label{rmq:wellposedness} Since the scheme \eqref{scheme:epsilon} is a coupled system of two implicit equations,
a nonlinear equation has to be solved to compute $(u^{n+1}_i)_{i\in\Z}$ and $I^{n+1}$. 
The fact that $I^{n+1}$ is well-defined is straightforward. 
 Indeed, it is solution of the equation $\phi(I)=0$, where 
 \begin{equation}
 \label{eq:def_Inp1}
  \phi(I)=I-\dx \sum\limits_{i\in\Z} \psi(x_i)\e^{-\tilde{u}^{n+1}_i/\ep} \e^{\dt R(x_i,I)/\ep},
 \end{equation}
 with 
 \[
  \tilde{u}^{n+1}_i \ds = u^n_i +\ep\dt \frac{u^n_{i+1}-2u^n_i + u^n_{i-1}}{\dx^2}- \dt H\left( \frac{u^n_i-u^n_{i-1}}{\dx}, \frac{u^n_{i+1}-u^n_i}{\dx}  \right).
 \]
  It is worth remarking that since $\phi$ is a difference between an increasing and a decreasing function, there exists an unique $I^{n+1}\in\R$ such that $\phi(I^{n+1})=0$. This property is independent of $\ep$, therefore the scheme \eqref{scheme:epsilon} is well-defined for all $\ep\in(0,1]$.
In practice, $I^{n+1}$ is computed first, with Newton's method for $\phi$. However, it must be implemented with care, to ensure that it is properly solved for all $\ep\in(0,1]$, with constant computational cost. 
The solution of equation \eqref{eq:def_Inp1} is uniformly bounded with respect to $\ep$. 
Indeed, we prove in Section \ref{sec:AP} that it is bounded by $2I_M$ when $\ep$ is small enough, with $I_M$ defined in  \eqref{assumption:R_extrema}, and a bound is given in Section \ref{sec:scheme_ep} for larger $\ep$, see Remark \ref{rmq:bound_I}.
However, \eqref{eq:def_Inp1} is very sensitive to approximations in the arguments of the exponentials, that are dramatically increased when $\ep$ is small. As a consequence, the numerical resolution of \eqref{eq:def_Inp1} can collapse during Newton's iterations. To avoid such a phenomenon, Newton's iterations are computed as analytically as possible, and implemented with special care of the compensations between terms. When it is not enough, $\phi(I)=0$ with $\phi$ defined in \eqref{eq:def_Inp1} is replaced by the equivalent equation 
\[
 \ln(I) = \ln(\dx) + \ln\left( \sum\limits_{i\in\Z} \psi(x_i) \e^{-\tilde{u}^{n+1}_i/\ep} \e^{\dt R(x_i,I)/\ep} \right),
\]
that is also solved with Newton's method. We refer to \cite{CEMRACS2018}, and to the codes available at \cite{Git}, 
for more details.
 \end{rmq}

 \begin{rmq} \label{rmq:truncation_ep}
  Since it is defined for indices $i\in\Z$, the scheme \eqref{scheme:epsilon} cannot be implemented exactly as it is defined. However, Prop. \ref{prop:scheme_ep} also holds for a truncated version of the scheme \eqref{scheme:epsilon}, in which the sum defining $I^{n+1}$ is considered on a finite number of indices.  
 We refer to Section \ref{sec:scheme_ep} for details.
 \end{rmq}

\begin{prop}[Convergence of the scheme \eqref{scheme:epsilon} to the scheme \eqref{scheme:limit}]
\label{prop:AP}
Under assumptions
\eqref{assumption:psi}-\eqref{assumption:R_extrema}-\eqref{assumption:R_bounded_decreasing}-\eqref{assumption:u_increasing_infty}-\eqref{assumption:u0Lipschitz}-\eqref{assumption:min_u0}, and supposing that $\dt$ and $\dx$ are fixed, such that
\begin{equation}
 \label{eq:CFL_AP}
 \tag{$\mathrm{CFL}_{\ep\to 0}$}
  2\ep \frac{\dt}{\dx^2}+ C_H(L_0+TK) \frac{\dt}{\dx}\le 1,
\end{equation}
is satisfied for all $\ep\in(0,1]$, where $L_0$ is defined in \eqref{assumption:u0Lipschitz}, $K$ in \eqref{assumption:R_bounded_decreasing} and $C_H$ in \eqref{eq:CH_L}.
Let  
$(u^{n+1}_i)_{n, i}$
and 
$(I^{n+1})_{n}$ 
be the $\ep$-dependent sequences defined by \eqref{scheme:epsilon}.  Then, for all $n\in\ccl 0,N_t-1\ccr$ and for all $i\in\Z$, 
\[
 u^{n+1}_i\underset{\ep\to 0}\longrightarrow v^{n+1}_i, \;\; I^{n+1} \underset{\ep\to 0}\longrightarrow J^{n+1} 
\]
where the sequences 
$(v^{n+1}_i)_{n,i}$  
and 
$(J^{n+1})_{n}$
satisfy the scheme
\begin{equation}
 \tag{$\mathrm{S}_0$}
 \label{scheme:limit}
 \left\{ 
\begin{array}{l l}
\ds \frac{v^{n+1}_i-v^n_i}{\dt}+ H\left(\frac{v^n_i-v^n_{i-1}}{\dx},\frac{v^n_{i+1}-v^n_i}{\dx}\right) = - R(x_i,J^{n+1}), & \;  n\in \ccl 0,N_t-1\ccr,\; i\in\Z \vspace{4pt} \\
\ds\min\limits_{i\in\Z} v^{n+1}_i=0, & \; n\in\ccl 0,N_t-1\ccr,
\end{array}
 \right.
\end{equation}
initialized with $v^0_i=v^\tin(x_i)$, for all $i\in\Z$. 
\end{prop}

\begin{rmq}
 As in Remark \ref{rmq:truncation_ep}, Prop. \ref{prop:AP} also holds for the truncated scheme that is implemented in practice.
\end{rmq}

\begin{rmq}
The well-posedness of \eqref{scheme:limit} is a consequence of Prop. \ref{prop:AP}. Indeed, the convergence of $u^{n+1}_i$ and $I^{n+1}$ when $\ep\to 0$ gives the existence of a solution of the implicit scheme \eqref{scheme:limit}.
The fact that $(v^{n+1}_i)_{i\in\Z}$ and $J^{n+1}$ are uniquely defined follows from the proof of Prop. \ref{prop:AP}. Discussion about the direct implementation of \eqref{scheme:limit} is postponed to Section \ref{sec:scheme_0}. 
\end{rmq}

The next proposition states that the solution of the scheme \eqref{scheme:limit} converges to the solution of the limit equation \eqref{eq:limit} when the discretization parameters go to $0$.
To this end, we extend the definition of the scheme \eqref{scheme:limit}, in order to make it coincide at  the grid points with a function defined over $[0,T]\times\R$, and we reformulate it, so that the monotonic component of the scheme is taken apart. 
It can be seen as an operator, denoted by $\mathcal{M}_s$, acting on functions defined on $\R$. Namely, for all $s\in(0,\dt]$ and $f:\R\rightarrow\R$, $\mathcal{M}_s(f):\R\rightarrow\R$ is defined by
\begin{equation}
 \label{scheme:Monotonic}
 \tag{$\mathcal{M}_s$}
 \forall x\in\R, \;\mathcal{M}_s (f)(x) = f(x) - s\; H\left( \frac{f(x)-f(x-\dx)}{\dx},\frac{f(x+\dx)-f(x)}{\dx} \right).
\end{equation}
Suppose now that the ratio $\lambda=\dt/\dx$ is fixed. Let us define 
$(t,x)\mapsto v_\dt(t,x)$ on $[0,T]\times\R$,  
and $t\mapsto {J}_\dt(t)$ on $(0,T]$, such that 
 for all $n\in\ccl 0,N_t-1\ccr$, $s\in(0,\dt]$, 
 and $x\in\R$,
\begin{subequations}
\label{scheme_0:MonotonicFormulation}
\begin{empheq}[left=\empheqlbrace]{alignat=3}
 \label{scheme_0:MonotonicFormulation_v} 
& \ds v_\dt(t_n+s, x) = \mathcal{M}_s \left( v_\dt(t_n,\cdot)\right)(x) - s R\left (x, {J}_\dt(t_n+s)\right) 
\vspace{4pt}\\
&\ds J_\dt(t_n+s)=J^{n+1} \vspace{4pt} \label{scheme_0:MonotonicFormulation_J} 
\\
& \ds \min\limits_{i\in\Z} v_\dt(t_{n+1},x_i) = 0, 
\label{scheme_0:MonotonicFormulation_min} 
 \end{empheq}
\end{subequations}
 and initialized with $v_\dt(0,\cdot)=v^\tin$.  
The function $J_\dt$ is piecewise constant, with $J^{n+1}$ defined in \eqref{scheme:limit}. 
It is easy to remark,
that  $v_\dt$ and ${J}_\dt$ coincide with the solution of the scheme \eqref{scheme:limit} at the grid points
\[\forall n\in\ccl 0,N_t\ccr, \;\forall i\in\Z, \;v_\dt(t_n,x_i)=v^n_i, \text{\;and\;} \forall n\in\ccl 0,N_t-1\ccr, \; {J}_\dt(t_{n+1})= J^{n+1}.\]
This is due to the fact that the constraint $\min v_\dt=0$ is only considered on the grid points in \eqref{scheme_0:MonotonicFormulation_min}.

For the sake of simplicity, let us denote
\begin{equation}
\label{eq:CH}
\mathcal{C}_H=C_H(14(L_0+KT)+1),
\end{equation} 
in what follows, where $C_H$ is defined in \eqref{eq:CH_L}, $L_0$ in \eqref{assumption:u0Lipschitz}, and $K$ in \eqref{assumption:R_bounded_decreasing}.

\begin{prop}[Convergence of the scheme \eqref{scheme:limit}]
\label{prop:scheme_0}
 Suppose that the assumptions of Theorem \ref{thm:continuous} are satisfied, and that $v^\tin$ satisfies \eqref{assumption:u_increasing_infty}-\eqref{assumption:u0Lipschitz}-\eqref{assumption:min_u0} for $\ep=0$. 
 Suppose that the ratio $\dt/\dx$ is fixed such that
 \begin{equation}
 \tag{$\mathrm{CFL}_0$}
 \label{eq:CFL_0}
\mathcal{C}_H\frac{\dt}{\dx}\le 1, \;\;\;\mathrm{and}\;\;\; \frac{\dt}{\dx}\sqrt{(L_0+TK)^2+K} \le 1,
 \end{equation}
 with $\mathcal{C}_H$ defined in \eqref{eq:CH}, $L_0$ in \eqref{assumption:u0Lipschitz}, and $K$ in \eqref{assumption:R_bounded_decreasing}.
 Then for all $t\in(0,T]$ and for all $x\in\R$, 
 \[
  |v_\dt(t,x)-v(t,x)|\underset{\dt\to 0}\longrightarrow 0,
\]
and the convergence is locally uniform on $(0,T]\times\R$.
Moreover, for almost all $t\in(0,T]$, 
\[
 |J_\dt(t) - J(t)| \underset{\dt\to 0}\longrightarrow 0,
\]
where $v$ and $J$ are uniquely determined as  the viscosity solution of \eqref{eq:limit}, and with $v_\dt$, $J_\dt$ defined by \eqref{scheme_0:MonotonicFormulation}.
\end{prop}

\begin{rmq}
 Remark that condition \eqref{eq:CFL_0} contains two items. Although they both show linear relations between $\dt$ and $\dx$, they are of very different nature. Indeed, the first one is a classical stability condition, which yields in particular the monotonicity of scheme \eqref{scheme:Monotonic}, with $\mathcal{C}_H$ in \eqref{eq:CH} taken a few larger than necessary for technical reasons.  On the other hand, the second condition makes $J_\dt$ nondecreasing. This is crucial in the compactness argument used to prove that $J_\dt$ converges when $\dt\to 0$. We refer to Section \ref{sec:scheme_0} for details.  
\end{rmq}

\begin{rmq}
 Contrary to Prop. \ref{prop:scheme_ep},
 Prop. \ref{prop:scheme_0} does not give any convergence rate for scheme \eqref{scheme:limit}. This comes from the lack of regularity of the viscosity solution $v$ and $J$ of \eqref{eq:limit}. Indeed, $J\in BV(0,T)$, while 
 $v$ enjoys Lipschitz regularity in $[0,T]\times\R$. This property is a consequence of the definition of $v$ as the variational solution of \eqref{eq:limit}, but it is also obtained in Section \ref{sec:scheme_0}, where $v$ is shown to be a limit of Lipschitz functions.
 According to this observation, one can come back to Prop. \ref{prop:scheme_ep}, and remark that no uniform bound in $\ep$ is to be expected for $C(\ep)$, even if the estimates of the proof of Prop. \ref{prop:scheme_ep} were made sharper.
\end{rmq}

\begin{rmq}
 Scheme \eqref{scheme:limit} can be implemented on a truncated domain, that is reduced at each time step, but with no more approximation. Hence, Prop. \ref{prop:scheme_0} holds for the scheme that is implemented in practice.
\end{rmq}

\section{Convergence of \eqref{scheme:epsilon} to the limiting scheme \eqref{scheme:limit}}
\label{sec:AP}

In this section, we prove that \eqref{scheme:epsilon} enjoys stability properties with respect to $\ep\in(0,1]$, thus Prop. \ref{prop:AP} follows. Prop. \ref{prop:AP} states that, when $\ep$ goes to $0$ with fixed discretization parameters, the solution of \eqref{scheme:epsilon} converges to the solution of \eqref{scheme:limit}. It
 relies on a convenient reformulation of the scheme \eqref{scheme:epsilon}, for all $n\in\ccl 0,N_t-1\ccr$ and for all $i\in\Z$

\begin{subequations}
\label{scheme_ep:MonotonicFormulation}
\begin{empheq}[left=\empheqlbrace]{alignat=2}
 \ds  u^{n+1}_i & \ds =  
 \mathcal{M}^\ep_\dt\left(u^n\right)_i
 - \dt R(x_i,I^{n+1}) \vspace{4pt}
 \label{scheme_ep:MonotonicFormulation_u}
 \\
  \ds I^{n+1}&\ds =\dx \sum\limits_{i\in\Z} \psi(x_i) \e^{-u^{n+1}_i/\ep},
  \label{scheme_ep:MonotonicFormulation_I}
 \end{empheq}
\end{subequations}
where $u^n=(u^n_i)_{i\in\Z}$, and $\mathcal{M}^\ep_\dt\left(u^n\right)\in\R^\Z$ is defined for all $i\in\Z$ by
\begin{equation}
 \mathcal{M}^\ep_\dt\left(u^n\right)_i= 
  u^n_i +\ep\dt \frac{u^n_{i+1}-2u^n_i + u^n_{i-1}}{\dx^2}- \dt H\left( \frac{u^n_i-u^n_{i-1}}{\dx}, \frac{u^n_{i+1}-u^n_i}{\dx}  \right).
 \label{scheme_ep:MonotonicFormulation_utilde}
\end{equation}
As it has been announced in Section \ref{sec:results}, the scheme \eqref{scheme:epsilon} enjoys monotonicity properties. More precisely, it is a consequence of  the first step \eqref{scheme_ep:MonotonicFormulation_utilde}. Indeed, the following properties hold (see \cite{CrandallLions}):

\begin{lem}
 \label{lem:Monotonic_ep} 
 Let $u=(u_i)_{i\in\Z}$ and $v=(v_i)_{i\in\Z}\in\R^\Z$, and $\mathcal{M}_\dt^\ep$ defined as in \eqref{scheme_ep:MonotonicFormulation_utilde}.
 Let $L>0$,
 and suppose that $2\ep \dt/\dx^2 + \dt C_H(L)/\dx\le 1$, with $C_H(L)$ defined in \eqref{eq:CH_L}. 
Then the following results hold true 
\begin{itemize}
\item 
If there exists $i\in\Z$ such that,   
$\left|  u_i-u_{i\pm1} \right| \le L\dx$, $\left|  v_i-v_{i\pm1} \right| \le L\dx$, 
and $\forall j\in\ccl i-1,i+1\ccr$, 
$u_j\le v_j$,
then $\mathcal{M}_\dt^\ep\left(u\right)_i \le \mathcal{M}_\dt^\ep\left(v\right)_i$. 
\item If for all $i\in\Z$, $|u_i-u_{i-1}|\le L\dx$, then for all $i\in\Z$, 
$ \left| \mathcal{M}_\dt^\ep\left(u\right)_i-\mathcal{M}_\dt^\ep\left(u\right)_{i-1} \right| \le L\dx. $
\item If $u-v=(u_i-v_i)_{i\in\Z}\in\ell^\infty(\Z)$, and if, for all $i\in\Z$, $|u_i-u_{i-1}|\le L\dx$, and $|v_i-v_{i-1}|\le L\dx$,  then $\mathcal{M}_\dt^\ep\left(u\right) - \mathcal{M}_\dt^\ep\left(v\right)\in\ell^\infty(\Z)$ and 
$
 \left\| \mathcal{M}_\dt^\ep\left(u\right)-  \mathcal{M}_\dt^\ep\left(v\right)\right\|_\infty \le\|u-v\|_\infty.
$
\end{itemize}
\end{lem}
Using this lemma and the reformulation \eqref{scheme_ep:MonotonicFormulation} of the scheme \eqref{scheme:epsilon}, stability properties of the scheme \eqref{scheme:epsilon} when $\ep$ goes to $0$ are proved. The following lemma is stated:

\begin{lem}
\label{lem:AP}
 Suppose that assumptions \eqref{assumption:psi}-\eqref{assumption:R_extrema}-\eqref{assumption:R_bounded_decreasing}-\eqref{assumption:u_increasing_infty}-\eqref{assumption:u0Lipschitz}-\eqref{assumption:min_u0} hold true, and  that $\dt$ and $\dx$ are fixed such that the inequality \eqref{eq:CFL_AP} is satisfied.  
 Then, there exists an $\ep_0>0$, depending only on the constants arising in the assumptions and on $\dx$ and $\dt$, such that for all $\ep\in(0,\ep_0)$, the sequence $(u^n_i)_{n,i}$ defined by the scheme \eqref{scheme:epsilon} satisfies:
 \begin{enumerate}
  \item Uniform Lipschitz continuity in trait: \label{lem:AP_Lipschitz} For all $n\in\ccl 0,N_t\ccr$, there exists a constant $L_n=L_0+n\dt K\le L_0+TK$, with $L_0$ defined in \eqref{assumption:u0Lipschitz} and $K$ in \eqref{assumption:R_bounded_decreasing}, such that the sequence $u^n=(u^n_i)_{i\in\Z}$ enjoys $L_n$-Lipschitz property:
  \[
   \forall i\in\Z, \; \left| \frac{u^n_{i+1}-u^n_i}{\dx}\right| \le L_n.
  \]
  \item Uniform bound from below for $u^n$: \label{lem:AP_LowEstimate} For all $n\in\ccl 0,N_t\ccr$, there exists $\underline{b}_n\in\R$, such that
  $ \underline{b}_n\ge 
   \underline{b}_{N_t}=\underline{b}-TH(\underline{a},\underline{a})-TK$,
    and that for all $i\in\Z$, 
  \[ 
  u^{n}_i \ge   \underline{a}|x_i-x_0|+\underline{b}_n,
  \]
  where $\underline{a}$ and $\underline{b}$ have been defined in \eqref{assumption:u_increasing_infty}, $H$ in \eqref{eq:H}, $K$ in \eqref{assumption:R_bounded_decreasing}, and $T$ is the fixed final time. 
\item \label{lem:AP_I} Uniform bounds for $(I^n)_{n\in\ccl 1,N_t\ccr}$:
For all $n\in\ccl 0,N_t-1\ccr$,
\[
 I_m/2\le I^{n+1}\le 2I_M.
 \]
\item \label{lem:AP_minu} Estimate for $\min u^n$: There exist  $c_m$ and $c_M$ such that for all $n\in\ccl 0,N_t\ccr$,
\[
c_m \ep \le \min\limits_{i\in\Z} u^n_i \le c_M \ep,
\]
and $c_m\le c_m^\tin$ and $c_M\ge c_M^\tin$ depend only on the constants defined in the assumptions and on $\dx$ and $\dt$. 
 \end{enumerate}

\end{lem} 

\begin{proof}
We proceed by induction. 
 Thanks to the assumptions, the initial data $u^0=(u^0_i)_{i\in\Z}$ satisfies the properties \ref{lem:AP_Lipschitz}-\ref{lem:AP_LowEstimate} and \ref{lem:AP_minu} of Lemma \ref{lem:AP}.
 Let us suppose that the items \ref{lem:AP_Lipschitz}-\ref{lem:AP_LowEstimate}-\ref{lem:AP_minu} of Lemma \ref{lem:AP} hold true for a given $n\in\ccl 0,N_t-1\ccr$, and prove  that $u^{n+1}=(u^{n+1}_i)_{i\in\Z}$ enjoys these properties, while $I^{n+1}$ satisfies \ref{lem:AP_I}.
 
 \begin{itemize}
\item First of all, 
 we recall that $I^{n+1}$ is well defined for all $\ep\in(0,1]$, see Remark \ref{rmq:wellposedness}.
We now prove that, if $\ep\le \ep_0$, then $I^{n+1} \ge I_m/2$. 
 For a given $j\in\Z$, the following inequality holds 
 \[
  I-\dx \sum\limits_{i\in\Z} 
  \psi(x_i)\e^{-\mathcal{M}_\dt^\ep\left(u^n \right)_i/\ep} \e^{\dt R(x_i,I)/\ep}
  \le I-\dx \psi_m \e^{-
  \mathcal{M}_\dt^\ep\left( u^n \right)_j
  /\ep} \e^{\dt R(x_j,I)/\ep},
 \]
 where we used \eqref{assumption:psi}. 
An upper bound for  
$\mathcal{M}_\dt^\ep\left( u^n \right)_j$
is obtained  thanks to the positivity of $H$  and to  \ref{lem:AP_Lipschitz}
\begin{align*}
 \mathcal{M}_\dt^\ep\left( u^n \right)_j&= u^n_j + \ep\frac{\dt}{\dx} \frac{u^{n}_{j+1}-2u^n_i + u^n_{j-1}}{\dx} -\dt H\left( \frac{u^n_i-u^n_{i-1}}{\dx}, \frac{u^n_{i+1}-u^n_i}{\dx} \right)
 \\&\le u^n_j + 2\ep \frac{\dt}{\dx}L_n,
 \end{align*}
using the upper bound $L_0+TK$ of $L_n$, 
and with the choice of $j$ such that $u^n_j=\min_{i\in\Z} u^n_i$, property \ref{lem:AP_minu} provides
\[
 \mathcal{M}_\dt^\ep\left( u^n \right)_j
 \le \ep\left( c_M + 
 2\frac{L_0 +TK}{C_H(L_0+TK)}
 \right),
\]
where the estimate independent of $\dt$ and $\dx$ comes from \eqref{eq:CFL_AP}. 
This estimate yields
\[
 \phi(I) \le I - \dx \psi_m \e^{
 -c_M -2(L_0+TK)/C_H(L_0+TK)
 } \e^{\dt R(x_j,I)/\ep},
\]
and the bound from below for $I^{n+1}$ in \ref{lem:AP_I} is then obtained by a contradiction argument. 
Indeed, since $\phi$ is an increasing function, for all $I<I_m/2$, the following inequality holds true 
\begin{equation}
 \label{eq:prooflem_phi}
 \phi(I)\le \phi(I_m/2) \le \frac{I_m}{2} - \dx \psi_m \e^{
 -c_M -2(L_0+TK)/C_H(L_0+TK)
 } \e^{\dt R(x_j,I_m/2)/\ep},
\end{equation}
and 
thanks to the strict monotonicity of $R$ with respect to its second argument, together with \eqref{assumption:R_extrema}, one can show that $R(x_j,I_m/2)$ is uniformly positive with respect to $j\in\Z$.
Indeed, it writes
\[
 R(x_j,I_m/2) \ge R(x_j,I_m) + 
 \frac{ I_m}{2K}
 ,
\]
thanks to \eqref{assumption:R_bounded_decreasing}, and assumption \eqref{assumption:R_extrema} eventually yields 
\[
 R(x_j,I) 
  \ge \frac{I_m}{2K}
 .
\]
Coming back to \eqref{eq:prooflem_phi}, one has for all $I\le I_m/2$
\[
 \phi(I)\le \phi(I_m/2) \le \frac{I_m}{2} - \dx\psi_m 
 \e^{-c_M-2(L_0+TK)/c_H(L_0+TK)}\e^{\dt  I_m/{2K\ep}}
 \underset{\ep\to 0}\longrightarrow -\infty,
\]
hence there exists an $\ep_1>0$, depending only on $I_m$, $\psi_m$, $c_M$, $C_H$, $L_0$, $T$, $K$, $\dt$, and $\dx$, such that 
\[
 \forall \ep\in(0,\ep_1),\; \forall I\le I_m/2, \; \phi(I)\le -1.
\]
Since $I^{n+1}$ is defined as the solution of $\phi(I^{n+1})=0$, the first inequality in \ref{lem:AP_I} holds true.

\item The bound from below \ref{lem:AP_LowEstimate} of $(u^{n+1}_i)_{i\in\Z}$, is  a consequence of the monotonicity of 
the first step
of the scheme \eqref{scheme:epsilon}. Indeed, if we denote $v^n_i=  \underline{a}|x_i-x_0|+\underline{b}_n$, the scheme \eqref{scheme_ep:MonotonicFormulation_utilde} applied to $v^n=(v^n_i)_{i\in\Z}$ gives 
\[
 \mathcal{M}_\dt^\ep\left( v^n \right)_i
 = \left\{
 \begin{array}{l l}
 \ds  \underline{a} |x_i-x_0|+\underline{b}_n -\dt H(\underline{a},\underline{a}) & \ds \;\;\mathrm{if}\; i\neq 0 \vspace{4pt} \\
  \ds \underline{b}_n + 2\underline{a} \frac{\ep \dt}{\dx} & \ds \;\;\mathrm{if}\; i =0,
 \end{array}
 \right. 
\]
since $H(\underline{a},\underline{a})=H(-\underline{a},-\underline{a})$. Therefore,  
$\mathcal{M}_\dt^\ep\left( u^n \right)_i \ge \mathcal{M}_\dt^\ep\left( v^n \right)_i$,
for all $i\in\Z$, thanks to Lemma \ref{lem:Monotonic_ep},
so that
\[
 \left\{ 
\begin{array}{l l} 
\ds u^{n+1}_i\ge  \underline{a}|x_i-x_0|+\underline{b}_n -\dt H(\underline{a},\underline{a}) -\dt R(x_i,I^{n+1}) & \ds \;\;\mathrm{if}\; i \neq 0 
\vspace{4pt} \\
  \ds u^{n+1}_0 \ge \underline{b}_n + 2\underline{a} \frac{\ep \dt}{\dx} -\dt R(x_0,I^{n+1}). &
\end{array} 
 \right.
\]
Since $I^{n+1}\ge I_m/2$, and $I\mapsto R(x_i,I)$ is decreasing for all $i\in\Z$, 
the choice $\underline{b}_{n+1}= \underline{b}_n - \dt H(\underline{a},\underline{a}) - \dt K$ yields \ref{lem:AP_LowEstimate} thanks to \eqref{assumption:R_bounded_decreasing}.

\item The second inequality in \ref{lem:AP_minu} is a consequence of the bound from below of $(u^{n+1}_i)_{i\in\Z}$, as well as the one for $I^{n+1}$. Indeed, considering $\ep\in(0,\ep_1)$, the definition of $I^{n+1}$ yields
\[
 I_m/2\le I^{n+1} \le \dx \sum\limits_{i\in\Z} \psi(x_i) \e^{-u^{n+1}_i/\ep},
\]
and for an integer $N\ge 1$, which will be determined later, the following inequality holds true 
\[
 I_m/2\le \dx\sum\limits_{|i|< N} \psi(x_i)\; \e^{-u^{n+1}_i/\ep}+ \dx \sum\limits_{|i|\ge N} \psi(x_i)\; \e^{-\left(\underline{b}_{N_t}+ \underline{a} |x_i-x_0|\right)/\ep},
\]
because of \ref{lem:AP_LowEstimate}. In both terms, we use \eqref{assumption:psi}, and since $x_i-x_0=i\dx$, we have 
\[
 I_m/2 \le (2N-1)\dx\;\psi_M \;\e^{-\min\limits_{i\in\Z} u^{n+1}_i/\ep} + 2\dx\; \psi_M\; \e^{-\left(\underline{b}_{N_t}+\underline{a}N\dx\right)/\ep}\sum\limits_{i\ge 0} \e^{-\underline{a} i \dx/\ep}. 
\]
Therefore, $N$ is chosen such that $\underline{b}_{N_t}+\underline{a}N\dx\ge 1$. Note that this choice is independent of $n$, and that it depends only on the assumptions, and $\dx$. Hence, the previous inequality can be simplified as
\[
 I_m/2\le (2N-1) \dx\;\psi_M\; \e^{-\min\limits_{i\in\Z}u^{n+1}_i/\ep} + \frac{2\dx\;\psi_M \;\e^{-1/\ep}}{1-
 \e^{-\underline{a}\dx/\ep}
 },
\]
and $\ep_2$ can be defined, as a function of the parameters arising in the assumptions and of $\dx$, but independently of $n$ such that 
\[
 \forall \ep\in(0,\ep_2),\;  \frac{2\;\dx \;\psi_M\;\e^{-1/\ep}}{1-
 \e^{-\underline{a}\dx/\ep}
 }\le I_m/4,  
\]
so that for $0<\ep<\min(\ep_1,\ep_2)$, the second inequality in \ref{lem:AP_minu} is satisfied, with 
\[
c_M = \max \left\{ -\ln\left( \frac{I_m}{4(2N-1) \dx\psi_M}\right), c_M^\tin\right\}.
\]
Once again, it is worth noticing that this choice is independent of $n$.

\item The previous results yield the inequality $I^{n+1}\le 2I_M$ in \ref{lem:AP_I}. Indeed, thanks to \ref{lem:AP_minu}, $u^n_i\ge c_m \ep$ for all $i\in\Z$, and because of the monotonicity of \eqref{scheme_ep:MonotonicFormulation_utilde}, this implies that 
\[
 \forall i\in\Z, \; 
 \mathcal{M}_\dt^\ep\left( u^n \right)_i
 \ge c_m \ep. 
\]
Moreover the bound from below \ref{lem:AP_LowEstimate} satisfied by $u^{n+1}=(u^{n+1}_i)_{i\in\Z}$, ensures that there exists an index $k\in\Z$ such that $u^{n+1}_k=\min_{i\in\Z} u^{n+1}_i\le c_M \ep$. The definition of $u^{n+1}_i$ at line  \eqref{scheme_ep:MonotonicFormulation_u} yields
\[
 \dt R(x_k, I^{n+1}) = 
 \mathcal{M}_\dt^\ep\left( u^n \right)_k
 - u^{n+1}_k\ge -(c_M  - c_m)\ep \ge \dt R(x_k,I_M) -(c_M-c_m) \ep,
\]
where the last inequality has been obtained considering that $R(x_k,I_M)\le 0$, according to \eqref{assumption:R_extrema}.
One can conclude similarly as above, using \eqref{assumption:R_bounded_decreasing}  to write
\[
 \dt R(x_k,I^{n+1}) \ge \dt R(x_k,2I_M) +\frac{\dt I_M}{K}-\ep(c_M-c_m).
\]
Denoting $\ep_3 = \dt I_M / (K(c_M-c_m))$, the previous inequality states that for all $\ep\in (0,\ep_3)$, $R(x_k,I^{n+1})\ge R(x_k,2I_M)$. Once again, we emphasize the fact that $\ep_3$ is defined once for all since it depends only on the assumptions and on $\dt$, and it is independent of $n$.  The second inequality in \ref{lem:AP_I} follows, since $R$ is decreasing with respect to $I$.  

\item The first inequality in \ref{lem:AP_minu} is a consequence of the previous result. Indeed, thanks to the definition of $I^{n+1}$ in \eqref{scheme_ep:MonotonicFormulation_I} and to \eqref{assumption:psi}, one has
\[
 \forall i\in\Z,\; \psi_m \;\dx\; \e^{-u^{n+1}_i/\ep} \le \dx\sum\limits_{i\in\Z} \psi(x_i) \e^{-u^{n+1}_i/\ep}=I^{n+1}\le 2I_M.
\]
It gives that $\forall i\in\Z$, $u^{n+1}_i\ge c_m\ep$, where 
$
 c_m=\min\left\{ -\ln\left( \frac{2I_M}{\psi_m \dx} \right), c_m^\tin\right\}$, 
depends only on the constants defined in the assumptions and on $\dx$, and is independent of $n$.

\item 
Eventually, Lemma \ref{lem:Monotonic_ep} yields that $\mathcal{M}_\dt^\ep(u^n)$ enjoys $L_n$-Lipschitz property. 
The $L_{n+1}$-Lipschitz bound \ref{lem:AP_Lipschitz} of $u^{n+1}$ is then a consequence of \ref{lem:AP_I} and \eqref{assumption:R_bounded_decreasing}.
\end{itemize}
Eventually, we denote $\ep_0=\min(\ep_1,\ep_2,\ep_3)$, such that Lemma \ref{lem:AP} holds. 
\end{proof}

This technical lemma provides the necessary tools to prove the convergence of the sequences defined by the scheme \eqref{scheme:epsilon} to the sequences defined by the scheme \eqref{scheme:limit}: 

\begin{proof}[Proof of Prop. \ref{prop:AP}]
 As for the proof of Lemma \ref{lem:AP}, we proceed by induction. Thanks to the assumptions, there exists a sequence $(v^0_i)_{i\in\Z}$ such that $u^0_i\underset{\ep\to 0}\longrightarrow v^0_i$ for all $i\in\Z$. We suppose that it is true for a given $n\in\ccl 0,N_t-1\ccr$ and we prove that there exist $(v^{n+1}_i)_{i\in\Z}$ and $J^{n+1}\in\R$ such that 
 \[
  \forall i\in\Z,\; u^{n+1}_i\underset{\ep\to 0}\longrightarrow v^{n+1}_i,\;\;\mathrm{and}\;\; I^{n+1}\underset{\ep\to 0}\longrightarrow J^{n+1}.
 \]
 First of all, $(u^n_i)_{i\in\Z}$ enjoys the Lipschitz property \ref{lem:AP_Lipschitz} of Lemma \ref{lem:AP}, and  \eqref{eq:CFL_AP} holds. These properties are uniform with respect to $\ep$ small enough, thus the convergence of the first step \eqref{scheme_ep:MonotonicFormulation_utilde} of the scheme \eqref{scheme:epsilon} follows immediately
 \[
  \forall i\in\Z,\; 
  \mathcal{M}^\ep_\dt(u^n)_i
  \underset{\ep\to 0}\longrightarrow 
  \mathcal{M}^0_\dt(v^n)_i
  = v^n_i -\dt H\left( \frac{v^n_i-v^n_{i-1}}{\dx}, \frac{v^{n}_{i+1}-v^n_i}{\dx}\right).
 \]
Moreover, Lemma \ref{lem:AP} gives that  $(I^{n+1})_{\ep\in(0,\ep_0)}$ is uniformly bounded with respect to $\ep$, so that $I^{n+1}\longrightarrow J^{n+1}$ when $\ep\to 0$, up to an extraction. It provides an extraction of $(u^{n+1}_i)_{\ep\in(0,\ep_0)}$ such that 
\[
 \forall i\in\Z, \; u^{n+1}_i \underset{\ep\to 0}{\longrightarrow} v^{n+1}_i =  
 \mathcal{M}^0_\dt(v^n)_i
 -\dt R(x_i, J^{n+1}), 
\]
and such that $\min_{i\in\Z} v^{n+1}_i = 0$, thanks to the point \ref{lem:AP_minu} of Lemma \ref{lem:AP}. Hence, $(v^{n+1}_i)_{i\in\Z}$ satisfies the scheme \eqref{scheme:limit}. 

To conclude the proof, one has to prove that all extractions of $(I^{n+1})_{\ep\in(0,\ep_0)}$ converge to the same limit. We proceed by contradiction, supposing that there are two extractions which converge respectively to $J^{n+1}_a$ and $J^{n+1}_b$, with $J^{n+1}_a< J^{n+1}_b$. As previously, it provides two extractions of $(u^{n+1}_i)_{\ep\in(0,\ep_0)}$ which converges respectively to $v^{n+1}_{i,a}$ and $v^{n+1}_{i,b}$ when $\ep\to 0$, where 
\[
 \left\{ 
\begin{array}{l l}
\ds \forall i\in\Z,\; v^{n+1}_{i,a}= \tilde{v}^{n+1}_i -\dt R\left(x_i,J^{n+1}_a\right) & \ds \;\;\mathrm{and}\;\;\;\;\min\limits_{i\in\Z} v^{n+1}_{i,a}=0, \vspace{4pt} \\
\ds\forall i\in\Z, \; v^{n+1}_{i,b}= \tilde{v}^{n+1}_i -\dt R\left(x_i,J^{n+1}_b\right) &  \ds \;\;\mathrm{and}\;\;\;\;\min\limits_{i\in\Z} v^{n+1}_{i,b}=0,
\end{array}
 \right.
\]
and as $R$ is decreasing with respect to its second variable \eqref{assumption:R_bounded_decreasing}, 
\[
 \forall i\in\Z, \; v^{n+1}_{i,a} -v^{n+1}_{i,b}=\dt \left( R\left(x_i, J^{n+1}_{i,b}\right)-R\left(x_i,J^{n+1}_{i,a}\right)\right) < 0. 
\]
Eventually, we remark that $(v^{n+1}_{i,b})_{i\in\Z}$ is increasing enough at infinity, since the inequality \ref{lem:AP_LowEstimate} of Lemma \ref{lem:AP} is 
uniform with respect to $\ep$
when $\ep\to 0$. As a consequence, 
\[
\exists j\in\Z,\;v^{n+1}_{j,b}=\min_{i\in\Z}v^{n+1}_{i,b}=0.
\] 
But the previous inequality then gives that $v^{n+1}_{j,a} <0$, which  contradicts with the fact that $\min_{i\in\Z} v^{n+1}_{i,a}=0$.
\end{proof}

\section{Convergence of the limiting scheme \eqref{scheme:limit}}
\label{sec:scheme_0}

As in Section \ref{sec:AP}, 
stability estimates for the scheme \eqref{scheme:limit} are obtained using
the convenient reformulation \eqref{scheme_0:MonotonicFormulation} of the scheme \eqref{scheme:limit}, in which the monotonic component of the scheme \eqref{scheme:Monotonic} is taken apart. 
We start by recalling useful properties of the monotonic scheme \eqref{scheme:Monotonic} (see \cite{CrandallLions}):

\begin{lem}
 \label{lem:Monotonic} Let $s\in(0,\dt]$ and $\mathcal{M}_s$ defined as in \eqref{scheme:Monotonic}. Let $L>0$,
 and suppose that $\dt C_H(L)\le \dx$, with $C_H(L)$ defined in \eqref{eq:CH_L}. 
Then the following results hold true 
\begin{enumerate}
\item If $f(x)$, $f(x\pm\dx)$, $g(x)$ and $g(x\pm\dx)$ are such that $f(x)\le g(x)$, $f(x\pm\dx)\le g(x\pm\dx)$, and
\[
 \left|\frac{f(x)-f(x\pm\dx)}{\dx} \right|\le L, \;\;\; \left|\frac{g(x)-g(x\pm\dx)}{\dx}\right| \le L,
\]
then $\mathcal{M}_s(f)(x) \le \mathcal{M}_s(g)(x)$. 
\item In particular, if $f$ and $g$ are two $L$-Lipschitz functions such that,  $f\le g$, then $\mathcal{M}_s(f)\le\mathcal{M}_s(g)$. 
Moreover, both $\mathcal{M}_s(f)$ and $\mathcal{M}_s(g)$ are $L$-Lipschitz continuous. 
\end{enumerate}
\end{lem}

\begin{rmq}
In particular, using the notations and assumptions of Lemma \ref{lem:Monotonic}, 
if $f:\R\mapsto\R$ is a $L$-Lipschitz function such that $\forall i\in\Z$, $f(x_i)\ge 0$, then 
$
 \forall i\in\Z, \;\mathcal{M}_s(f)(x_i)\ge 0.
$
Moreover, 
if
\[
 \exists (\underline{a},\overline{a})\in[0,L]^2, \; \exists (\underline{b},\overline{b})\in\R^2, \; \forall x\in\R, \underline{a}|x-x_0|+\underline{b} \le f(x) \le \overline{a}|x-x_0| + \overline{b}, 
\]
then 
\[
 \forall x\in\R, \; \underline{a}|x-x_0|+\underline{b}-sH(\underline{a},\underline{a}) \le \mathcal{M}_s(f)(x)\le \overline{a}|x-x_0|+\overline{b}.
\]
\end{rmq}

Using these notations, we prove   
the following lemma, which establishes stability properties of the scheme \eqref{scheme:limit}, as well as the fact that $J_\dt$ is nondecreasing.

\begin{lem}
 \label{lem:stability_scheme0}
 Suppose that the assumptions of Prop. \ref{prop:scheme_0} are satisfied, and that $v_\dt$ and $J_\dt$ are defined in  \eqref{scheme_0:MonotonicFormulation}. The following results hold: 
 \begin{enumerate}
  \item  \label{lem:LimitSpace_Lipschitz}
  Uniform Lipschitz continuity in trait: for all $t\in[0,T]$, there exists a constant $L_t=L_0+tK\le L_T$, with $L_0$ defined in \eqref{assumption:u0Lipschitz} and $K$ in \eqref{assumption:R_bounded_decreasing}, such that $v_\dt(t,\cdot)$ 
is $L_t$-Lipschitz continuous.
\item \label{lem:LimitTime_Lipschitz}
Uniform Lipschitz continuity in finite time: for all $x\in\R$, $v_\dt(\cdot,x)$ is $(L_T^2+K)$ Lipschitz continuous on $[0,T]$, 
where $L_T$ is defined in \ref{lem:LimitSpace_Lipschitz} and $K$ in \eqref{assumption:R_bounded_decreasing}.
  \item \label{lem:LimitSpace_vBoundedBelow} 
  Uniform bounds for $v_\dt$: 
  for all $t\in[0,T]$, there exist $\underline{b}_t, \overline{b}_t\in\R$ such that $\underline{b}_t=\underline{b}-tH(\underline{a},\underline{a})-tK$ and $\overline{b}_t=\overline{b}+tK$, such that 
  \[
   \forall x\in\R, \; \underline{a}|x-x_0|+\underline{b}_t \le v_\dt(t,x)\le \overline{a}|x-x_0|+\overline{b}_t.
  \]
where $\underline{a},\; \overline{a},\; \underline{b}$ and $\overline{b}$ are defined in \eqref{assumption:u_increasing_infty}, $H$ in \eqref{eq:H}, and K in \eqref{assumption:R_bounded_decreasing}. 
    \item \label{lem:LimitSpace_J}
    Uniform bounds for  $J_\dt$:
    $
    \forall t\in(0,T]$, $I_m\le J_\dt(t) \le I_M.
    $
\item
\label{lem:LimitTime_Jincreasing} 
Monotonicity of $J_\dt$: $J_\dt$ is nondecreasing on $(0,T]$. 
 \end{enumerate}
\end{lem}

\begin{rmq}
 This lemma is similar to the stability properties stated in Section \ref{sec:AP}, but it is important to notice that all the constants are independent of the discretization. Moreover, the result  \ref{lem:LimitTime_Jincreasing} only holds for the limit scheme \eqref{scheme:limit}.
\end{rmq}

\begin{rmq}
Lemma \ref{lem:stability_scheme0} gives hints for the implementation of
\eqref{scheme:limit} independently of \eqref{scheme:epsilon}. Indeed thanks to the properties above, and to \eqref{assumption:R_extrema}-\eqref{assumption:R_bounded_decreasing}, 
\begin{equation}
\label{eq:def_Jnp1}
 [I_m,I_M] \ni J \mapsto \min\limits_{i\in\Z}\left\{v^n_i - \dt H\left(\frac{v^n_i-v^n_{i-1}}{\dx},\frac{v^n_{i+1}-v^n_i}{\dx}\right) -\dt R(x_i, J)\right\},
\end{equation}
is increasing, takes a negative value at $I_m$, a positive one at $I_M$, and it is equal to $0$ at $J^{n+1}$. One can also notice that it is continuous, as the minimum in \eqref{eq:def_Jnp1} is taken on a finite number of indices, thanks to Lemma \ref{lem:stability_scheme0}-\ref{lem:LimitSpace_vBoundedBelow} and \eqref{assumption:R_bounded_decreasing}.   Even with no further result on the regularity of $R$, and hence on \eqref{eq:def_Jnp1}, $J^{n+1}$ can be approximated, for instance by dichotomy.   In practice, an approximated Newton's method works, and it is more efficient in terms of computational time.
\end{rmq}

\begin{proof}
 We start by proving \ref{lem:LimitSpace_Lipschitz}-\ref{lem:LimitSpace_vBoundedBelow} and \ref{lem:LimitSpace_J} by induction.  
 Let $n\in\ccl 0,N_t-1\ccr$. Suppose that \ref{lem:LimitSpace_Lipschitz}-\ref{lem:LimitSpace_vBoundedBelow} are true for $v_\dt(t_n,\cdot)$, as it is the case for the initial data $v^\tin$ thanks to \eqref{assumption:u_increasing_infty}-\eqref{assumption:u0Lipschitz}. Let $n\in\ccl 0,N_t-1\ccr$ and $s\in(0,\dt]$. In what follows, we show that $v_\dt(t_n+s,\cdot)$ satisfies \ref{lem:LimitSpace_Lipschitz}-\ref{lem:LimitSpace_vBoundedBelow}, and that ${J}_\dt(t_{n+1})$ satisfies \ref{lem:LimitSpace_J}:
 \begin{itemize}
 \item  Let $j$ realize the minimum of $(v_\dt(t_n,x_i))_{i\in\Z}$. Hence $v_\dt(t_n,x_j)=0$ thanks to \eqref{scheme_0:MonotonicFormulation_J}, and   the definition of $H$ in \eqref{eq:H} yields
\[
 \mathcal{M}_\dt\left(v_\dt(t_n,\cdot)  \right)(x_j)=0.
\]
Coming back to \eqref{scheme_0:MonotonicFormulation_v}, we have
\[
 -\dt R\left(x_j, {J}_\dt(t_{n+1})\right) = v_\dt(t_{n+1},x_j) \ge 0,
\]
and 
 we obtain that $I_m\le {J}_\dt(t_{n+1})$ thanks to \eqref{assumption:R_extrema} and \eqref{assumption:R_bounded_decreasing}. 
 \item Notice that since \eqref{eq:CFL_0} is satisfied, 
 the first step \eqref{scheme:Monotonic} of the scheme \eqref{scheme:limit} is monotonic. Hence, Lemma \ref{lem:Monotonic} gives
    \[
 \forall x\in\R, \; \underline{a}|x-x_0|+\underline{b}_{t_n}-sH(\underline{a},\underline{a}) \le \mathcal{M}_s\left((v_\dt(t_n,\cdot)\right)(x),
\]
so that $v_\dt(t_n+s,x) \ge \underline{a}|x-x_0| + \underline{b}_{t_n} - sH(\underline{a},\underline{a}) - sR(x,J_\dt(t_n+s))$. The monotonicity of $R(x,\cdot)$ as well as the lower bound for $J_\dt(t_n+s)$ yield
\[
 \underline{a}|x-x_0| +\underline{b}_{t_n+s} \le v_\dt(t_n+s,x),
\]
with $\underline{b}_{t_n+s} =\underline{b}_{t_n}-sH(\underline{a},\underline{a}) - sK$, thanks to \eqref{assumption:R_bounded_decreasing}. 
\item Since $v_\dt(t_n,x_i)\ge 0$ for all $i\in\Z$, Lemma \ref{lem:Monotonic} yields that  
 \[
 \forall i\in\Z,\; \mathcal{M}_\dt\left( v_\dt(t_n,\cdot)  \right)(x_i)\ge 0. 
 \]
Consider then $k\in\Z$ such that 
\[
 v_\dt(t_{n+1}, x_k)=\min\limits_{i\in\Z}v_\dt(t_{n+1},x_i) = 0. 
\]
Note that such a $k$ exists, thanks to the previous step of the proof. 
We obtain 
\[
 \dt R\left(x_k, {J}_\dt(t_{n+1})\right) = \mathcal{M}_\dt\left( v_\dt(t_n,\cdot)  \right)(x_k)\ge 0, 
\]
thanks to \eqref{scheme_0:MonotonicFormulation_v}.  The inequality ${J}_\dt(t_{n+1})\le I_M$ is then a consequence of  assumptions \eqref{assumption:R_extrema} and \eqref{assumption:R_bounded_decreasing}. The bounds for $J_\dt$ in \ref{lem:LimitSpace_J} follow, since it is constant on $(t_n,t_{n+1}]$. 
  \item Once again, the monotonicity of the first step \eqref{scheme:Monotonic} of scheme \eqref{scheme:limit}, yields
    \[
 \forall x\in\R,  \mathcal{M}_s\left(v_\dt(t_n,\cdot)\right)(x)\le \overline{a}|x-x_0|+\overline{b}_{t_n}, 
\]
so that
property \ref{lem:LimitSpace_vBoundedBelow} is proved with $\overline{b}_{t_n+s} = \overline{b}_{t_n}+sK,$ thanks to \eqref{assumption:R_bounded_decreasing}. 
\item Similarly, $v_\dt(t_n+s,\cdot)$ is $L_{t_n+s}$-Lipschitz continuous thanks to Lemma \ref{lem:Monotonic} and \eqref{assumption:R_bounded_decreasing}. 
 \end{itemize}
 
The Lipschitz-in-time property \ref{lem:LimitTime_Lipschitz} is a consequence of \ref{lem:LimitSpace_Lipschitz}. 
We now show that ${J}_\dt(t_{n+1})\ge {J}_\dt(t_n)$. 
Recalling  that $J_\dt$ is constant on $(0,\dt]$ and that is not defined at $t=0$, we then  suppose  that $n\in\ccl 1,N_t-1\ccr$. 
Considering an index $j$ such that $v_\dt(t_n,x_j)=\min_{i\in\Z} v_\dt(t_n,x_j)=0$, \eqref{scheme_0:MonotonicFormulation_v} yield 
 \[
  R(x_j, {J}_\dt(t_{n+1})) \le 0,
 \]
as previously. 
Let us now consider the previous step of the scheme, at the same index $j$. As this part of the proof only uses points of the grid, we use rather the formulation \eqref{scheme:limit} for the sake of simplicity. We have, 
\begin{equation}
\label{eq:lem_LimitTime_J_step1}
 \left\{ 
 \begin{array}{l}
 \ds \frac{v^n_j-v^{n-1}_j}{\dt} + H\left( \frac{v^{n-1}_j-v^{n-1}_{j-1}}{\dx}, \frac{v^{n-1}_{j+1}-v^{n-1}_j}{\dx}\right)+ R(x_j,J^n)=0 \vspace{4pt} \\
 \ds v^n_j=\min\limits_{i\in\Z} v^n_i=0.
 \end{array}
 \right.
\end{equation}
Since all the $v^{n-1}_i$ for $i\in\Z$ are nonnegative, one has 
\begin{equation*}
 \frac{v^{n-1}_j-v^{n-1}_{j-1}}{\dx} \le \frac{v^{n-1}_j}{\dx}, \;\;\;\;\mathrm{and}\;\;\;\; 
 \frac{v^{n-1}_{j+1}-v^{n-1}_j}{\dx}\ge \frac{-v^{n-1}_j}{\dx}.
\end{equation*}
Moreover, because $H$ is increasing with respect to its first variable and decreasing with respect to the second one, the following inequality holds
 \[
  H\left( \frac{v^{n-1}_j-v^{n-1}_{j-1}}{\dx}, \frac{v^{n-1}_{j+1}-v^{n-1}_j}{\dx}\right) \le H\left( \frac{v^{n-1}_j}{\dx}, \frac{-v^{n-1}_j}{\dx}\right) = \left( \frac{v^{n-1}_j}{\dx} \right)^2,
 \]
where the last equality comes from the expression of $H$, see \eqref{eq:H}. Once injected in \eqref{eq:lem_LimitTime_J_step1}, we obtain 
\[
 R(x_j,J^n) \ge 
 \frac{v^{n-1}_j}{\dt}\left( 1-\frac{\dt}{\dx^2}\;v^{n-1}_j \right),
\]
and the right hand side of this inequality is positive. Indeed, thanks to the Lipschitz-in-time property \ref{lem:LimitTime_Lipschitz}, we have
\[
 \left|\frac{v^n_j-v^{n-1}_j}{\dt}\right| = \frac{v^{n-1}_j}{\dt}\le L_T^2+K,
\]
and the condition \eqref{eq:CFL_0} yields the result.
To conclude, let us remark that 
\begin{equation}
\label{eq:argmin}
 R\left(x_j, {J}_\dt(t_{n+1})\right)\le 0 \le R\left(x_j,J^n\right)=R\left(x_j,{J}_\dt(t_n)\right),
\end{equation}
and use the fact that $R$ is decreasing with respect to its second variable. The monotonicity of $J_\dt$ in \ref{lem:LimitTime_Jincreasing} follows immediately since it is constant on the interval $(t_n,t_{n+1}]$.
\end{proof}

\begin{rmq}
  Let us emphasize the fact that the above proof strongly relies on considerations on the minimum of $(v^n_j)_j$. This bears similarities with \cite{BarlesPerthame}, where the  relation  $R(\overline{x}(t), J(t))=0$, with
 $\overline{x}(t)= \arg\min v(t,\cdot)$, is used to study $J$. In the discrete setting, \eqref{eq:argmin} is the equivalent of this relation. 
\end{rmq}

The next step consists in establishing the convergence of $v_\dt$ and $J_\dt$ defined in  \eqref{scheme_0:MonotonicFormulation}  when $\dt$ and $\dx$ go to $0$ with $\dt/\dx$ fixed. The following results hold

\begin{lem}
 \label{lem:PointwiseCVofv_dtJ_dt} Suppose that the assumptions of Prop. \ref{prop:scheme_0} are satisfied, and that $v_\dt$, $J_\dt$ are defined by \eqref{scheme_0:MonotonicFormulation}. Then, 
 \begin{enumerate}
 \item \label{lem:PointwiseCV_vdt} Convergence of $(v_\dt)_{\dt >0}$: there exists $ v_0\in\mathcal{C}^0([0,T]\times\R), $ such that
 \[
  \forall (t,x)\in[0,T]\times \R, \; v_\dt(t,x) \underset{\dt\to 0}\longrightarrow v_0(t,x) \text{\;up to a subsequence},
 \]
and with $\min v_0(t, \cdot)=0$, for all $t\in[0,T]$. Moreover, the convergence is locally uniform on $[0,T]\times\R$.
 \item \label{lem:PointwiseSV_Jdt} Convergence of $(J_\dt)_{\dt >0}$: 
 there exists $J_0\in BV(0,T)$, lower semi-continuous, 
 such that 
 \[
  \text{for\;almost\;all\;} t\in (0,T], \; J_\dt(t) \underset{\dt\to 0}\longrightarrow J_0(t) \text{\;up to a subsequence}.
 \]
 Moreover, $J_0$ is nondecreasing, and $\forall t\in(0,T]$,  $I_m\le J_0(t)\le I_M$. 
 \end{enumerate}
\end{lem}

\begin{proof}
 Thanks to
 Lemma \ref{lem:stability_scheme0}, the family
 $(v_\dt)_{\dt>0}$ is composed of Lipschitz functions, having the same Lipschitz constant. 
 Considering $R>0$, one can notice that since $v_\dt$ enjoys Lipschitz-in-time regularity
\[
 \|v_\dt\|_{\LL^\infty([0,T]\times[-R,R])}\le 
 (L_T^2+K)T+ \|v^\tin\|_{\LL^\infty([-R,R])},
\]
hence the family $(v_\dt)$ satisfies the hypothesis of Ascoli's theorem for $(t,x)\in[0,T]\times[-R,R]$. 
Then, there exists a function $v_0\in\mathcal{C}^0([0,T]\times\R)$ such that $v_\dt\longrightarrow_{\dt\to0} v_0$ uniformly on $[0,T]\times[-R,R]$. 
Moreover, $(v_\dt)_\dt$ is a sequence of uniformy coercive and Lipschitz functions, such that 
\[
\forall n\in\ccl 0,N_t\ccr, \;\min\limits_{i\in\Z}v_\dt(t_n,x_i)=0.
\]
Hence, there exists a constant $c$ such that for all $t\in[0,T]$, 
\[
 \left|\min v_\dt (t,\cdot)\right|\le c(\dt+\dx),
\]
and $\min v_0(t,\cdot)=0$ is a consequence of the local uniform convergence of  $(v_\dt)_\dt$  to  $v_0$. 
This proves \ref{lem:PointwiseCV_vdt}.

The second point \ref{lem:PointwiseSV_Jdt} is a consequence of Helly's selection theorem. Indeed,  Lemma \ref{lem:stability_scheme0}-\ref{lem:LimitSpace_J}-\ref{lem:LimitTime_Jincreasing} states that $(J_\dt)_{\dt>0}$ is a sequence of uniformly bounded BV functions with uniformly bounded total variation. Hence, there exists a BV function $\tilde{J}_0$ such that 
 \[
  J_\dt\underset{\dt\to 0}\longrightarrow\tilde{J}_0, \text{\;pointwise in\;} (0,T]\text{\;up to a subsequence}.
 \]
Moreover, $\tilde{J}_0$ is nondecreasing, and $I_m\le \tilde{J}_0\le I_M$, since these properties hold for all $J_\dt$. Considering a lower semi-continuous function $J_0$ such that $J_0=\tilde{J}_0$ almost everywhere in $(0,T]$ yields the result. 
\end{proof}

\begin{rmq}
\label{rmq:subsequence}
  In what follows, the mention $\dt\to 0$ will always refer to a subsequence for which the convergences of Lemma \ref{lem:PointwiseCVofv_dtJ_dt} hold true. 
\end{rmq}

\begin{rmq}
 Note that, although it is not defined by the scheme, a value for $J_\dt(0)$ is needed in what follows, because of the compactness argument used below. When it is necessary, we define $J_\dt(0)=J^1$. This choice consists in extending continuously $J_\dt$ at $0$, but it has no meaning from the point of view of the constraint of the scheme. However, it 
 is well-suited to the fact that $J_\dt$ is bounded and nondecreasing.
\end{rmq}

To complete the proof of Prop. \ref{prop:scheme_0}, it remains to identify
$v_0=v$ and $J_0=J$ almost everywhere, where $(v,J)$ is the viscosity solution of \eqref{eq:limit}. However, $J_0$ enjoys only BV regularity, and in particular it is not expected to be continuous (we refer to Section \ref{sec:tests}, where numerical tests show that $J_0$ can have jumps).   As a consequence, general convergence results of numerical schemes for Hamilton-Jacobi equations such as \cite{CrandallLions} cannot be applied directly. To the best of our knowledge, there is no general framework for finite-differences numerical schemes for Hamilton-Jacobi equation when the Hamiltonian is not continuous in time. In what follows, we propose a proof of the convergence of the scheme \eqref{scheme:limit} to the viscosity solution of \eqref{eq:limit}. The key ingredient of the proof is an appropriate regularization of $J_\dt$ and $J_0$, used in \cite{AmbrosioFuscoPallara2000}, and also in  \cite{CalvezLam2020} for the study of the uniqueness of viscosity solution of constrained Hamilton-Jacobi equation. For $k>0$ and $\dt\ge 0$, let us define
\begin{align}
 \label{eq:J_delta_k} 
 \forall t\in[0,T], \;&J_\dt^k(t)=\inf\limits_{s\in[0,T]}(J_\dt(s) + k|t-s| ).
\end{align}
The following results hold true

\begin{lem}
 \label{lem:UniformCVJdeltak} 
 Suppose that the assumptions of Prop. \ref{prop:scheme_0} are satisfied. Let  $J_\dt^k$ and $J_0^k$ defined by \eqref{eq:J_delta_k}. Then,
 \begin{enumerate}
 \item For all $\dt\ge 0$, and for all $k>0$, $I_m\le J_\dt^k\le I_M$, and $J_\dt^k$ is a nondecreasing function on $[0,T]$.
  \item For fixed $\dt\ge 0$, and for all $t\in[0,T]$, $J_\dt^k(t) \nearrow J_\dt(t)$ when $k\to +\infty$. 
\item For fixed $\dt\ge 0$,  $J_\dt^k$ is a $k$-Lipschitz function on $[0,T]$.
\item For fixed $k>0$,  $\|J^k_\dt - J^k_0\|_\infty \underset{\dt\to 0}\longrightarrow 0$. \label{lem:UniformCVJdeltak_iii}
 \end{enumerate}
\end{lem}

\begin{proof}
 We only detail the proof of \ref{lem:UniformCVJdeltak_iii}.  
 Let $k>0$ be fixed.
 From Lemma \ref{lem:PointwiseCVofv_dtJ_dt}, 
 $J_\dt(t)\to J_0(t)$ almost everywhere in $[0,T]$ when $\dt\to 0$. We first remark that 
 \begin{equation}
 \label{eq:lem_JdtkCV_CVae}
  J_\dt^k(t) \underset{\dt\to 0}\longrightarrow J_0^k(t) \;\text{a.\;e.\;in\;} [0,T].
 \end{equation}
Indeed, let us consider $t\in[0,T]$ such that $J_\dt(t)\to_{\dt\to 0} J_0(t)$. Since $(J_\dt^k(t))_{\dt>0}$ is a bounded sequence, it admits a converging subsequence, once again denoted by $(J_\dt^k(t))_{\dt>0}$. Let us denote by $\ell$  its limit. Since
\[
 \forall s\in[0,T], \; J^k_\dt (t) \le J_\dt(s)+k|t-s|, 
\]
then letting $\dt\to 0$ in the previous inequality yields 
\[
 \forall s\in [0,T], \; \ell \le J_0(s) + k|t-s|, 
\]
so that $\ell\le J_0^k(t)$. Moreover, as $J_\dt^k(t)$ is defined as an infimum,
\[
 \forall n\in\N^*, \; \exists s^*_n\in[0,T], \;J_\dt(s_n^*) +k|t-s_n^*|-\frac{1}{n} \le J^k_\dt(t).
\]
Since $(s^*_n)_{n\ge 1}$ converges (up to an extraction) to $s^*\in[0,T]$ when $n\to +\infty$, taking the $\liminf$ in the previous inequality gives
\[
 J_0(s^*)+k|t-s^*|\le \ell,
\]
because $J_0$ is lower semi-continuous. As a consequence $J_0^k(t)\le \ell$. The only adherence value of $(J_\dt^k(t))_{\dt>0}$ is then $J_0^k(t)$, which yields  \eqref{eq:lem_JdtkCV_CVae}. To conclude, the uniform convergence in \ref{lem:UniformCVJdeltak_iii} is a consequence of the convergence almost everywhere of a family of Lipschitz functions defined on a compact domain. 
\end{proof}

\begin{rmq} The uniform convergence in Lemma
 \ref{lem:UniformCVJdeltak}-\ref{lem:UniformCVJdeltak_iii} does not generally hold true in the limit $k\to \infty$. This result will only be used for fixed $k>0$. 
\end{rmq}

Now that $J_\dt^k$ and $J_0^k$ are defined, we consider them as a source term respectively in the scheme and in the equation. Namely, let us define $v^k$ the viscosity solution of the Hamilton-Jacobi equation 
\begin{equation}
 \label{eq:eqlimit_regularized} 
 \partial_t v^k + |\nabla_x v^k |^2 = - R(x,J_0^k), \; x\in\R, \; t>0,
\end{equation}
with initial data $v^\tin$. 
Thanks to the Lipschitz properties of the right-hand side of \eqref{eq:eqlimit_regularized}, $v^k$ exists, is uniquely determined, and enjoys Lipschitz-regularity properties. Moreover, 
the following lemma establishes that,
because of the construction of $J^k_0$, $v^k$ converges when $k\to +\infty$ to the viscosity solution of \eqref{eq:eqlimit_regularized} with $J_0$ instead of $J^k_0$. 
Similarly, let us define
$v_\dt^k$ by 
\begin{equation}
 \label{eq:schemelimit_regularized} 
 v_\dt^k(t_n+s,x) = \mathcal{M}_s(v_\dt^k(t_n,\cdot))(x)-s\;R(x, J_\dt^k(t_n+s)),
\end{equation}
for all $n\in\ccl 0,N_t-1\ccr$, $s\in(0,\dt]$,
and $x\in\R$,
with $\mathcal{M}_s$ defined in \eqref{scheme:Monotonic},
and initialized with $v_\dt^k(0,\cdot)=v^\tin$. The properties of $v^k$ and $v^k_\dt$ are summarized in the following lemma: 
\begin{lem}
\label{lem:vdtk}
 Suppose that the assumptions of Prop. \ref{prop:scheme_0} are satisfied. Let $k>0$, $v^k$ and  $v^k_\dt$ defined by \eqref{eq:eqlimit_regularized} and \eqref{eq:schemelimit_regularized}. Then, $v^k$ and $v_\dt^k$ enjoy the following properties 
 \begin{enumerate}
  \item \label{lem:vdtk_space} Uniform Lipschitz continuity in trait: 
  for all $t\in[0,T]$, $v^k(t,\cdot)$ and $v^k_\dt(t,\cdot)$ are $L_T$-Lipschitz continuous, with $L_T$ defined in Lemma \ref{lem:stability_scheme0}-\ref{lem:LimitSpace_Lipschitz}.
\item \label{lem:vdtk_time} Uniform Lipschitz continuity in finite time: 
for all $x\in\R$, $v^k(\cdot,x)$ and $v^k_\dt(\cdot,x)$ are 
$(L_T^2+K)$-Lipschitz 
continuous on $[0,T]$, where 
$L_T$ is defined in Lemma \ref{lem:stability_scheme0}-\ref{lem:LimitSpace_Lipschitz},
and $K$ in \eqref{assumption:R_bounded_decreasing}.
\item \label{lem:vk_bounds} Uniform bounds for $v^k$: for all $(t,x)\in[0,T]\times\R$, 
 \[
  \underline{b}-KT \le v^k(t,x) \le \overline{a}|x-x_0| + \overline{b}+TK,
 \]
 where $K$, $\underline{b}$, $\overline{b}$, and $\overline{a}$ are defined in \eqref{assumption:R_bounded_decreasing} and \eqref{assumption:u_increasing_infty}.
\item \label{lem:vdtk_below} Uniform bounds for $v_\dt^k$: for all $(t,x)\in[0,T]$,
\begin{equation}
 \label{eq:vdtk_below} 
 \underline{a}|x-x_0| +\underline{b}_t \le v_\dt^k(t,x) \le \overline{a}|x-x_0|+\overline{b}_t,
\end{equation}
where $\underline{a}$, $\underline{b}_t$, $\overline{a}$ and $\overline{b}_t$ are defined in \eqref{assumption:u_increasing_infty} and in Lemma \ref{lem:stability_scheme0}-\ref{lem:LimitSpace_vBoundedBelow}. 
 \item \label{lem:vk_Monotonicity} Monotonicity of the approximation: $v^k\nearrow v^\infty$  when $k\to +\infty$, pointwise in $[0,T]\times\R$ , where $v^\infty$ is the viscosity solution of 
\begin{equation}
 \label{eq:limit_lsc}
 \partial_t v^\infty + |\nabla_x v^\infty|^2 = -R(x, J_0), \; x\in\R, t\in(0,T],
\end{equation}
initialized with $v^\tin$.
\item \label{lem:vdtk_Monotonicity} Monotonicity of the approximation: 
$v_\dt^k\le v_\dt$,
where $v_\dt$ is defined in \eqref{scheme_0:MonotonicFormulation}.
 \end{enumerate}
\end{lem}

\begin{proof} Concerning the properties of $v^k$,
the points \ref{lem:vdtk_space} and \ref{lem:vdtk_time} are 
natural properties of viscosity solution of \eqref{eq:eqlimit_regularized}, while \ref{lem:vk_bounds} is a consequence of the comparison principle.  Point \ref{lem:vk_Monotonicity} is proved in \cite{CalvezLam2020}. 

 Concerning $v^k_\dt$, since we suppose that \eqref{eq:CFL_0} is satisfied, the proofs of the first points of  Lemma \ref{lem:stability_scheme0} can be applied. This yields immediately \ref{lem:vdtk_space}-\ref{lem:vdtk_time} and \ref{lem:vdtk_below}. 
 The last point of the Lemma is a consequence of the monotonicity of the scheme (Lemma \ref{lem:Monotonic}), and is done by induction. Indeed, the inequality \ref{lem:vdtk_Monotonicity} holds true at $t=0$.
 Moreover,
$v_\dt$ and $v_\dt^k$ 
 enjoy the Lipschitz properties 
 of Lemmas \ref{lem:stability_scheme0}-\ref{lem:LimitSpace_Lipschitz}-\ref{lem:LimitTime_Lipschitz} and \ref{lem:vdtk}-\ref{lem:vdtk_space}-\ref{lem:vdtk_time},
 and \eqref{eq:CFL_0} is satisfied. As a consequence, the first step \eqref{scheme:Monotonic} of the reformulation of the scheme \eqref{scheme:limit} is monotonic. 
Hence,
if $v_\dt^k(t_n,\cdot)\le v_\dt(t_n,\cdot)$, one has $\mathcal{M}_s(v_\dt^k(t_n,\cdot))\le \mathcal{M}_s(v_\dt(t_n,\cdot))$ for all $s\in(0,\dt)$. 
 Eventually, we use Lemma \ref{lem:UniformCVJdeltak} and the fact that $R$ is noincreasing in its second variable to conclude that 
 $v_\dt^k(t_n+s,\cdot)\le v_\dt(t_n+s,\cdot)$ for all $s\in(0,\dt]$.
\end{proof}

\begin{rmq}
 Note that, contrary to the non-regular problems \eqref{eq:limit} and \eqref{scheme:limit},  there is no constraint neither on $\min v^k(t,\cdot)$, nor on $\min_{i\in\Z} 
 v_\dt^k(t_{n+1},x_i)$.  
\end{rmq}

Now that the problem is regularized, we can use viscosity procedures to show that scheme \eqref{scheme:limit} converges to the viscosity solution of \eqref{eq:limit}. Following the ideas developped in \cite{CrandallLions}, let us define an auxiliary function
\begin{align}
\label{eq:psi} 
 \psi(t,x,\tau,\xi)&= v^k(t,x) - v_\dt^k(\tau,\xi) -\frac{(x-\xi)^2}{2\dx^{1/2}} - \frac{(t-\tau)^2}{2\dt^{1/2}} -\left(\sigma +4\mathcal{C}_H^2\alpha\e^T\right) t 
\\ & - \alpha \frac{\e^{ t}}{2}(x^2+\xi^2) - \frac{\alpha}{T-t}  
 \nonumber
\end{align}
for all $(t,x,\tau,\xi)\in[0,T[\times\R\times[0,T]\times\R$. Here, $\alpha\in(0,1)$, and $\sigma$ is positive and will be determined later. The functions $v^k$ and $v_\dt^k$ are defined in \eqref{eq:eqlimit_regularized} and \eqref{eq:schemelimit_regularized},  Then, $\psi$ satisfies the following properties.

\begin{lem}
 \label{lem:psi} 
 Suppose that the assumptions of Prop. \ref{prop:scheme_0} hold, and that $\psi$ is defined by \eqref{eq:psi}. Then 
 \begin{itemize}
\item For all $\alpha$, and $\sigma$ positive, $\psi$ admits a global maximum. It is reached at $(t^*,x^*,\tau^*,\xi^*)\in[0,T[\times\R\times[0,T]\times\R$. 
\item There exists  $\sigma=\sigma(\dt,k)$ with $\sigma(\dt,k)\to_{\dt\to 0}0$ (when $k>0$ is fixed), and $\dt_0>0$, such that for all $\alpha\in(0,1)$ and for all $\dt<\dt_0$, 
$t^*\le 2(L_T^2+K)\dt^{1/2}$, 
where $K$ and $L_T$ are defined in  \eqref{assumption:R_bounded_decreasing} and Lemma \ref{lem:stability_scheme0}-\ref{lem:LimitSpace_Lipschitz}.
 \end{itemize}
\end{lem}

\begin{proof} The first point of the Lemma is immediate, thanks to Lemma \ref{lem:vdtk}-\ref{lem:vdtk_below}-\ref{lem:vk_bounds}. 
The idea of the proof of the second point is very similar to what is done in \cite{CrandallLions}, where monotonic schemes for bounded solutions of Hamilton-Jacobi equations are studied. However, it is worth noticing that, in our framework the boundedness hypothesis is lacking, since it would contradict with the definition of $I_\ep$ in \eqref{eq:u_epsilon}. Moreover, the proof we propose below spies the influence of the regularizations $J^k_0$ and $J^k_\dt$ of $J_0$ and $J_\dt$ through the parameter $k$. 
Indeed, it is necessary  to come back to the non-regularized problem.

\medskip\noindent\textbf{Step} \textbf{(i)}. 
Since $\psi(t^*,\xstar,\tau^*,\xistar) \ge \psi(t^*,0,\tau^*,0)$, we have
\[
 \alp\frac{\e^{ t^*}}{2}(\xstar^2+\xistar^2) \le v^k_\dt(\tau^*,0)-v^k_\dt(\tau^*,\xistar) +v^k(t^*,\xstar)-v^k(t^*,0).
\]
Then,  Lemma \ref{lem:vdtk}-\ref{lem:vdtk_space} gives
\[
 \alp \frac{\e^{ t^*}}{2}\max \left \{|\xstar|,|\xistar| \right \}^2 \le L_T  |\xistar|+L_T|\xstar|, 
\]
which yields 
 \begin{equation}
 \label{eq:proof_psi_compact}
  \alp\e^{ t^*}\max \left \{|\xistar|,|\xstar| \right \}\le 4L_T.
 \end{equation}

\medskip\noindent\textbf{Step} \textbf{(ii)}. 
We proceed as in the previous step. Comparing the values of $\psi$ at $(t^*,\xstar,\tau^*,\xistar)$ and $(t^*,\xstar,\tau^*,\xstar)$, we obtain 
\[
 \frac{(\xstar-\xistar)^2}{2\dx^{1/2}}\le v^k_\dt(\tau^*,\xstar) - v^k_\dt(\tau^*,\xistar) + \alpha\frac{\e^{t^*}}{2}(\xstar^2-\xistar^2),
\]
and Lemma \ref{lem:vdtk}-\ref{lem:vdtk_space} and \eqref{eq:proof_psi_compact} give 
 \begin{equation}
 \label{eq:psiproof_Lipschitz}
   \; |t^*-\tau^*|\le 2(L_T^2+K)\dt^{1/2}, \; \; |x^*-\xi^*|\le 10L_T\dx^{1/2}.
 \end{equation}
Note that the bound for $|\tau^*-t^*|$ is obtained similarly, starting from $\psi(t^*,\xstar,\tau^*,\xistar)\ge \psi(t^*,\xstar,t^*,\xistar)$ and using Lemma \ref{lem:vdtk}-\ref{lem:vdtk_time}.

\medskip \noindent\textbf{Step} \textbf{(iii)}. 
We aim to show that $t^*\le 2(L_T^2+K)\dt^{1/2}$, provided that $\sigma$ is  appropriately chosen. We argue by contradiction, and 
suppose that $t^*> 2(L_T^2+K)\dt^{1/2}$.  It implies that $\tau^*>0$, thanks to \eqref{eq:psiproof_Lipschitz}. Let us start by considering 
\[
 (t,x) \mapsto \psi(t,x,\tau^*,\xi^*) =v^k(t,x) - \phi(t,x),
\]
on $[0,T[\times \R$, with 
\[
 \phi(t,x)= v^k_\dt(\tau^*,\xi^*)+ \left(\sigma +4\mathcal{C}_H^2\alpha\e^T\right) t + \frac{(x-\xi^*)^2}{2\dx^{1/2}} + \frac{(t-\tau^*)^2}{2\dt^{1/2}}  +\alpha\frac{\e^{ t}}{2}(x^2+\xistar^2) + \frac{\alpha}{T-t}.
\]
It admits a maximum, precisely at $(t^*,x^*)$, with $t^*\in(0,T)$. Since $v^k$ is the viscosity solution of \eqref{eq:eqlimit_regularized}, we deduce
\[
 \partial_t \phi(t^*,x^*) + \mathcal{H}\left( \nabla_x \phi(t^*,x^*)\right) + R(x^*, J^k_0(t^*)) \le 0,
\]
that is
\begin{align}
\label{eq:proof_psi_InequationGauche}
\sigma +4\mathcal{C}_H^2\alpha\e^T & + \frac{t^*-\tau^*}{\dt^{1/2}}  +\alpha\frac{\e^{ t^*}}{2}(\xstar^2+\xistar^2)+\frac{\alpha}{(T-t^*)^2} 
\\&+ \mathcal{H}\left( \frac{x^*-\xi^*}{\dx^{1/2}} + \alpha \e^{ t^*}x^*\right) + R(x^*,J^k_0(t^*))\le 0, \nonumber
\end{align}
where $\mathcal{H}$ is defined in \eqref{eq:Hamiltonian_continuous}.
Next, let us consider 
\[
 (\tau,\xi)\mapsto \psi(t^*,x^*,\tau,\xi),
\]
on $[0,T]\times\R$.  As previously, it admits a maximum, precisely at $(\tau^*,\xi^*)$, so that for all $(\tau,\xi)\in[0,T]\times\R$
\begin{equation}
\label{eq:proof_psi_ineq_everywhere}
 v^k_\dt(\tau,\xi) \ge w(\tau,\xi) + k^*,
\end{equation}
with 
\begin{align*}
 w(\tau,\xi) &= -\frac{(x^*-\xi)^2}{2\dx^{1/2}}-\frac{(t^*-\tau)^2}{2\dt^{1/2}} -\alpha\frac{\e^{ t^*}}{2}\xi^2,
 \\
 k^* &= v_\dt^k(\tau^*,\xi^*) + \frac{(x^*-\xi^*)^2}{2\dx^{1/2}} +\frac{(t^*-\tau^*)^2}{2\dt^{1/2}} +
 \alpha \frac{\e^{ t^*}}{2}\xistar^2.
 \end{align*}
Remark that $\tau^*=t_{n^*}+s^*$ with $n^*\in\ccl 0,N_t-1\ccr$ and $s^*\in(0,\dt]$. The previous inequality yields 
\begin{equation}
\label{eq:proof_psi_ineq}
 v^k_\dt(t_{n^*},\xi^*) \ge w(t_{n^*},\xi^*) + k^*.
\end{equation}
The next step consists in applying the scheme \eqref{eq:schemelimit_regularized} to this inequality. To do so, one has to make sure that 
\begin{equation}
\label{eq:proof_psi_Lipw}
  |w(t_{n^*},\xi^*\pm\dx)-w(t_{n^*}, \xi^*)|\le (14L_T+1)\dx,
\end{equation}
 so that \eqref{eq:CFL_0} ensures that scheme \eqref{scheme:Monotonic} enjoys monotonicity.  From the expression of $w(\tau,\xi)$, we have
\[
 \left| \frac{w(t_{n^*},\xistar) - w(t_{n^*},\xistar\pm\dx)}{\dx} \right| \le \alpha \e^{ t^*} |\xistar| + \frac{\dx^{1/2}}{2} + \frac{|\xstar-\xistar|}{\dx^{1/2}} +\alpha \frac{\e^{ T}}{2}\dx.
\]
Hence, if  $\dx$ is chosen small enough, \eqref{eq:proof_psi_Lipw} holds, thanks to \eqref{eq:proof_psi_compact} and \eqref{eq:psiproof_Lipschitz}.  Since the ratio $\dt/\dx$ is fixed, this condition on $\dx$ implies that the result holds for all $\dt\le\dt_0$, for some $\dt_0>0$.  Since \eqref{eq:CFL_0} is satisfied, the first step of the scheme \eqref{scheme:Monotonic} is monotonic and can hence be applied to the inequality \eqref{eq:proof_psi_ineq_everywhere}, using \eqref{eq:proof_psi_ineq}. As $\mathcal{M}_{s^*}$ commutes with constants, it gives
\[
 \mathcal{M}_{s^*}(v_\dt^k(t_{n^*},\cdot))(\xi^*) - s^*R(\xi^*, J_\dt^k(t_{n^*}+s^*)) \ge \mathcal{M}_{s^*}(w(t_{n^*},\cdot))(\xi^*) +k^* - s^* R(\xi^*,J_\dt^k(t_{n^*}+s^*)),
\]
that is 
\begin{align*}
 v_\dt^k(\tau^*, \xi^*) &\ge w(t_{n^*},\xi^*) -s^* H\left(\frac{w(t_{n^*},\xi^*) - w(t_{n^*},\xi^*-\dx)}{\dx},\frac{w(t_{n^*},\xi^*+\dx)-w(t_{n^*},\xi^*)}{\dx}\right)
 \\&+k^* - s^* R(\xi^*, J_\dt^k(\tau^*)).
\end{align*}
The latter yields counterpart of \eqref{eq:proof_psi_InequationGauche}
\begin{align}
 \label{eq:proof_psi_InequationDroite}
 0 &\le  H\left( \frac{x^*-\xistar}{\dx^{1/2}} - \alpha\e^{ t^*} \xistar + \frac{\dx^{1/2}}{2}+\alpha\frac{\e^{ t^*}}{2} \dx, \frac{x^*-\xistar}{\dx^{1/2}}-\alpha\e^{ t^*}\xistar-\frac{\dx^{1/2}}{2} -\alpha\frac{\e^{ t^*}}{2}\dx \right) 
 \\ & + R(\xi^*, J_\dt^k(\tau^*)) +\frac{t^*-\tau^* + s^*/2}{\dt^{1/2}}. \nonumber
\end{align}
Inequalities \eqref{eq:proof_psi_InequationGauche} and \eqref{eq:proof_psi_InequationDroite} are now gathered, so that 
\begin{align}
 \label{eq:proof_psi_lefthandside}
 &\sigma +4\mathcal{C}_H^2\alpha\e^T   +\alpha\frac{\e^{ t^*}}{2}(\xstar^2+\xistar^2) + H\left( \frac{x^*-\xi^*}
{\dx^{1/2}} + \alpha \e^{ t^*}x^*, \frac{x^*-\xi^*}
{\dx^{1/2}} + \alpha \e^{ t^*}x^*\right) 
\\ &- H\left( \frac{x^*-\xistar}{\dx^{1/2}} - \alpha\e^{ t^*} \xistar + \frac{\dx^{1/2}}{2}+\alpha\frac{\e^{ t^*}}{2} \dx, \frac{x^*-\xistar}{\dx^{1/2}}-\alpha\e^{ t^*}\xistar-\frac{\dx^{1/2}}{2} -\alpha\frac{\e^{ t^*}}{2}\dx \right) \nonumber
 \\ \le &  \; 
R(\xi^*, J_\dt^k(\tau^*))-R(x^*,J^k_0(t^*))  + \frac{\dt^{1/2}}{2}, \nonumber
\end{align}
since $\mathcal{H}$ defined in \eqref{eq:Hamiltonian_continuous} and the numerical Hamiltonian satisfy $H(p,p)=\mathcal{H}(p)$ for any $p\in\R$. 
An upper bound for the right hand side is obtained from \eqref{assumption:R_bounded_decreasing}, and from the $k$-Lipschitz regularity of $J_0^k$ in Lemma \ref{lem:UniformCVJdeltak}
\[
 R(\xi^*, J_\dt^k(\tau^*))-R(x^*,J^k_0(t^*))   \le K|\xistar-\xstar| + K \|J_\dt^k-J_0^k\|_\infty + Kk|t^*-\tau^*|,
\]
that can, once again, be estimated using \eqref{eq:psiproof_Lipschitz} and the fact that the ratio $\dt/\dx$ is fixed.
On the other hand, the Lipschitz property of $H$ gives a lower bound for the left hand side of \eqref{eq:proof_psi_lefthandside}. Indeed, all the arguments of the functions $H$ in the inequality are bounded in absolute value by $14L_T+1$.
It yields
\begin{align*}
 &\sigma +4\mathcal{C}_H^2\alpha\e^T   +\alpha\frac{\e^{ t^*}}{2}(\xstar^2+\xistar^2) -2\mathcal{C}_H \left(\alpha\e^{ t^*} (|\xstar|+|\xistar|) + \frac{\dx^{1/2}}{2} + \alpha \frac{\e^{ T}}{2}\dx  \right)
  \\ \le & \;\mathcal{C}(k)\left( \dt^{1/2} +\dx^{1/2}+ \|J_\dt^k-J_0^k\|_\infty\right), 
\end{align*}
where $\mathcal{C}(k)$ is a constant depending on $k$, and on the parameters $K$ and $L_T$.
We remark now that the left-hand side of the inequality is bounded from below independently of $|\xstar|$ and $|\xistar|$.  Hence, 
\begin{align*}
 \sigma \le \sigma +4\mathcal{C}_H^2\alpha\e^T - 4\mathcal{C}_H^2 \alpha \e^{ t^*} \le \mathcal{C}(k) \left( \dt^{1/2} +\dx^{1/2}+ \|J_0^k-J_\dt^k\|_\infty \right)+\alpha \mathcal{C}_H \e^{ T} \dx, 
\end{align*}
and $\sigma=\sigma(\dt,k)$ is defined so that the previous inequality cannot hold, and that $\sigma(\dt,k)\rightarrow_{\dt\to 0} 0$ when $k$ is fixed, as does the right hand side of the inequality. 
Because of $\|J^k_0 - J^k_\dt\|_\infty$, 
 there is no indication for the speed of
the convergence $\sigma(\dt,k)\rightarrow_{\dt\to 0} 0$ when $k$ is fixed. Indeed, Lemma \ref{lem:UniformCVJdeltak}-\ref{lem:UniformCVJdeltak_iii} is obtained by using a compactness argument, which does not give a quantitative estimate.
\end{proof}

We are now able to gather all these preliminary results to prove Prop. \ref{prop:scheme_0}:

\begin{proof}[Proof of Prop. \ref{prop:scheme_0}]
Consider a choice of  $\sigma=\sigma(\dt,k)$  as in Lemma \ref{lem:psi}. Then, the function $\psi$ defined in \eqref{eq:psi} reaches its maximum at $(t^*,x^*,\tau^*,\xi^*)$, therefore
\[
 \forall (t,x)\in[0,T[\times\R, \; \psi(t,x,t,x)\le \psi(t^*,x^*,\tau^*,x^*),
\]
hence, for all $\alpha\in(0,1)$, 
\begin{align*}
 v^k(t,x) -v_\dt^k(t,x) & \le 
 \sigma(\dt,k) t + 4\mathcal{C}_H^2\alpha \e^T t + \alpha \e^t  x^2 + \frac{\alpha}{T-t} 
 \\ &+  L_T |\xstar-\xistar| + (L_T^2+K)(t^*+\tau^*),
\end{align*}
thanks to Lemma  \ref{lem:vdtk}-\ref{lem:vdtk_space}-\ref{lem:vdtk_time}.
Let us start by letting $\alpha\to 0$ in the previous inequality, to get 
\[
 v^k(t,x) - v^k_\dt(t,x) \le \sigma(\dt,k) t + C\dt^{1/2},
\]
where $C>0$ can be determined using the fact that
$(x^*,\xi^*)\in\R^2$ satisfy \eqref{eq:psiproof_Lipschitz}, that the ratio $\dt/\dx$ is fixed, and that $t^*\le 2(L_T^2+K)\dt^{1/2}$. Thanks to \eqref{eq:psiproof_Lipschitz},  $\tau^*\le 4(L_T^2+K)\dt^{1/2}$ also holds. It is worth noticing that since $v^k$ and $v^k_\dt$ are continuous, this inequality also holds if $t=T$. 
Then,  Lemma \ref{lem:vdtk}-\ref{lem:vdtk_Monotonicity} yields
\[
 v^k(t,x) \le v_\dt(t,x) + \sigma(\dt,k) t + C\dt^{1/2},
\]
 Still considering a fixed $k>0$, let now $\dt\to 0$. As $\sigma(\dt,k)\to_{\dt\to 0}0$, and $v_\dt$ converges pointwise to $v_0$ (see Lemma \ref{lem:PointwiseCVofv_dtJ_dt}-\ref{lem:PointwiseCV_vdt}),
\[
 v^k(t,x) \le v_0(t,x).
\]
We conclude by noticing that this inequality is true for all $(t,x)\in[0,T]\times\R$. 
Finally, we let $k\to +\infty$, to get the following inequality
\[
 v^\infty\le v_0.
\]

Now, we have to prove the reverse inequality. The proof is very similar to what was done previously, but some modifications are necessary. We list here the modifications that are to be done in the steps of the proof:
\begin{itemize}
 \item In Lemma \ref{lem:PointwiseCVofv_dtJ_dt}-\ref{lem:PointwiseSV_Jdt}, 
 an upper semi-continuous representative should be opted for ($\mathfrak{J}_0$ instead of $J_0$, say).
 Note that $\mathfrak{J}_0=J_0$ almost everywhere.
 \item The functions $J_\dt$ and $\mathfrak{J}_0$ should be regularized \emph{from above} 
 instead of \eqref{eq:J_delta_k}. Namely, for $k>0$ and $\dt> 0$, let us define
 \begin{align*}
  \forall t\in[0,T],\;\mathfrak{J}_0^k (t) &= \sup\limits_{s\in[0,T]} \left( \mathfrak{J}_0(s) - k|t-s|\right),
  \\
  \mathfrak{J}_\dt^k (t) &= \sup\limits_{s\in[0,T]} \left( {J}_\dt(s) - k|t-s|\right).
 \end{align*} 
Most of the properties of Lemma \ref{lem:UniformCVJdeltak} still hold true, except that for fixed $\dt\ge 0$ and for all $t\in[0,T]$, $\mathfrak{J}^k_\dt(t) \searrow \mathfrak{J}_\dt(t)$ as $k\to +\infty$.  
Similarly, \ref{lem:UniformCVJdeltak_iii} has to be replaced by 
\[
 \text{for\;fixed\;} k>0, \; \|\mathfrak{J}_\dt^k - \mathfrak{J}_0^k\|_\infty\underset{\dt\to 0}\longrightarrow 0.
\]
\item The viscosity solution $w^k$ of the following Hamilton-Jacobi equation should be defined accordingly
\[
 \partial_t w^k + |\nabla_x w^k|^2 = -R(x,\mathfrak{J}^k_0), \; x\in\R, \; t>0,
\]
initialized with $v^\tin$.
The properties of Lemma \ref{lem:vdtk} are still true, except \eqref{eq:limit_lsc}. We have instead: $w^k\searrow w^\infty$ when $k\to+\infty$ pointwise in $[0,T]\times\R$, where $w^\infty$ is the viscosity solution of 
\begin{equation}
\label{eq:limit_usc}
 \partial_t w^\infty + |\nabla_x w^\infty|^2 = -R(x,\mathfrak{J}_0), \;x\in\R, t\in[0,T],
\end{equation}
initialized with $v^\tin$.
\item 
The regularized
scheme associated to $\mathfrak{J}^k_\dt$ 
should be defined as well, namely
\begin{equation}
 \label{eq:schemelimit_regularized_w}
 w^k_\dt(t_n+s) = \mathcal{M}_s(w^k_\dt(t_n,\cdot))(x)-s\;R(x,\mathfrak{J}_\dt^k(t_n+s)),
\end{equation}
for all $n\in\ccl 0,N_t-1\ccr$, $s\in(0,\dt]$ and $x\in\R$, with $\mathcal{M}_s$ defined in \eqref{scheme:Monotonic}, and initialized with $v^\tin$. The properties of Lemma \ref{lem:vdtk} are still true, except \ref{lem:vdtk_Monotonicity} that has to be replaced by $w^k_\dt\ge v_\dt$. 
\item Lemma \ref{lem:psi} should also be adapted. Instead of $\psi$, let us define 
\begin{align*}
 \Psi(t,x,\tau,\xi) &= w^k(t,x)-w^k_\dt(\tau,\xi) + \frac{(x-\xi)^2}{2\dx^{1/2}} +\frac{(t-\tau)^2}{2\dt^{1/2}} +(\sigma+4\mathcal{C}_H^2\alpha \e^T) t 
 \\&+ \alpha \frac{\e^t}{2}(x^2+\xi^2) + \frac{\alpha}{T-t}.
\end{align*}
Then, Lemma \ref{lem:psi} still holds true, but with a minimum instead of a maximum. 
\item As it has been done in the first part of this proof, we obtain eventually $w^\infty \ge v_0.$
\end{itemize}

To conclude, remark that $v^\infty$ and $w^\infty$ are respectively viscosity solution of \eqref{eq:limit_lsc} and \eqref{eq:limit_usc}, that are recalled here 
\begin{align*}
& \partial_t v^\infty + |\partial_x v^\infty|^2=-R(x,J_0), \;\;\;\; \partial_t w^\infty + |\partial_x w^\infty|^2=-R(x,\mathfrak{J}_0),
\end{align*}
both initialized with $v^\tin$, and with the source terms being such that $J_0=\mathfrak{J}_0$ $a. e$. Thanks to Theorem \ref{thm:uniqueness}-\ref{thm:uniqueness_unconstrained},
 it implies that  $v^\infty=w^\infty$ $a. e$. Then, the equality 
\[
v_0=v^\infty=w^\infty,
\]
comes immediately from  $v^\infty\le v_0\le w^\infty$, and because all these functions are continuous. Indeed, one can notice that the Lipschitz constants of Lemma \ref{lem:vdtk} do not depend on $k$ or $\dt$. Hence, $v_0$ enjoys the same Lipschitz regularity as $v_\dt^k$, and is, in particular, continuous.  

The next step consists in identifying $v_0$ and $J_0$ to the viscosity solution $v$ of \eqref{eq:limit}, and to the associated constraint $J$. It is a consequence of  $\min v_0=0$, proved in Lemma \ref{lem:PointwiseCVofv_dtJ_dt}-\ref{lem:PointwiseCV_vdt} and of Theorem \ref{thm:uniqueness}, meaning that
\[
v_0=v, \text{\;\;and\;\;}J_0=\mathfrak{J}_0=J\; a.e.
\]
Indeed, thanks to the assumptions made on the problem, $v$ is also continuous (see \cite{BarlesMirrahimiPerthame2009}), so the equality $v_0=v$ is  true pointwise in $(t,x)\in[0,T]\times\R$.

To conclude, remark that the only limit of the subsequences $(v_\dt)$ and $J_\dt$, defined in Lemma \ref{lem:PointwiseCVofv_dtJ_dt} are $v$ and $J$. The restriction \emph{up to a subsequence} stated in Remark \ref{rmq:subsequence} can then be removed, and Prop. \ref{prop:scheme_0} is proved.
\end{proof}

\section{Convergence of the scheme \eqref{scheme:epsilon}}
\label{sec:scheme_ep}

In this section, we fix $\ep>0$, and we prove that \eqref{scheme:epsilon} approximates properly \eqref{eq:u_epsilon} when the discretization parameters $\dt$ and $\dx$ go to $0$. We start with a technical lemma, which states properties of the sequences $(I^{n+1})_{n\in\ccl 0,N_t-1\ccr}$  and $(u^{n+1})_{n\in\ccl 0,N_t-1\ccr}$, with $u^{n+1}=(u^{n+1}_i)_{i\in\Z}$, defined by the scheme \eqref{scheme:epsilon}.

\begin{lem}
 \label{lem:stability_epfixed}
 Suppose that the assumptions \eqref{assumption:psi}-\eqref{assumption:R_extrema}-\eqref{assumption:R_bounded_decreasing}-\eqref{assumption:initial_data}-\eqref{assumption:u_increasing_infty}-\eqref{assumption:u0Lipschitz} are satisfied, and that $\ep>0$ is fixed. There exist  $\dx_0>0$, and  $I_{M'}>0$, depending on $\ep$, and on the constants arising in the assumptions, such that 
 if $\dt$ and $\dx<\dx_0$, are fixed such that \eqref{eq:epfixed_CFL} holds,
  the scheme \eqref{scheme:epsilon} is well-defined. 
 Moreover, the sequence $(u^n_i)_{n,i}$ defined by the scheme \eqref{scheme:epsilon} satisfies:
 \begin{enumerate}
  \item \label{lem:epfixed_Lipschitz} For all $n\in\ccl 0,N_t\ccr$, there exists a constant $\lambda_n=L_0+n\dt \kappa \le L_0 + T\kappa=\lambda_{N_t}$, with $L_0$ defined in \eqref{assumption:u0Lipschitz} and $\kappa$ in \eqref{eq:kappa}, such that the sequence $(u^n_i)_{i\in\Z}$ enjoys $\lambda_n$-Lipschitz property 
  \[
   \forall i\in\Z, \; \left|\frac{u^n_i-u^n_{i-1}}{\dx}\right|\le\lambda_n.
  \]
\item \label{lem:epfixed_boundedbelow} For all $n\in\ccl 0,N_t\ccr$, there exists $\underline{\beta}_n\in\R$, with $\underline{\beta}_n\ge \underline{\beta}_{N_t}=\underline{b}-TH(\underline{a},\underline{a}) - T\|R(\cdot,0)\|_\infty$, and 
$\overline{\beta}_n\le \overline{\beta}_{N_t} = \overline{b} +\overline{a}N_t\dx_0 + T\kappa$, such that
for all $i\in\Z$, 
\[
 \underline{a}|x_i-x_0| + \underline{\beta}_n\le u^n_i  
 \le \overline{a}|x_i-x_0|+ \overline{\beta}_n
 ,
\]
where $\underline{a}$ and $\underline{b}$ have been defined in \eqref{assumption:u_increasing_infty}, $\kappa$ in \eqref{eq:kappa} and $T$ is the fixed final time. 
\item \label{lem:epfixed_boundsI} For all $n\in\ccl 0,N_t-1\ccr$,
$
 0\le I\le I_{M'}.
$
\end{enumerate}
\end{lem}

\begin{rmq}
 At first sight, this lemma is similar to Lemma \ref{lem:AP}. However, it  holds here for a fixed $\ep>0$, and it states uniform estimates in $\dt$ and $\dx<\dx_0$ such that \eqref{eq:epfixed_CFL} is satisfied. On the contrary, Lemma \ref{lem:AP} states uniform estimates in $\ep\in(0,\ep_0)$, where $\ep_0>0$ depends on the assumptions, and on $\dt$ and $\dx$. 
\end{rmq}

\begin{rmq}
 In what follows, $I_{M'}$ is chosen large enough, such that 
 \begin{equation}
 \label{eq:epfixed_boundIep} 
 \forall t\in[0,T], \;0< I_\ep(t) + 1 \le I_{M'},
 \end{equation}
 where $I_\ep$ is defined in \eqref{eq:u_epsilon}. We refer to  \cite{BarlesMirrahimiPerthame2009} for the existence of such a bound.
\end{rmq}

\begin{proof} 
 The proof is done by induction. The initial data $u^0=(u^0_i)_{i\in\Z}$ enjoys the properties of Lemma \eqref{lem:stability_epfixed}. Let $\ep>0$ be fixed, and let us suppose that the items \ref{lem:epfixed_Lipschitz}-\ref{lem:epfixed_boundedbelow} are satified by $u^n=(u^n_i)_{i\in\Z}$ for a given $n\in\ccl 0,N_t-1\ccr$, and prove that $I^{n+1}$ and $u^{n+1}=(u^{n+1}_i)_{i\in\Z}$ are well defined, and satisfy \ref{lem:epfixed_Lipschitz}-\ref{lem:epfixed_boundedbelow}-\ref{lem:epfixed_boundsI}. 
 
 First of all, let us remark that $I^{n+1}$ is solution of $\Phi(I)=0$, with
  \begin{equation}
  \label{eq:Phi}
   \Phi(I)=I-\dx\sum\limits_{i\in\Z} \psi(x_i) \e^{-\mathcal{M}^\ep_\dx(u^n)_i/\ep}\e^{\dt R(x_i,I)/\ep},
  \end{equation}
  where $\mathcal{M}^\ep_\dx$ is defined in \eqref{scheme_ep:MonotonicFormulation}. Thanks to  Lemma \ref{lem:Monotonic_ep}, and because of \eqref{eq:epfixed_CFL}, 
  \begin{equation}
  \label{eq:epfixed_proof_boundbelow}
   \forall i \in\Z, \; \mathcal{M}^\ep_\dt(u^n)_i\ge \underline{a}|x_i-x_0| + \underline{\beta}_n - \dt H(\underline{a},\underline{a})\ge \underline{a}|x_i-x_0| + \underline{\beta}_{N_t},
  \end{equation}
so that the sum in \eqref{eq:Phi} is well-defined for all $I\in\R$.  Since $\Phi$ is a difference between an increasing and a decreasing function, there exists a unique $I^{n+1}\in\R$ such that $\Phi(I^{n+1})=0$. Therefore, $u^{n+1}$ is uniquely determined too. Moreover, the inequality $\Phi(I)\le I$ immediately yields that $I^{n+1}\ge 0$. As $R(x,\cdot)$ is decreasing for all $x$,  
\[
 \forall i\in\Z, \; u^{n+1}_i\ge  \mathcal{M}^\ep_\dt(u^n)_i -\dt R(x_i,0),
\]
which gives the lower estimate in \ref{lem:epfixed_boundedbelow}, with $\underline{\beta}_{n+1}=\underline{\beta}_n - \dt H(\underline{a},\underline{a}) - \dt \|R(\cdot,0)\|_\infty.$

Now, consider  $I\ge I_M$, with $I_M$  defined in \eqref{assumption:R_extrema}. Then,  $R(x_i,I)\le 0$ for all $i\in \Z$, thanks to \eqref{assumption:R_extrema}-\eqref{assumption:R_bounded_decreasing}, and using \eqref{assumption:psi}-\eqref{eq:epfixed_proof_boundbelow} we have
\begin{equation}
\label{eq:bound_I}
 \Phi(I)\ge I-2 \psi_M \e^{-\underline{\beta}_{N_t}/\ep} \dx 
 \frac{1}{1-\e^{-a\dx/\ep}}  \underset{\dx\to 0}\longrightarrow I- \frac{2\ep}{a} \psi_M \e^{-\underline{\beta}_{N_t}/\ep}.
\end{equation}
Hence, there exists $\dx_0>0$ and $I_{M'}>0$ such that for all $\dx\le \dx_0$, $\Phi(I_{M'})>0$. Since $\Phi$ is increasing, $I^{n+1}\le I_{M'}$. 
Eventually, Lemma \ref{lem:Monotonic_ep} yields \ref{lem:epfixed_Lipschitz}, and
\[
 \forall i\in\Z, \;\mathcal{M}^\ep_\dt(u^n)_i\le \overline{a}|x_i-x_0|+ \overline{\beta}_n + 2\ep\overline{a} \frac{\dt}{\dx},
\]
where  \eqref{eq:epfixed_CFL} gives \ref{lem:epfixed_boundedbelow}.
\end{proof}

\begin{rmq}
\label{rmq:bound_I}
 It is worth noticing that $\dx_0$ and $I_{M'}$ are determined once for all and do not depend on $n\in\ccl 0,N_t\ccr$. Indeed, coming back to the definition of $\underline{\beta}_{N_t}$, one can remark that they can be fixed independently of the induction. However, they depend on $\ep$, which is fixed here. Their asymptotic behavior when $\ep\to 0$ is not satisfactory, since $\dx_0$ may vanish, and $I_{M'}$ grows to infinity, when $\ep\to 0$, as $\underline{\beta}_{N_t}$ might be negative. We refer to Lemma \ref{lem:AP} for a bound of $I^{n+1}$ independent of $\ep$, with fixed $\dt$ and $\dx$. Indeed, this bound is valid for small $\ep$, and the bound for $I^{n+1}$ outside of the asymptotic regime comes from \eqref{eq:bound_I}.
\end{rmq}

Going on with the proof of the convergence of scheme \eqref{scheme:epsilon}, its implicit character has to be dealt with. To this end, let us define,
\begin{equation}
 \label{eq:D_I} 
 \mathcal{D}_{I_{M'}} = \left\{ u=(u_i)_{i\in\Z}\in\R^\Z, \;u>\underline{u}, \; \dx\sum\limits_{i\in\Z}\psi(x_i)\e^{-u_i/\ep}< I_{M'}\right\},
\end{equation}
where
\[ \underline{u}=\left(\underline{u}_i\right)_{i\in\Z}=\left(\underline{a}|x_i-x_0| + \underline{\beta}_{N_t}\right)_{i\in\Z},
\]
 and
 $I_{M'}$ are defined in Lemma \ref{lem:stability_epfixed}. 
 Define then
 $\mathcal{S}_{I_{M'}}: \mathcal{D}_{I_{M'}}  \to \R^\Z$, such that 
\begin{equation}
 \label{eq:S_I} 
 \forall u=(u_i)_{i\in\Z}\in\mathcal{D}_{I_{M'}},\; \forall i\in\Z,\;
 \mathcal{S}_{I_{M'}}(u)_i= u_i+\dt R\left(x_i,\dx\sum\limits_{i\in\Z} \psi(x_i)\e^{-u_i/\ep}\right).
\end{equation}
Since $I\mapsto R(x,I)$ is smooth for all $x\in\R$, one can notice that $\mathcal{S}_{I_{M'}}\in\mathcal{C}^1\left(\mathcal{D}_{I_{M'}}\right)$. But  a stronger result holds:

\begin{lem}
\label{lem:scheme_ep_S_Lipschitz} 
Let $\dt >0$ and $\ep >0$. Suppose that $\dt K I_{M'}<\ep$, with $I_{M'}$  defined in Lemma \ref{lem:stability_epfixed}, and $K$ in \eqref{assumption:R_bounded_decreasing}. Then, $\mathcal{S}_{I_{M'}}: \mathcal{D}_{I_{M'}}\to \mathcal{S}_{I_{M'}} \left(\mathcal{D}_{I_{M'}}\right)$ is invertible. Moreover, its inverse enjoys Lispschitz regularity: for all $u,v\in\mathcal{S}_{I_{M'}} \left(\mathcal{D}_{I_{M'}}\right)$ such that $u-v\in\ell^\infty(\Z)$, $\mathcal{S}_{I_{M'}}^{-1} (u)-\mathcal{S}_{I_{M'}}^{-1}(v)\in\ell^\infty(\Z)$ and
\[
 \left\| \mathcal{S}_{I_{M'}}^{-1} (u)-\mathcal{S}_{I_{M'}}^{-1}(v)\right\|_\infty \le \frac{ 1}{ 1-\dt K I_{M'}/\ep}\|u-v\|_{\infty}.
\]
\end{lem}
As this Lemma is an elementary consequence of the implicit function theorem, its proof is not detailed here. These technical lemmas and Lemma \ref{lem:Monotonic_ep} yield Prop. \ref{prop:scheme_ep}.

\begin{proof}[Proof of Prop. \ref{prop:scheme_ep}]
Scheme \eqref{scheme:epsilon} can be rewritten using Lemma \ref{lem:stability_epfixed}, and notations \eqref{scheme_ep:MonotonicFormulation_utilde}-\eqref{eq:S_I}. Indeed, for all $n\in\ccl 0,N_t-1\ccr$, $u^{n+1}=(u^{n+1}_i)_{i\in\Z}\in\mathcal{D}_{I_{M'}}\cap\mathcal{S}_{I_{M'}}\left(\mathcal{D}_{I_{M'}}\right)$, and it is defined by induction with
\[
 \mathcal{S}_{I_{M'}}(u^{n+1})= \mathcal{M}^\ep_\dt(u^n).
\]
Considering $u_\ep$ and $I_\ep$ defined as the solution of \eqref{eq:u_epsilon}, the consistency error $E^{n+1}_i$ of the scheme \eqref{scheme:epsilon} at $(t_{n+1},x_i)$, with $n\in\ccl 0,N_t-1\ccr$ and $i\in\Z$, is defined by
\[
 E^{n+1}_i = u(t_{n+1},x_i) - \mathcal{M}_\dt^\ep\left(u(t_n,x_j\right)_{j\in\Z})_i + \dt R\left(x_i, I_\ep(t_{n+1})\right),
\]
and there exists a constant $C$ depending only on $\|\partial_t^2 u_\ep\|_{\infty,[0,T]\times\R}$, and $\|\partial_x^k u_\ep\|_{\infty,[0,T]\times\R}$ for $k=1,2,3$, such that
\begin{equation}
\label{eq:epfixed_ConsistencyError}
 \forall n\in\ccl 0,N_t-1\ccr, \;\forall i\in\Z, \; \left|E^{n+1}_i\right|\le C\dt(\dt+\dx).
\end{equation}
Apart from the finite-differences approximations of the derivatives, the scheme \eqref{scheme:epsilon} is constructed with a quadrature rule for the approximation of $I_\ep$. Its precision can be estimated, thanks to Lemma \ref{lem:stability_epfixed}. At first,  define 
a truncated version of $I_\ep(t)$, on a domain $[x_0-\mathcal{X},x_0+\mathcal{X}]$, by
\begin{equation}
 \label{eq:Iep_x}
 I_\ep^\mathcal{X}(t) = \int_{|x-x_0|\le \mathcal{X}} \psi(x) \e^{-u_\ep(t,x)/\ep} \dd x.
\end{equation}
Hence, $\mathcal{X}$ is determined such that, for all $t\in[0,T]$,
\begin{equation}
\label{eq:epfixed_approximationError}
 \left|I_\ep(t) - I^\mathcal{X}_\ep(t) \right| \le\dt
 , \text{\;\;and\;\;} \left| 
 \dx \sum\limits_{\genfrac{}{}{0pt}{1}{i\in\Z}{ |x_i-x_0|> \mathcal{X}}} \psi(x_i) \e^{-u_\ep(t,x_i)/\ep} 
 \right| \le \dt.
\end{equation}
Note that $\mathcal{X}$ can be chosen once for all, and independently of $\ep$,  remarking that, for all $t\in [0,T]$, $u_\ep(t,\cdot)$ is increasing at infinity. Indeed, thanks to \cite{BarlesMirrahimiPerthame2009}, the following estimate holds
\begin{equation}
\label{eq:min_uep}
 \forall t\in [0,T], \;\forall x\in\R, \; u_\ep(t,x) \ge \underline{a}|x-x_0|+ \underline{b}_{N_t}. 
\end{equation}
where we used 
the notations of Lemma \ref{lem:AP}.
Of course, such a choice makes $\mathcal{X}$ depend on $\dt$. Explicit computations using Lemma \ref{lem:stability_epfixed}-\ref{lem:epfixed_boundedbelow} and \eqref{eq:min_uep}, yield 
\begin{equation}
\label{eq:epfixed_AsymptoticBehaviorX}
 \mathcal{X}\underset{\dt\to 0}= \mathcal{O}(-\ln(\dt)),
\end{equation}
where we consider that $\ep>0$ is fixed. 
Note that $\mathcal{X}$ 
is 
such that for all $t\in[0,T]$,
\begin{equation}
\label{eq:epfixed_quadratureError}
 \left|I^\mathcal{X}_\ep(t)- \dx \sum\limits_{\genfrac{}{}{0pt}{1}{i\in\Z}{ |x_i-x_0|\le \mathcal{X}}} 
 \psi(x_i) \e^{-u_\ep(t,x_i)/\ep}
 \right| \le K\left( \mathcal{X} \dx^2+\dx\dt\right).
\end{equation}
Indeed, the approximation of the integral can be considered as if it were done with a trapezoidal rule, up to an error of order $\dt$ (adding half the sum of the two first neglected terms,
which are of size $\dt$). The error estimate of the trapezoidal rule yields that $K$ depends on the second derivative of $\psi \exp(-u_\ep(t,\cdot)/\ep)$, which is uniformly bounded with respect to $t\in[0,T]$. 
Suppose now that the ratio in  \eqref{eq:epfixed_CFL} is fixed. Then,  $\dx$ is uniquely determined for any given $\dt>0$, and
\[
 \dx \underset{\dt\to 0}= \mathcal{O}(\sqrt{\dt}). 
\]
Hence, thanks to \eqref{eq:epfixed_boundIep}-\eqref{eq:epfixed_approximationError}-\eqref{eq:epfixed_AsymptoticBehaviorX} and \eqref{eq:epfixed_quadratureError}, there exists $\dt_0>0$ such that for all $\dt<\dt_0$ and for all $t\in[0,T]$, 
$ \left(u(t,x_i)\right)_{i\in\Z}\in\mathcal{D}_{I_{M'}}$. 
Then, using  Lemmas \ref{lem:Monotonic_ep} and \ref{lem:scheme_ep_S_Lipschitz}, and the above estimates, there exists a constant, denoted $C(\ep)$, such that, for all $n\in\ccl 0,N_t-1\ccr$,
\[
 \left(1-\frac{\dt K I_{M'}}{\ep}\right)\|u_\ep(t_{n+1},x_j)_{j\in\Z} - u^{n+1}\|_\infty \le 
 \left\| (u_\ep(t_n,x_j)_{j\in\Z}   -  u^n  \right\|_\infty + C(\ep)\dt (\dt + \dx +  |\ln(\dt)| \dt). 
\]
Indeed, as $u^0_i=u_\ep(0,x_i)$ for all $i\in\Z$, the previous inequality yields that for all $n\in\ccl 0,N_t\ccr$, $(u_\ep(t_n,x_j)_{j\in\Z}   -  u^n\in\ell^\infty(\Z)$.
The first estimate of Prop. \ref{prop:scheme_ep} follows immediately. 
Eventually, one can notice that 
\[
 \mathcal{D}_{I_{M'}}\ni u=(u_i)_{i\in\Z}\mapsto \dx \sum\limits_{i\in\Z}\psi(x_i)\e^{-u_i/\ep}, 
\]
enjoys $I_{M'}/\ep$-Lipschitz regularity. This yields the second estimate of Prop. \ref{prop:scheme_ep}. 
\end{proof}

\begin{rmq}
\label{rmq:truncation}
Let us end with a remark about the implementation of \eqref{scheme:epsilon}. Its implicit character has been discussed in Section \ref{sec:results}, but another difficulty arises when coding it. Indeed, \eqref{scheme:epsilon} is defined for all indices $i\in\Z$, meaning that, in practice, the expressions have to be truncated. However, because of $I^{n+1}$, the expression of the scheme \eqref{scheme:epsilon} is nonlocal, in the sense that the whole distribution in trait $(u^n_i)_{i\in\Z}$ is needed to compute every single $u^{n+1}_i$ for $i\in\Z$. 
When implemented, the scheme \eqref{scheme:epsilon} uses an approximated value of $I^{n+1}$, with the truncation defined as previously. 

In addition, considering the scheme \eqref{scheme:epsilon} on a truncated domain raises questions about boundary conditions that are to be considered. Indeed, because of \eqref{scheme_ep:MonotonicFormulation_utilde}, the  $(u^n_i)$, for $|i|\le N+1$, are needed to compute $\mathcal{M}_\dt^\ep(u^n)_i$, for $|i|\le N$. In practice, $u^n_{-N-1}$ and $u^n_{N+1}$ can be approximated, we refer to Section \ref{sec:tests} for more details. Yet, to avoid more approximations, one can also define a truncation of $\mathcal{M}^\ep_\dt$
\[
\mathcal{M}_{\dt,N}^\ep: \R^{2(N+1)+1} \to \R^{2N+1},
\] 
such that 
\[
 \forall u=(u_j)_{j\in\Z}\in\R^\Z, \; \forall |i|\le N, \;  \mathcal{M}_{\dt,N}^\ep\left((u_j)_{|j|\le N+1}\right)_i = \mathcal{M}_\dt^\ep(u)_i. 
\]
Roughly speaking, this consists in avoiding the question of the boundary, by reducing the trait domain at each time step of the scheme. 
Note that $\mathcal{M}_{\dt,N}^\ep$ enjoys the same monotonicity properties as $\mathcal{M}_\dt^\ep$, and in particular the last point of Lemma \ref{lem:Monotonic_ep} can be easily adapted. Similarly, $\mathcal{S}_{I_{M'}}$ defined in \eqref{eq:S_I} can be defined on a truncated domain,  Lemmas \ref{lem:stability_epfixed} and \ref{lem:scheme_ep_S_Lipschitz} still hold, and Prop. \ref{prop:scheme_ep} is true in the truncated setting.
\end{rmq}

\section{Numerical tests}
\label{sec:tests}

In this section, we highlight and discuss the properties of the schemes \eqref{scheme:epsilon} and \eqref{scheme:limit} using numerical tests. Unless other choices are specified, we will consider the schemes in dimension $1$, with $\psi(x)\equiv1$ in \eqref{assumption:psi} and the initial data 
\begin{equation}
\label{eq:u_in}
 u^\tin(x)=v^\tin(x)= \frac{\min\left( (x-\beta)^2; (x-\alpha)^2 + \delta \right)}{\sqrt{1+x^2}},
\end{equation}
with $\alpha= 2$, $\beta=-0.2$ and $\delta =1$. This choice is adapted from \cite{BarlesPerthame} to satisfy the hypotheses \eqref{assumption:initial_data}-\eqref{assumption:u_increasing_infty}-\eqref{assumption:u0Lipschitz}-\eqref{assumption:min_u0}. We will also consider the function 
\begin{equation}
\label{eq:R}
 R(x,I)= \e^{-I} \frac{x^2}{1+x^2} - I,
\end{equation}
which satisfies \eqref{assumption:R_extrema} and \eqref{assumption:R_bounded_decreasing}. All the tests are done with final time $T=1$. In most cases, and if the discretization is not specified, we consider $\dt = 5\cdot10^{-4}$ and $\dx= 5\cdot10^{-2}$, such that \eqref{eq:epfixed_CFL}-\eqref{eq:CFL_AP}-\eqref{eq:CFL_0} are satisfied for all $\ep\in(0,1]$. First, the implementation of the schemes is done according to Remark \ref{rmq:truncation}. Namely, the iterations of the schemes are computed on a larger trait domain, that is reduced at each time iteration to avoid approximations at the boundary. 

The implementation of the schemes has been done using Matlab, the code is available at  \cite{Git}, where scripts for all the figures presented above are also provided. Note also that the solution of scheme \eqref{scheme:epsilon} will be denoted $u^\ep_\dt$ and $I^\ep_\dt$ in what follows. This choice is made to simplify the notations, and to be similar to $v_\dt$ and $J_\dt$ defined in scheme \eqref{scheme:limit}.

\subsection{Behavior of scheme \eqref{scheme:epsilon} when $\ep\to 0$}

\begin{figure}[!ht]
\begin{center}
\begin{tabular}{@{}c@{}c@{}}
\includegraphics[width=0.5\textwidth]{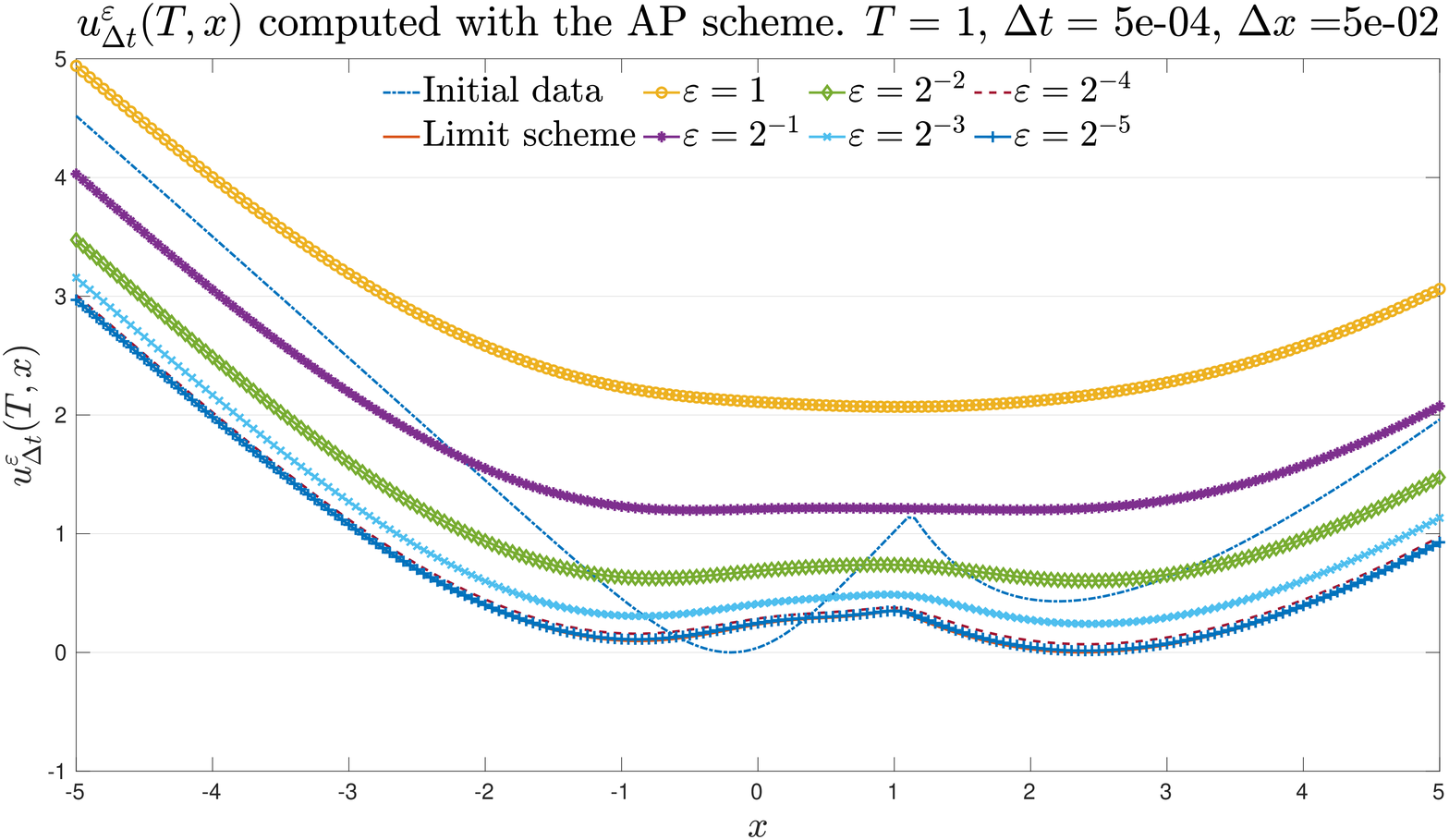}    &
\includegraphics[width=0.5\textwidth]{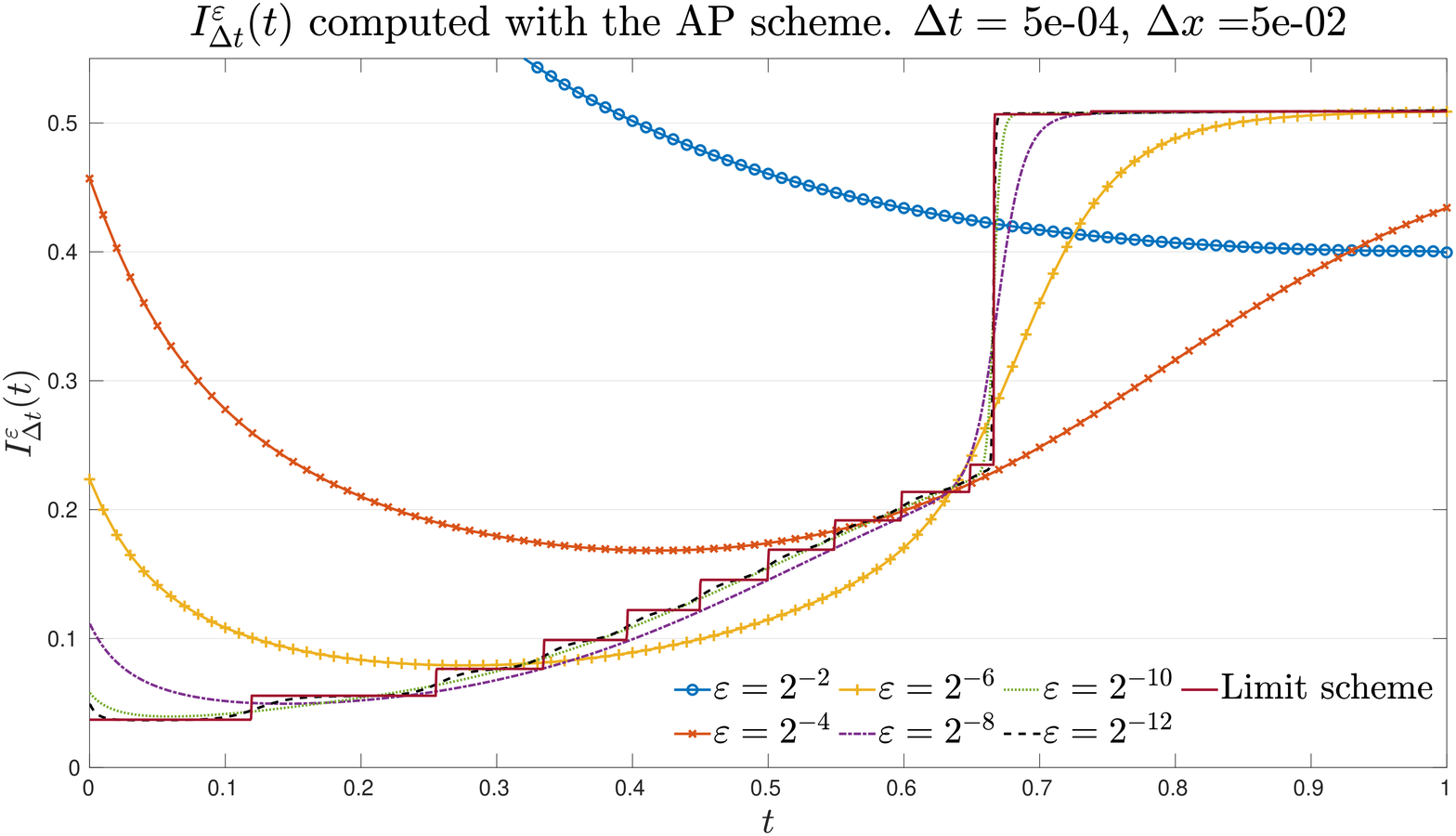}
\end{tabular}
\caption{$u^\ep_\dt$ (left) and $I^\ep_\dt$ (right) computed with \eqref{scheme:epsilon} for a series of $\ep$, and $v_\dt$ and $J_\dt$ computed with \eqref{scheme:limit}. Parameters: $T=1$, $\dx=5\cdot10^{-2}$, $\dt= 5\cdot 10^{-4}$, $u^\tin$ defined in \eqref{eq:u_in}, and $R$ in \eqref{eq:R}.}
\label{fig:stability}
\end{center}
\end{figure}

The behavior of scheme \eqref{scheme:epsilon} when $\ep\to 0$ is illustrated in Fig. \ref{fig:stability}, where $u^\ep_\dt$ and $I^\ep_\dt$, computed with \eqref{scheme:epsilon}, are displayed for a series of $\ep$. The limits $v_\dt$ and $J_\dt$ computed with \eqref{scheme:limit} are displayed on the same graph. As shown in Prop. \ref{prop:AP}, one can observe that the solution of scheme \eqref{scheme:epsilon} converges to the solution of scheme \eqref{scheme:limit} when $\ep\to 0$. It is also worth remarking that the solution $u^\ep_\dt$ of the problem \eqref{eq:u_epsilon} is smooth, and so is $I^\ep_\dt$ when $\ep>0$. Lipschitz singularities for $u^\ep_\dt$, and discontinuities for $I^\ep_\dt$, appear in the limit $\ep\to 0$. One can notice that $I^\ep_\dt$ is not necessarily increasing when $\ep>0$. Moreover, the convergence seems to be faster for $u^\ep_\dt$ than for $I^\ep_\dt$. 

\begin{figure}[!ht]
\begin{center}
\begin{tabular}{@{}c@{}c@{}}
\includegraphics[width=0.5\textwidth]{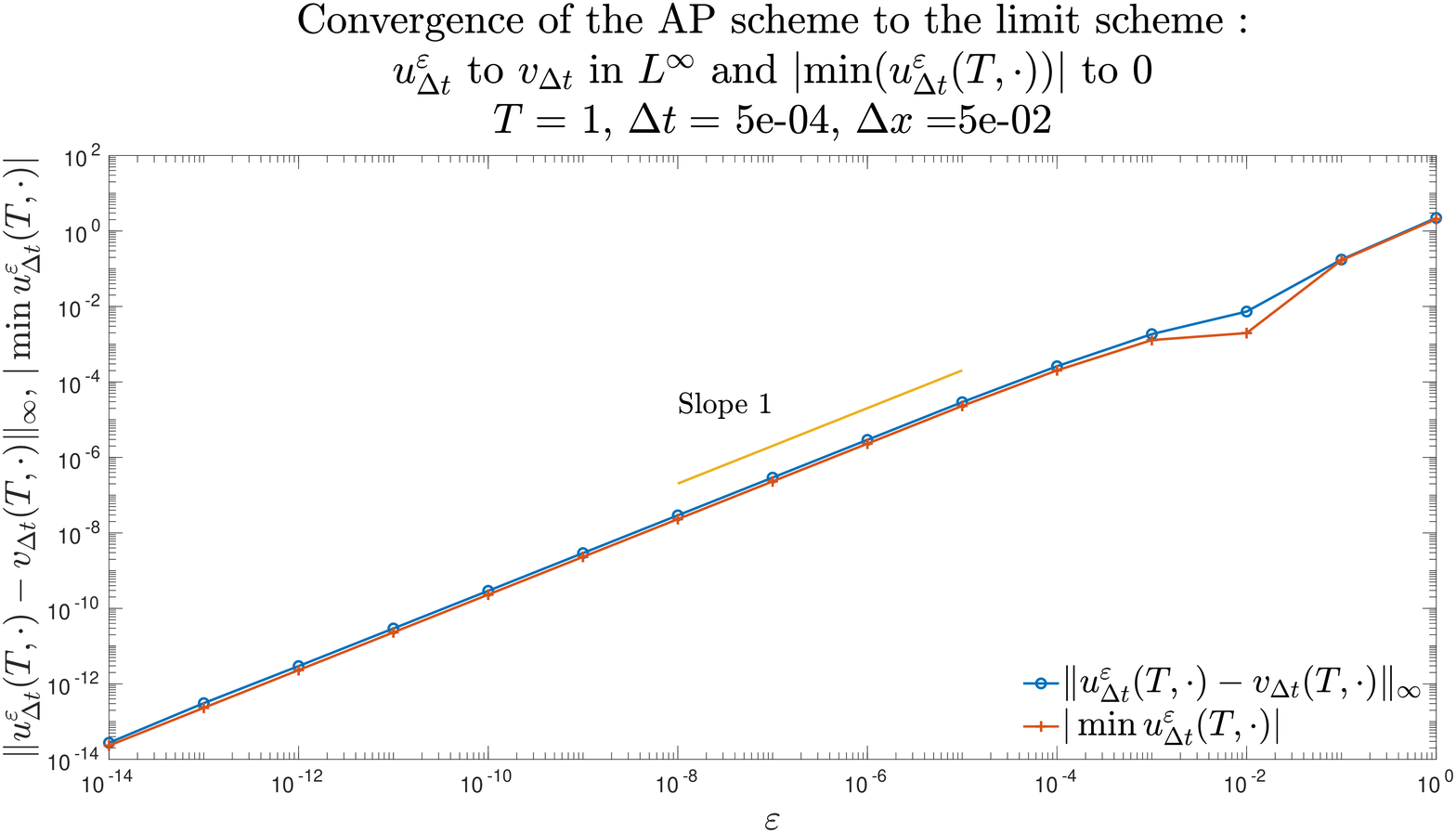}    &
\includegraphics[width=0.5\textwidth]{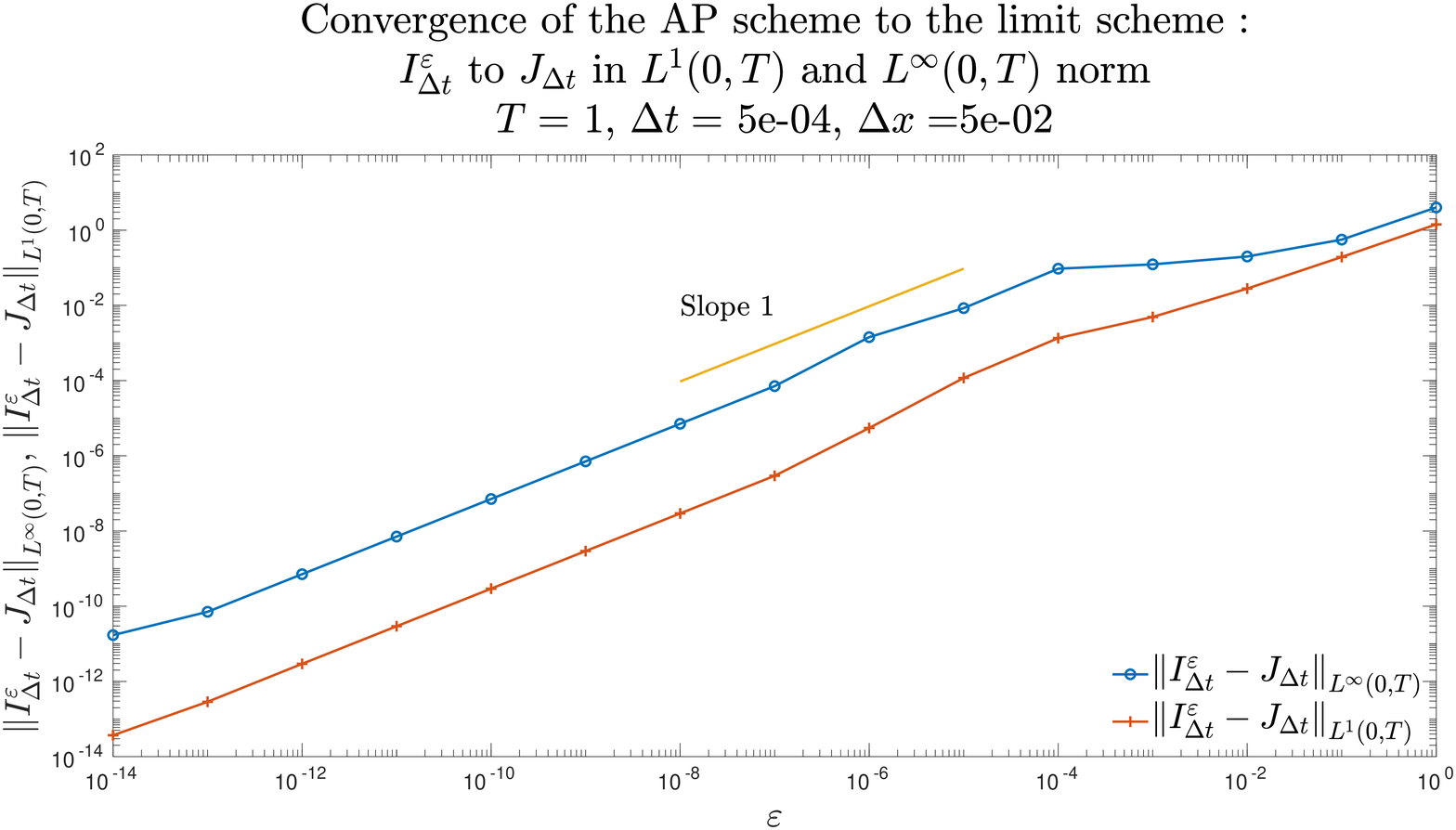}
\end{tabular}
\caption{Convergence of the solution of \eqref{eq:u_epsilon} to the solution of \eqref{eq:limit} . Left: $u^\ep_\dt$ to $v_\dt$ in $\LL^\infty$ norm, and $\min u^\ep_\dt$ to $0$ as functions of $\ep$ (logarithmic scale). Right: $I^\ep_\dt$ to $J_\dt$ in $\LL^1$ and $\LL^\infty$ norms, as functions of $\ep$ (logarithmic scale).
Parameters: $T=1$, $\dx=5\cdot10^{-2}$, $\dt= 5\cdot 10^{-4}$, $u^\tin$ defined in \eqref{eq:u_in}, and $R$ in \eqref{eq:R}.}
\label{fig:ConvergenceSpeed}
\end{center}
\end{figure}

More precisely, the convergence rate for $u^\ep_\dt$ and $I^\ep_\dt$ is numerically studied in Fig. \ref{fig:ConvergenceSpeed}. First of all, Lemma \ref{lem:AP}-\ref{lem:AP_minu} yields that the minimum of the approximation of $u^\ep_\dt$ with \eqref{scheme:epsilon} is of order $\ep$. This can indeed be observed on the left-hand side of Fig. \ref{fig:ConvergenceSpeed}, where the minimum of $u^\ep_\dt$ is plotted in logarithmic scale as a function of $\ep$. As expected, we observe a line which has slope $1$. This figure presents, on the same graph,  a numerical study of the convergence rate of the solution $u^\ep_\dt$ of \eqref{scheme:epsilon} to the solution $v_\dt$ of \eqref{scheme:limit} when $\ep\to 0$. The $\LL^\infty$ norm of $u^\ep_\dt(T,\cdot)-v_\dt(T,\cdot)$ is displayed in logarithmic scale as a function of $\ep$. This test suggests that the convergence of the solution $u_\ep$ of \eqref{eq:u_epsilon} to the solution $v$ of \eqref{eq:limit} is of order $1$ in $\ep$. Similarly, the convergence rate of $I^\ep_\dt$ to $J_\dt$ is studied in the right-hand side of Fig. \ref{fig:ConvergenceSpeed}, in discrete $\LL^1(0,T)$ and $\LL^\infty(0,T)$ norms. Once again, the rate of convergence seems to be $1$. 
However, we observe a discrepancy between the two tests in the regime $\ep\ge 10^{-4}$ which is  the order of the time step. 
The order of convergence is recovered in the regime $\ep\le 10^{-4}$, essentially because 
this convergence test is done for given $\dt$ and $\dx$, fixing the dimension of the problem.
We conclude from this observation that $\LL^1$  is more appropriate to capture the AP property due to the occurrence of true discontinuities of $J_\dt$. 
This behavior means a lack of uniform accuracy in $L^\infty$ norm, and we refer to Section \ref{sec:UA} for more details. Coming back to the continuous problems \eqref{eq:u_epsilon} and \eqref{eq:limit}, this   suggests that the convergence of $I_\ep$ to $J$ when $\ep\to0$ might be true in $L^1(0,T)$, but not in $L^\infty(0,T)$.

 \subsection{Behavior of scheme \eqref{scheme:limit}}
 
\begin{figure}[!ht]
\begin{center}
\begin{tabular}{@{}c@{}c@{}}
\includegraphics[width=0.5\textwidth]{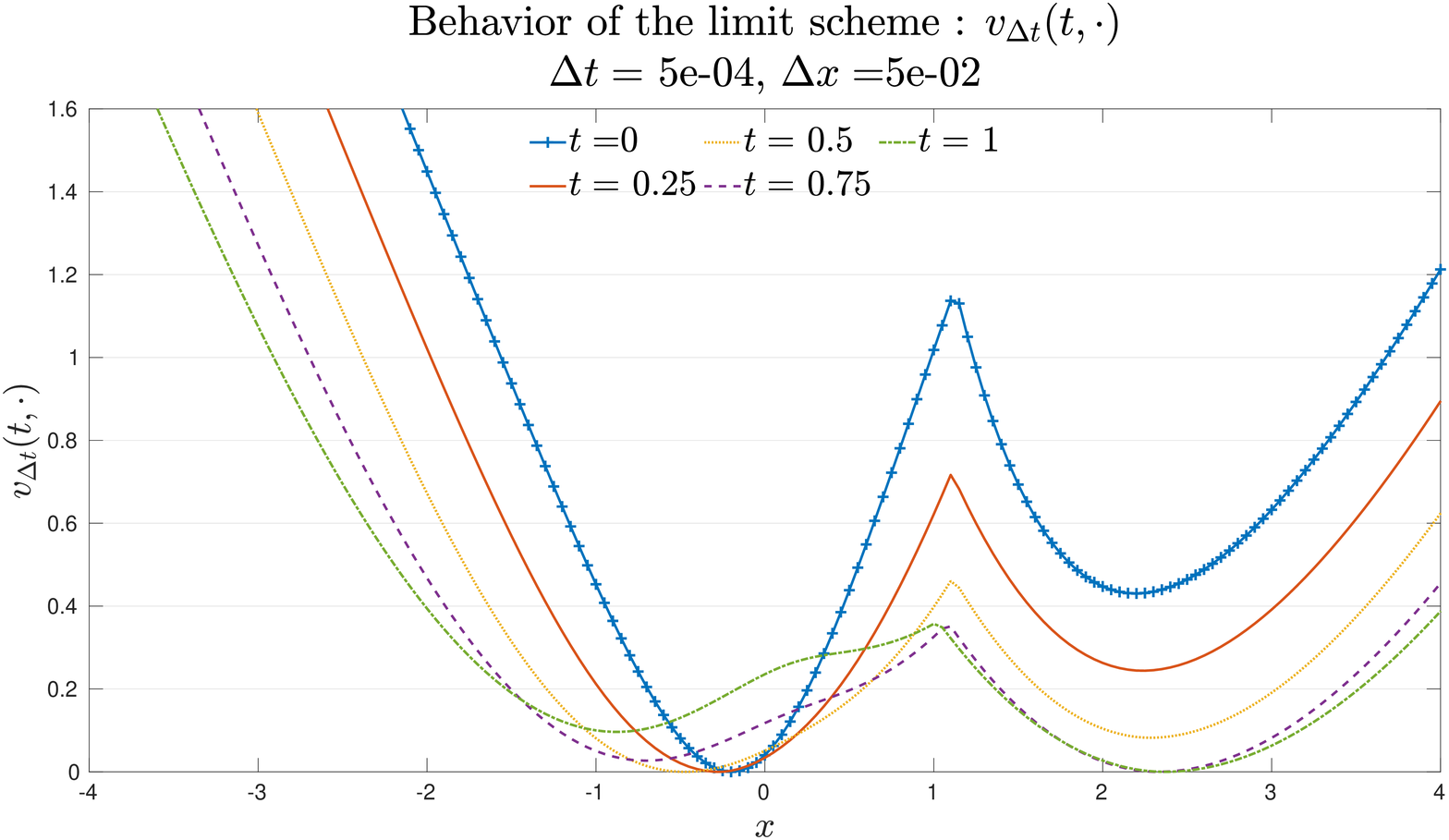}    &
\includegraphics[width=0.5\textwidth]{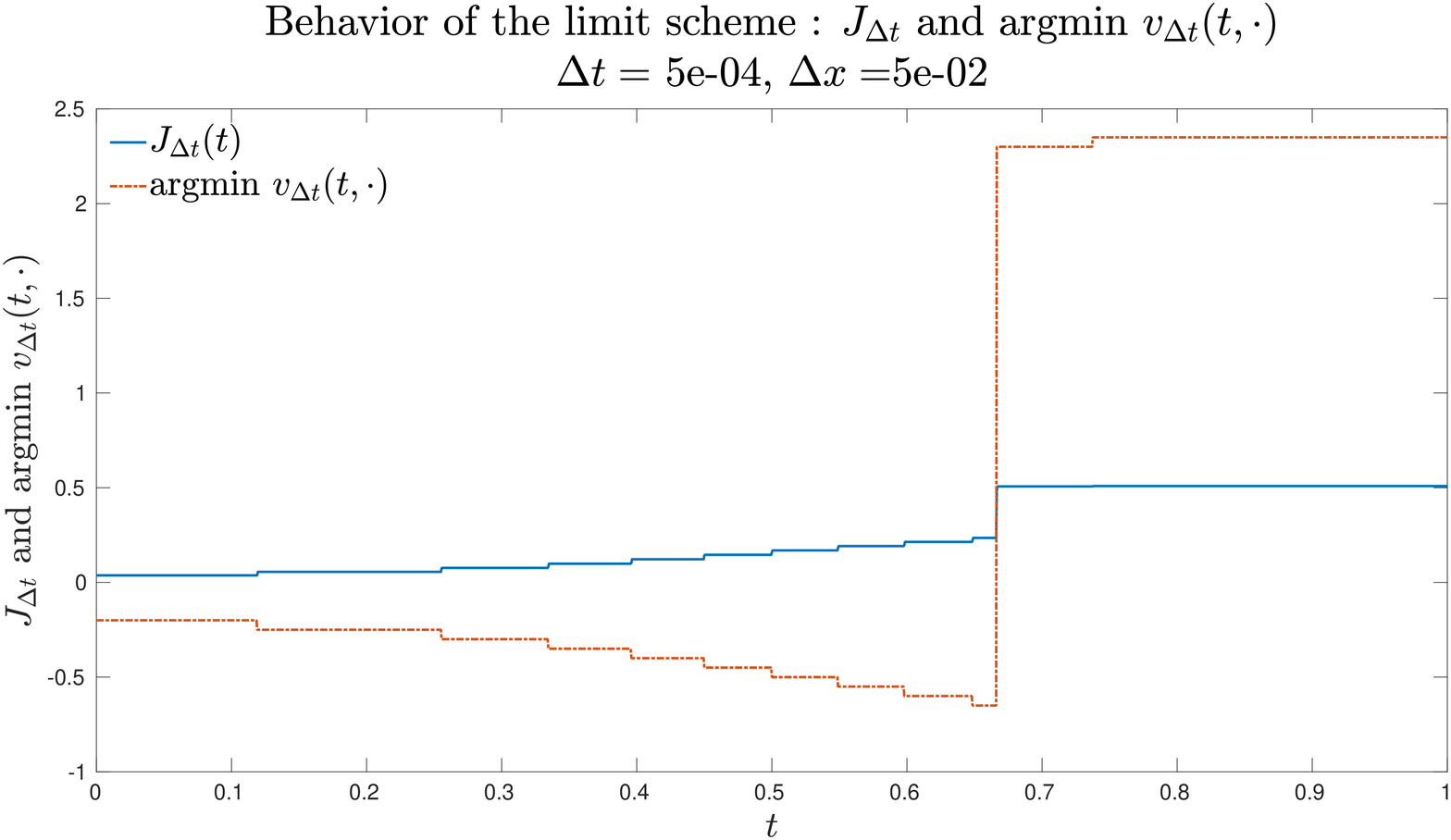}
\end{tabular}
\caption{Scheme \eqref{scheme:limit}. Left: $v_\dt$ computed with \eqref{scheme:limit} for a series of times. Right: $\arg\min_x v_\dt$ and $J_\dt$ as functions of $t$.
Parameters: $T=1$, $\dx=5\cdot10^{-2}$, $\dt= 5\cdot 10^{-4}$, $u^\tin$ defined in \eqref{eq:u_in}, and $R$ in \eqref{eq:R}.
}
\label{fig:LimitScheme}
\end{center}
\end{figure}

 We now discuss the behavior of scheme \eqref{scheme:limit}, regarding the lack of regularity of the solution of \eqref{eq:limit}. Indeed, $v$ enjoys Lipschitz regularity, while $J\in BV(0,T)$ can, in particular, have jumps. This behavior is highlighted in Fig. \ref{fig:LimitScheme}, where the left-hand side displays the solution $v_\dt$ of \eqref{scheme:limit} for some fixed times, as functions of $x$. 
 We emphasize the lack of diffusing effects, as shown by the $\mathcal{C}^1$ discontinuity of the solution $v_\dt$, which seems to be maintained as time grows. It is also interesting to notice that the function $J_\dt$ has true numerical jumps, where the solution varies considerably in a single time step, due to the implicit character of the scheme.
 Moreover, coming back to problem \eqref{eq:limit}, the selection makes the dominant trait, i.e.  $\bar{x}(t)$ such that of $v(t,\bar{x}(t))=\min v(t,\cdot)$, evolve in time. Left-hand side of Fig. \ref{fig:LimitScheme} 
 exhibits a case with a jump from the left local minimum to the right one.
 This is confirmed
 on the right-hand side of Fig. \ref{fig:LimitScheme}, where $\bar{x}_\dt$ and $J_\dt$ are displayed as functions of $t$. One can notice that the jumps occur simultaneously, which was to be expected since $J_\dt$ is a constraint that makes $\min v$ equal to $0$. Moreover, $J_\dt$ is nondecreasing on $[0,T]$, as it has been proved in Lemma \ref{lem:stability_scheme0}-\ref{lem:LimitTime_Jincreasing}.

\subsection{Truncated scheme}

As it has already been emphasized in Remark \ref{rmq:truncation}, the schemes \eqref{scheme:epsilon} and \eqref{scheme:limit} are nonlocal, meaning that the whole distribution in trait at time $t_n$ is needed to compute any single point at time $t_{n+1}$.
We proposed a way to restrict the schemes to a finite grid, by considering a larger trait domain at the initialization and removing points of the domain at each time iterations. The propositions of Section \ref{sec:results} hold true with this approximation, provided that the considered trait domain is large enough so that \eqref{eq:epfixed_approximationError} is satisfied. Thanks to this strategy, no approximation is required at the boundary. However, it is costly in terms of computational time, since $2N_t$ points in $x$ are to be added to the initial grid. This drawback can be dealt with in dimension $1$, but the cost increases with the dimension. Moreover, this stategy leads to complications when considering initial data which do not exactly satisfy \eqref{assumption:u_increasing_infty}-\eqref{assumption:u0Lipschitz}. Indeed, it would be natural to consider Gaussian distributions for the initial data of \eqref{eq:n_epsilon}, so that $u^\tin_\ep$ is quadratic. However, such distributions do not enjoy uniform Lipschitz property. When implemented, the conditions \eqref{eq:epfixed_CFL}-\eqref{eq:CFL_AP}-\eqref{eq:CFL_0} then have to be considered with the Lispschitz constant which is valid on the larger grid. It makes these stability conditions always more restrictive, as each point added in time makes the Lipschitz constant increase. 

To avoid this difficulty, we propose a truncated version of the schemes \eqref{scheme:epsilon}-\eqref{scheme:limit}. This consists in, once again, considering a truncated trait space $(x_i)_{i\in\ccl 1,N_x\ccr}$, such that \eqref{eq:epfixed_approximationError} is satisfied. However, this trait space is of constant size in all the time iterations. Since they are needed, the values at $x_0$ and $x_{N_x+1}$ are approximated. For all $n\in\ccl 0,N_t\ccr$, we propose the following approximation in \eqref{scheme:epsilon}-\eqref{scheme:limit}
\begin{equation} 
\label{eq:BoundaryApprox}
\begin{array}{l}
 \ds u^n_0= 4u^n_1 - 6 u^n_2 + 4 u^n_3 - u^n_4 \vspace{4pt} \\ 
 \ds u^n_{N_x+1} = 4u^n_{N_x} -6 u^n_{N_x-1} + 4u^n_{N_x-2} - u^n_{N_x-3},
 \end{array}
\end{equation}
which consists in extrapolating $(u^n_i)_{i\in\ccl 1,N_x\ccr}$ by a polynomial, whose derivatives coincide with the discrete derivatives of $(u^n_i)_{i\in\ccl 1,N_x\ccr}$. Namely, for the left point, we define
\[
 P(y) = u^n_1+ \frac{u^n_2-u^n_1}{\dx} y + \frac{u^n_1-2u^n_2+u^n_3}{\dx^2} y^2 + \frac{(u^n_2-2u^n_3 + u^n_4) - (u^n_1 - 2u^n_2 + u^n_3)}{\dx^3} y^3,
\]
is such that $P(0)$, $P'(0)$, $P''(0)$ and $P^{(3)}(0)$ coincide with the first discrete derivatives of $(u^n_i)_i$ that can be computed, and it satisfies $P(-\dx)=u^n_0$. A similar explanation holds for the right boundary. 

\begin{figure}[!ht]
\begin{center}
\begin{tabular}{@{}c@{}c@{}}
\includegraphics[width=0.5\textwidth]{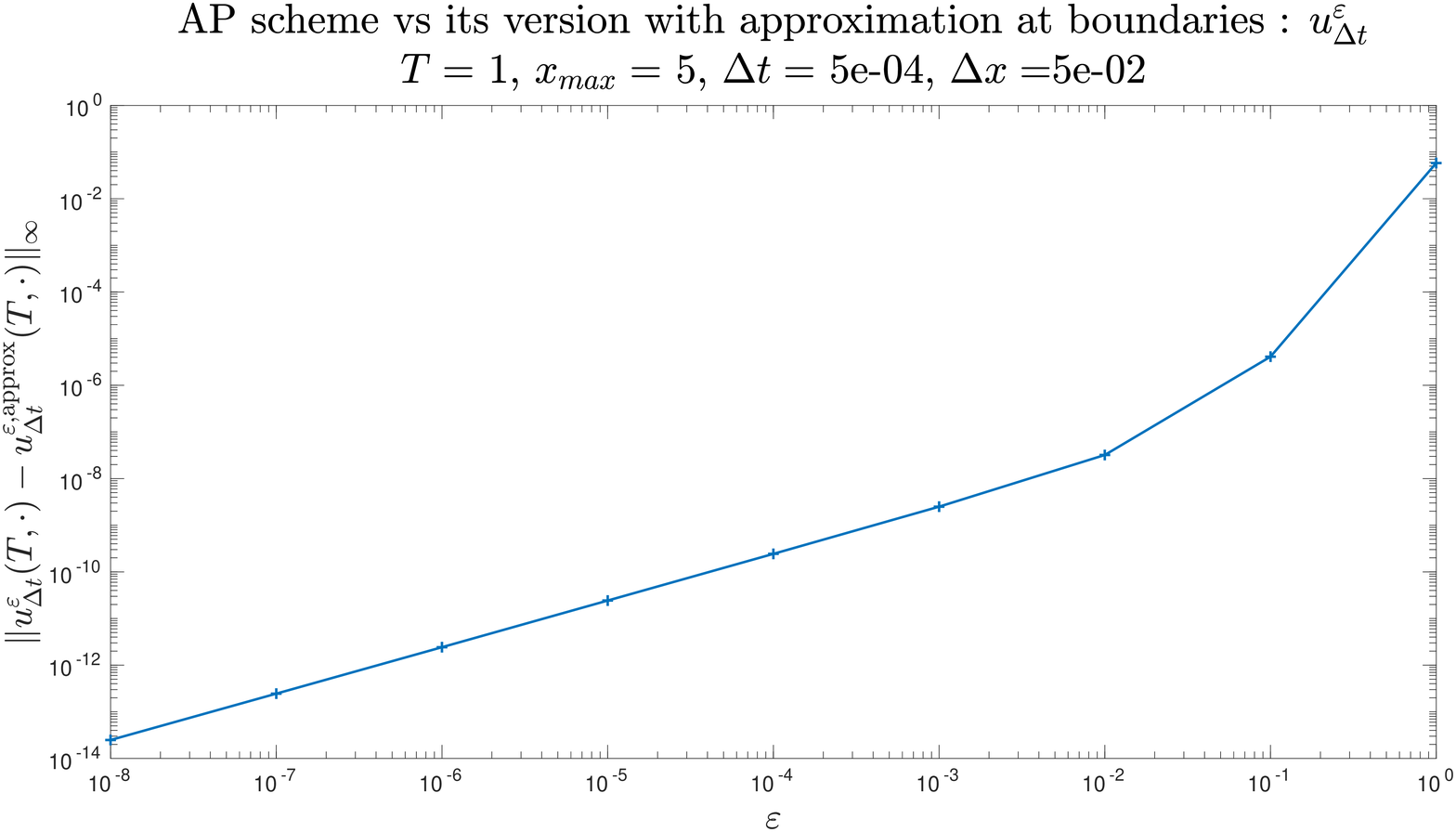}    &
\includegraphics[width=0.5\textwidth]{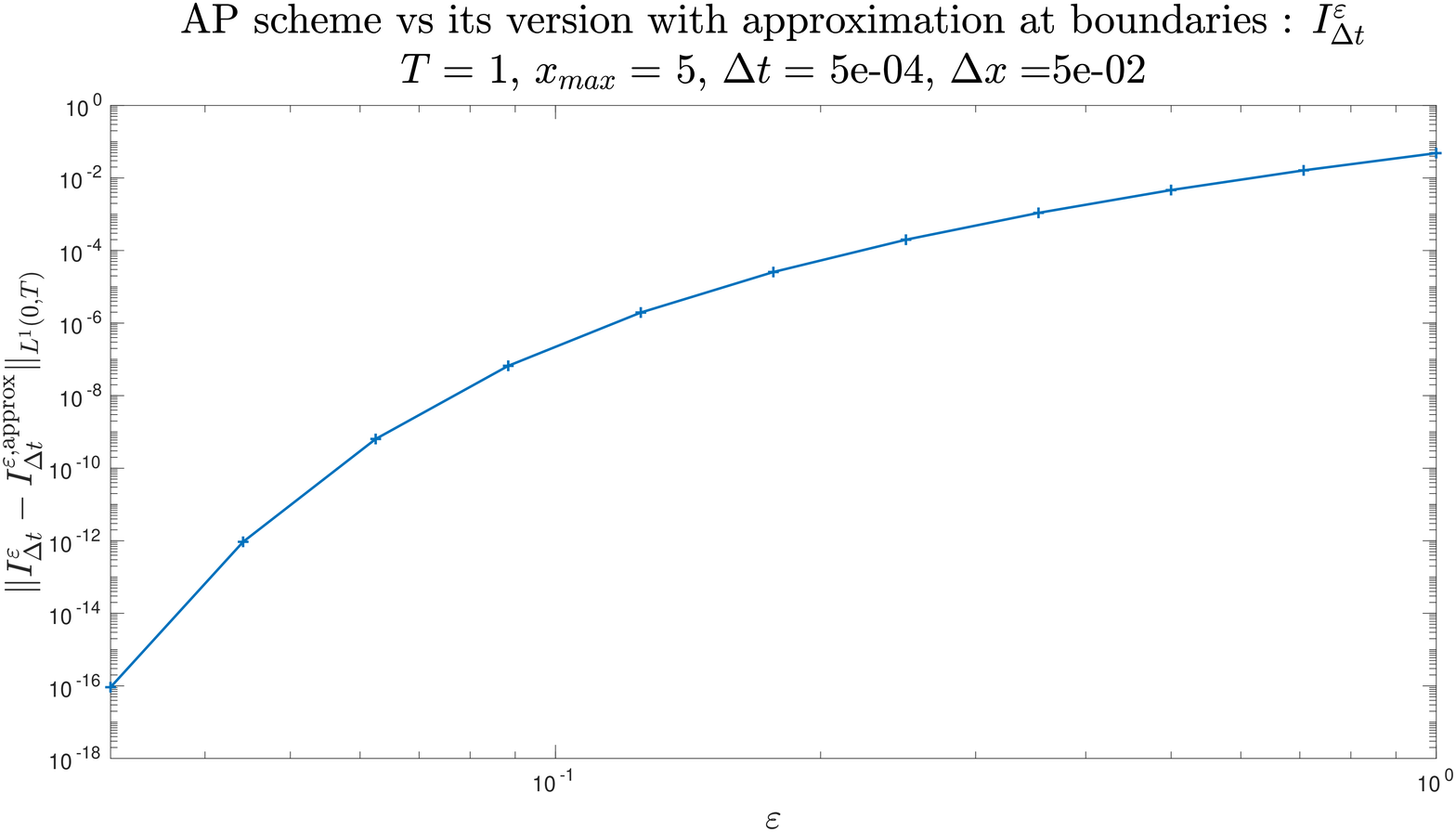}
\end{tabular}
\caption{Comparison between scheme \eqref{scheme:epsilon} and its version with approximations at boundaries, as function of $\ep$ (logarithmic scale). Left: $\LL^\infty$ norm in $x$ of the difference between the $u_\ep$ computed with the two schemes at final time. Right: difference between the two $I_\ep$ in $\LL^1(0,T)$ norm.  
Parameters: $T=1$, $\dx=5\cdot10^{-2}$, $\dt= 5\cdot 10^{-4}$, $u^\tin$ defined in \eqref{eq:u_in}, and $R$ in \eqref{eq:R}.
}
\label{fig:TruncatedScheme}
\end{center}
\end{figure}

This approximation is tested in Fig. \ref{fig:TruncatedScheme}, where results of scheme \eqref{scheme:epsilon} without and with the approximation at the boundary are compared. The left-hand side displays the $\LL^\infty$ norm in $x$ of the difference of  $u^\ep_\dt$ computed with the two versions of the scheme at time $T$, while the difference between the two $I^\ep_\dt$ in $\LL^1(0,T)$ norm is presented on the right-hand side.  In both cases, the results are presented as functions of $\ep$. Observe that the difference between \eqref{scheme:epsilon} and its version with approximations at boundaries is never greater than the discretization error. Moreover, this difference goes to $0$ when $\ep\to 0$, 
likely because the error due to the truncation in the quadrature step \eqref{eq:epfixed_quadratureError} is vanishing, combined with the fact that the characteristics lines are exiting the domain at $\ep=0$.
Thanks to the stability of \eqref{scheme:epsilon}, this validates numerically the approximation at boundary for \eqref{scheme:limit} as well.

\begin{rmq}
 As they are less expensive in terms of computational time, and since their results are very close to the results of schemes \eqref{scheme:epsilon}-\eqref{scheme:limit},  in what follows we will use 
 the corresponding schemes including the approximations \eqref{eq:BoundaryApprox} at the boundaries.
 \end{rmq}

\subsection{Accuracy of \eqref{scheme:limit}}

\begin{figure}[!ht]
\begin{center}
\begin{tabular}{@{}c@{}c@{}}
\includegraphics[width=0.5\textwidth]{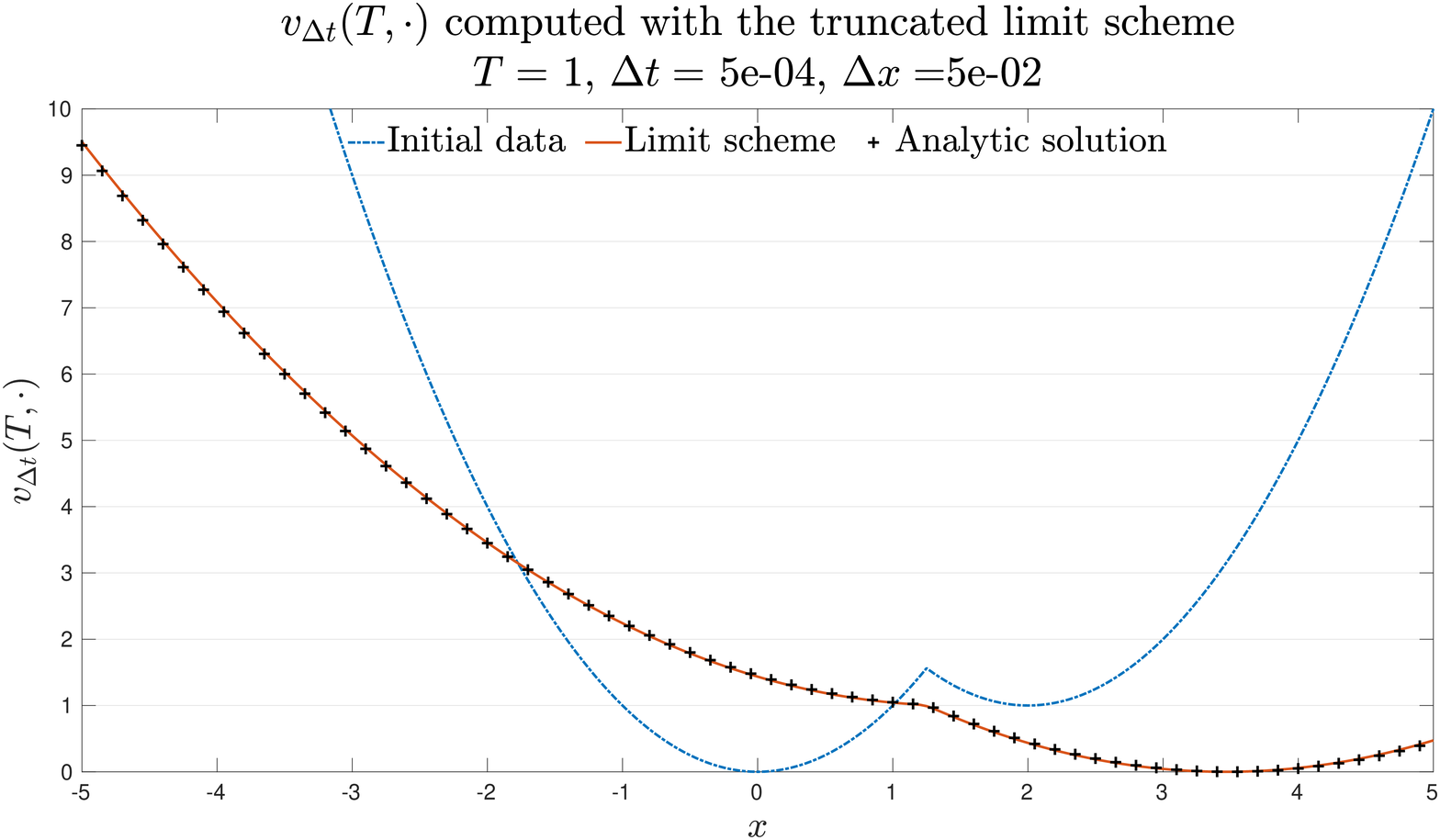}    &
\includegraphics[width=0.5\textwidth]{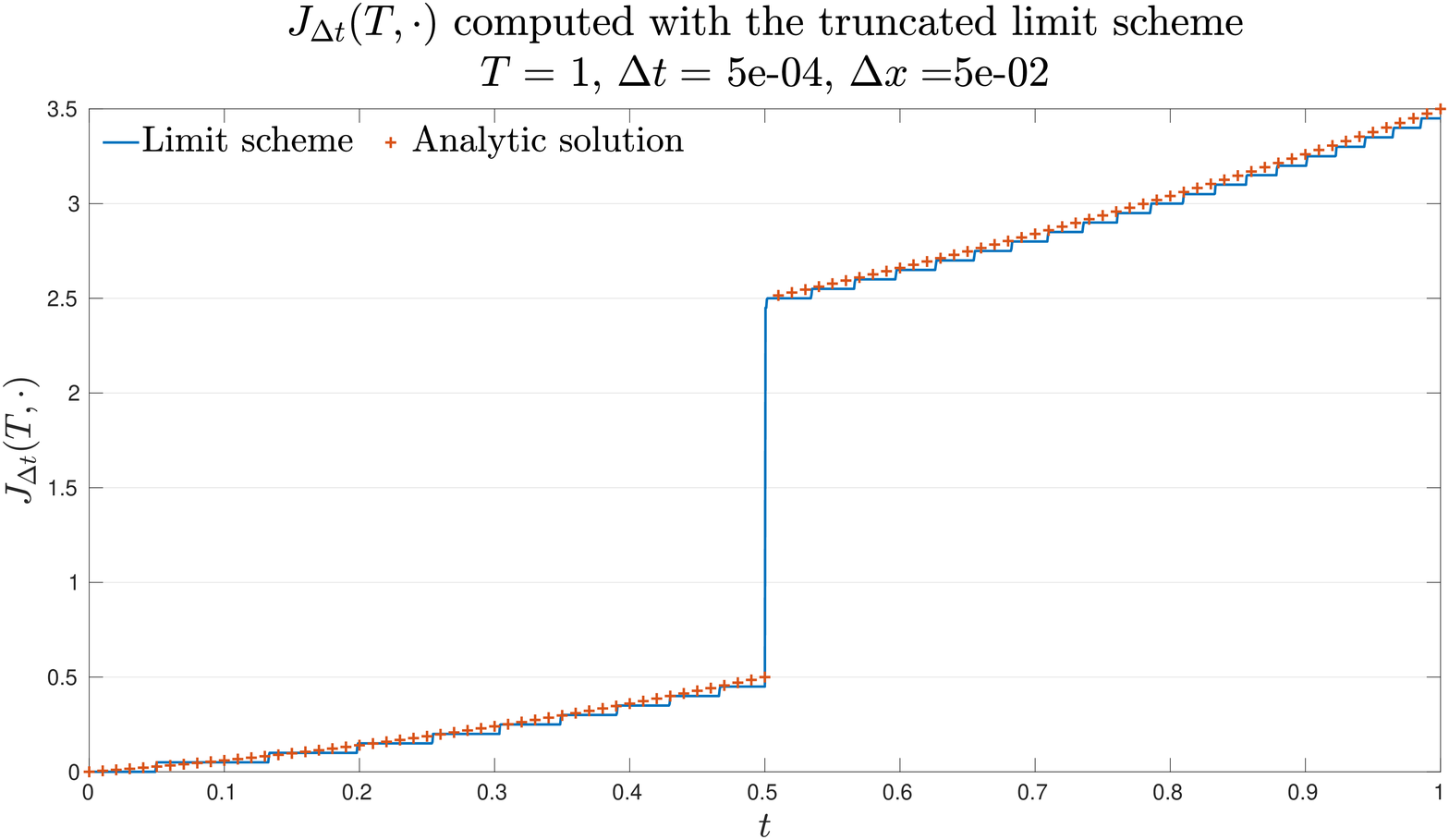}
\end{tabular}
\caption{Comparison between scheme \eqref{scheme:limit} and analytic solution. Left: $v$. Right: $J$.
Parameters: $T=1$, $\dx=5\cdot10^{-2}$, $\dt= 5\cdot 10^{-4}$, $v^\tin$ defined in \eqref{eq:PrecisionLimitScheme_vin}, and $R$ in \eqref{eq:PrecisionLimitScheme_R}.
}
\label{fig:PrecisionLimitScheme_1}
\end{center}
\end{figure}

Using its version with approximations at boundaries, we test the accuracy of \eqref{scheme:limit} with parameters which does not satisfy exactly the hypotheses \eqref{assumption:R_extrema}-\eqref{assumption:R_bounded_decreasing}-\eqref{assumption:u_increasing_infty}-\eqref{assumption:u0Lipschitz}. Indeed, we consider 
\begin{equation}
\label{eq:PrecisionLimitScheme_vin}
 v^\tin = \min\left( x^2; (x-\alpha)^2 + \delta\right),
\end{equation}
with $\alpha=2$, $\delta = 1$, and
\begin{equation}
 \label{eq:PrecisionLimitScheme_R} 
 R(x,I)=x-I.
\end{equation}
The solution of \eqref{eq:limit} is analytically known using these parameters, see \cite{BarlesPerthame}. Moreover, this explicit solution do not enjoy more regularity than what is expected. Indeed, $v$ enjoys Lispschitz regularity but is not $\mathcal{C}^1$, and $J$ jumps at $t=1/2$. 
The results of scheme \eqref{scheme:limit} are displayed in Fig. \ref{fig:PrecisionLimitScheme_1} together with the analytic solution. 
The agreement is visually very good, including the discontinuity of $J$ which is captured at the correct time point.

\begin{figure}[!ht]
\begin{center}
\begin{tabular}{@{}c@{}c@{}}
\includegraphics[width=0.5\textwidth]{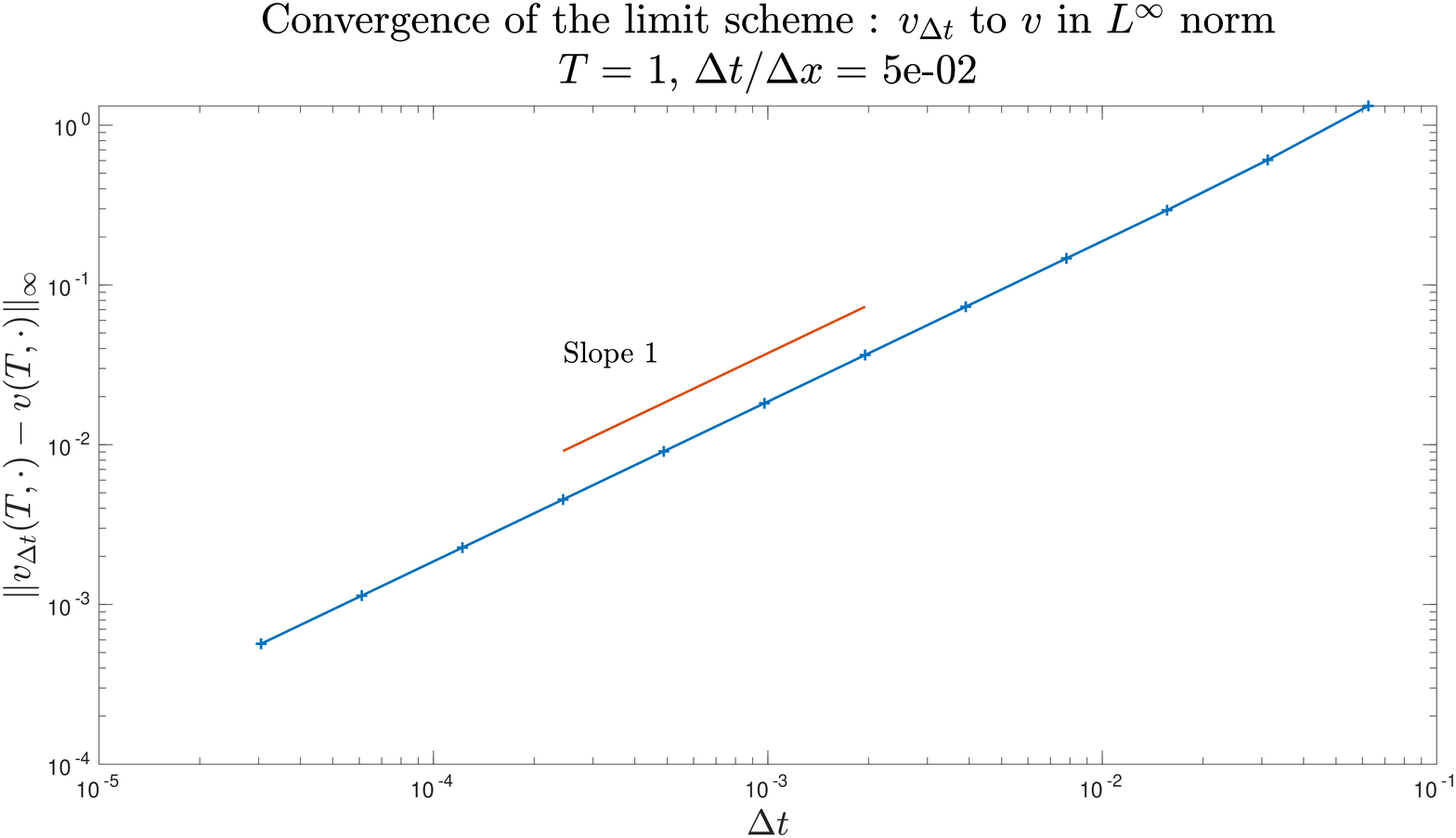}    &
\includegraphics[width=0.5\textwidth]{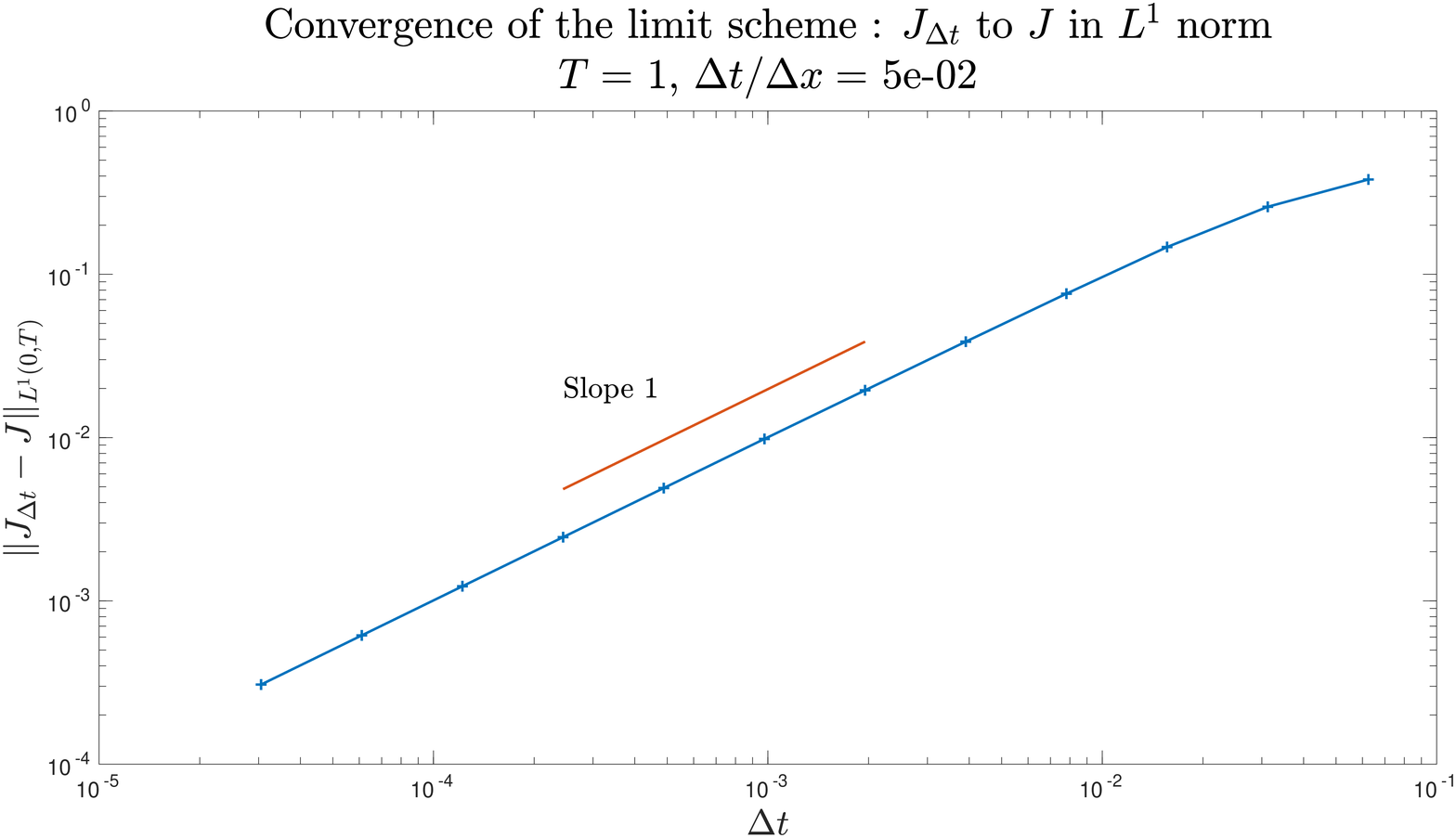}
\end{tabular}
\caption{Convergence rate of scheme \eqref{scheme:limit}. Left: $\|v_\dt(T,\cdot) - v(T,\cdot)\|_\infty$ as a function of $\dt$, with $\dt/\dx$ fixed (logarithmic scale). Right: $\|J_\dt - J\|_{\LL^1(0,T)}$ as a function of $\dt$, with $\dt/\dx$ fixed (logarithmic scale).
Parameters: $T=1$, $\dt/\dx=5\cdot10^{-2}$, $v^\tin$ defined in \eqref{eq:PrecisionLimitScheme_vin}, and $R$ in \eqref{eq:PrecisionLimitScheme_R}.
}
\label{fig:PrecisionLimitScheme_2}
\end{center}
\end{figure}

Although Prop. \ref{prop:scheme_0} states the convergence of \eqref{scheme:limit} to the solution of \eqref{eq:limit} when $\dt$ and $\dx$ go to $0$ with $\dt/\dx$ fixed, it does not give any convergence rate. Indeed, the lack of regularity of the solutions of \eqref{eq:limit} makes this problem difficult to address theoretically. 
To bypass this issue, we proposed a proof  based on compactness arguments and on an appropriate regularization of $J$. However, quantitative estimates cannot be expected using such arguments.
We propose a numerical study of the rate of convergence of \eqref{eq:limit} in Fig. \ref{fig:PrecisionLimitScheme_2}. For this numerical test, we compare the functions $v_\dt$ and $J_\dt$ computed with \eqref{scheme:limit}, to the solution $v$ and $J$ of \eqref{eq:limit} analytically computed in \cite{BarlesPerthame}. As in Prop. \ref{prop:scheme_0}, we fix $\dt/\dx$ and we make $\dt$ go to $0$.  The comparison is done in $\LL^\infty$ for $v_\dt(T,\cdot)-v(T,\cdot)$, while $J_\dt-J$ is estimated in $\LL^1(0,T)$ norm. The error is displayed in logarithmic scale. Remark that the numerical convergence rate of scheme \eqref{scheme:limit} is $1$, both for $v_\dt$ and $J_\dt$.

\subsection{Uniform accuracy of \eqref{scheme:epsilon}}
\label{sec:UA}

\begin{figure}[!ht]
\begin{center}
\begin{tabular}{@{}c@{}c@{}}
\includegraphics[width=0.5\textwidth]{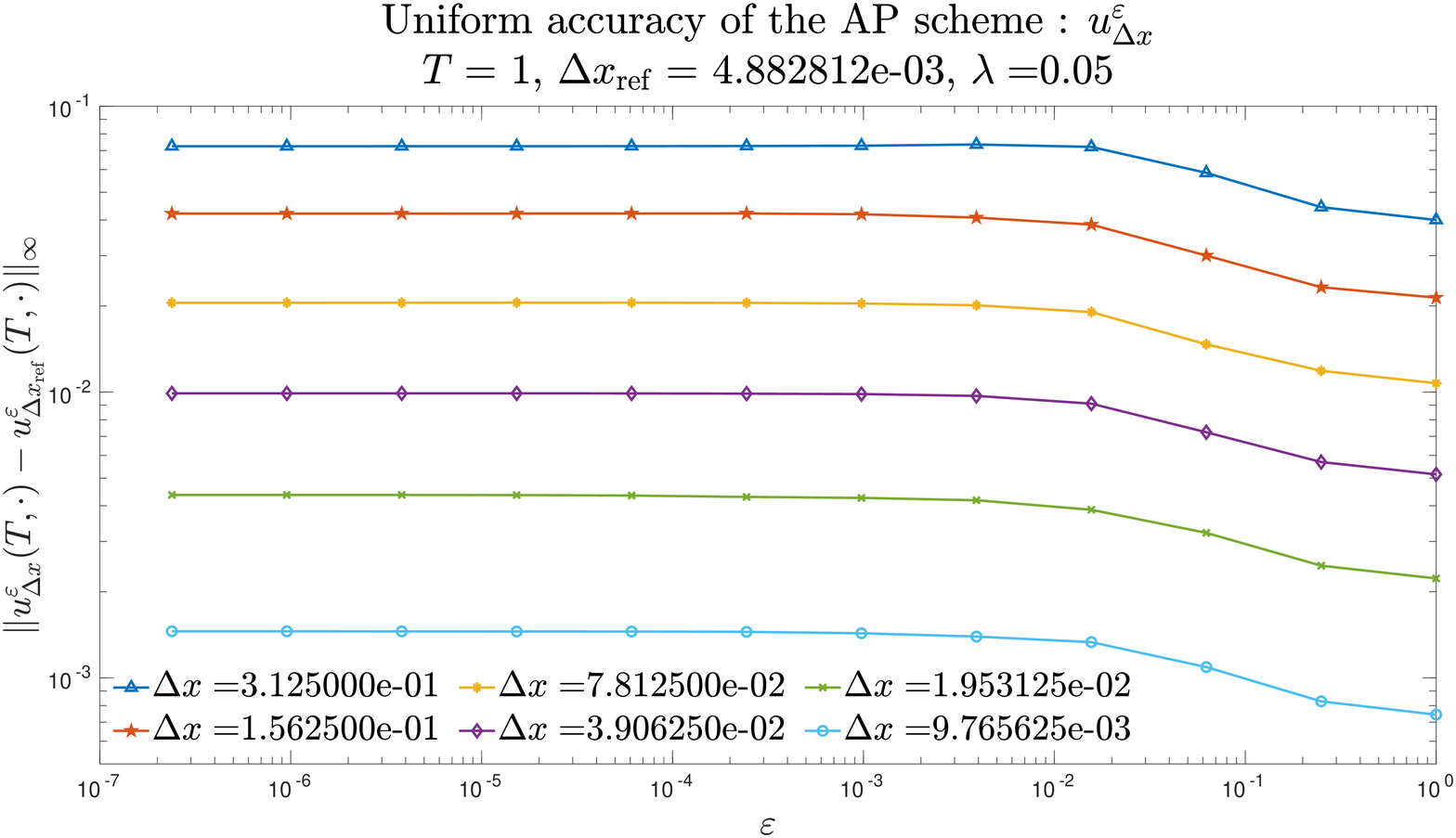}    &
\includegraphics[width=0.5\textwidth]{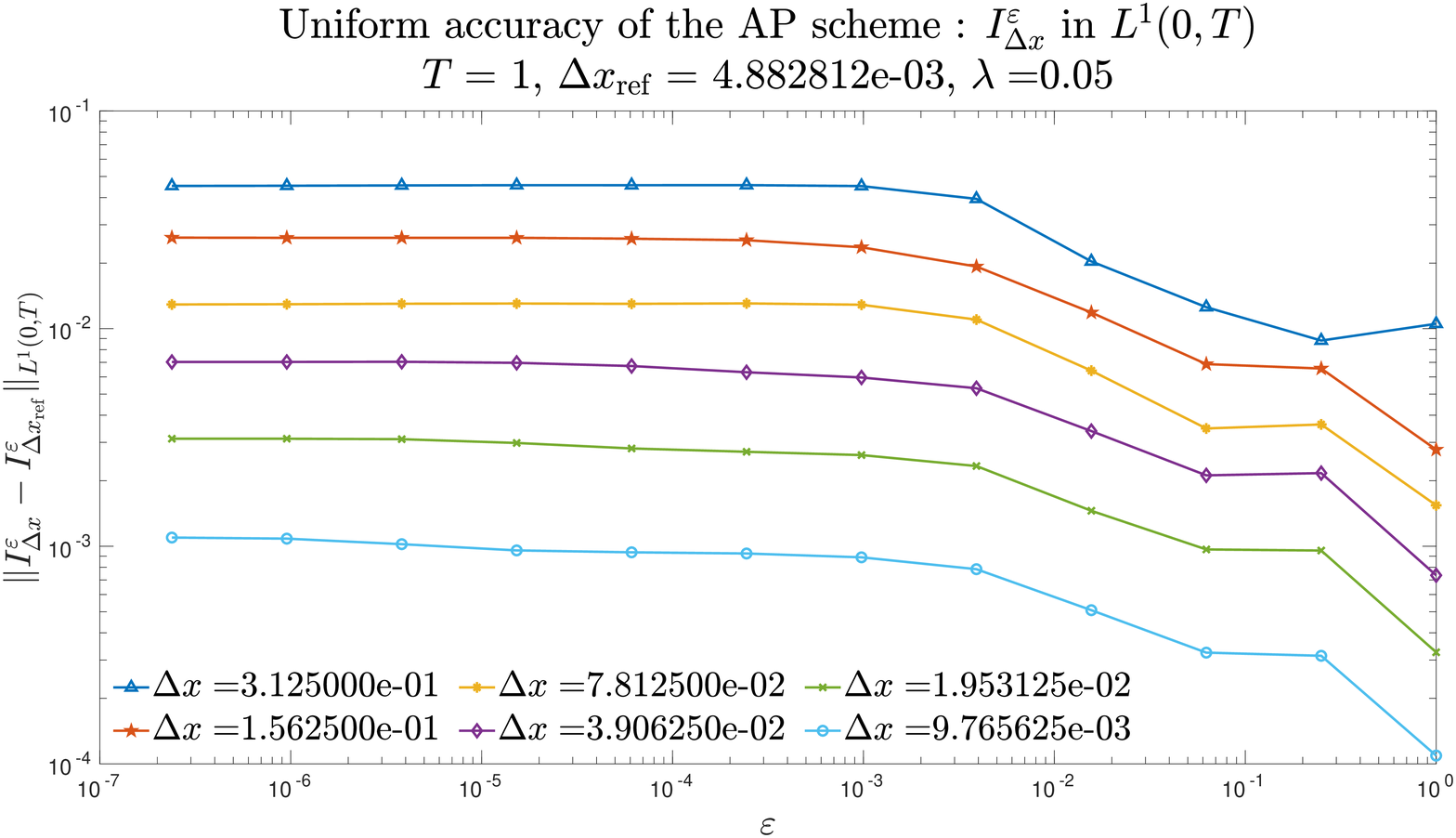}
\end{tabular}
\caption{Uniform accuracy of \eqref{scheme:epsilon}. Left:  $\|u_\dx^\ep(T,\cdot)-u_{\dx_{\mathrm{ref}}}^\ep(T,\cdot)\|_{\LL^\infty}$ for a series of $\dx$ and as functions of $\ep$ (logarithmic scale). Right: $\|I_\dx^\ep-I_{\dx_{\mathrm{ref}}}^\ep\|_{\LL^1(0,T)}$ for a series of $\dx$ and as functions of $\ep$ (logarithmic scale).
Parameters: $T=1$, $\lambda=5\cdot10^{-2}$, $u^\tin$ defined in \eqref{eq:u_in}, and $R$ in \eqref{eq:R}.
}
\label{fig:UA}
\end{center}
\end{figure}

In this section, the uniform accuracy of the scheme \eqref{scheme:epsilon}, in its version with approximations at boundaries, is tested.  Prop. \ref{prop:scheme_ep} establishes that, for all $\ep>0$, \eqref{scheme:epsilon} converges with rate $C(\ep)\left(|\ln(\dt)|\dt + \dx\right)$, with $\dt$ and $\dx$ satisfying \eqref{eq:epfixed_CFL}, and where $C(\ep)$ depends on $u^\tin$, $T$, and $\ep$. As it is emphasized in Remark \ref{rmq:UA}, this proposition does not give any clue on the order of the scheme uniformly in $\ep$, since $C(\ep)$ is expected to go to $+\infty$ when $\ep\to 0$. However, thanks to the stability properties of scheme \eqref{scheme:epsilon} stated in Prop. \ref{prop:AP}, a better behavior can be suspected. The uniform accuracy of scheme \eqref{scheme:epsilon} is tested by computing the results of \eqref{scheme:epsilon} for a series of $\ep$ and $\dx$. 
The solution of the corresponding scheme 
 will be denoted $u^\ep_\dx$ in what follows. Once $\dx$ is given, $\dt$ is fixed by 
$ \dt= \lambda\min(\dx; \dx^2/\ep)$, with $\lambda$ such that \eqref{eq:epfixed_CFL} and \eqref{eq:CFL_AP} hold. These $u^\dx_\ep$ are then compared to a reference solution. However, contrary to the previous section, no analytic solution of \eqref{eq:u_epsilon} is known, to the best of our knowledge, so that the reference solution has to be itself an approximation. A $\dx_{\mathrm{ref}}$ is introduced, smaller than all the $\dx$ previously considered, and $u_{\dx_{\mathrm{ref}}}^\ep$ is computed for all the $\ep$ considered. The $\LL^\infty$ norm of  $u_{\dx_{\mathrm{ref}}}(T,\cdot)^\ep-u_\dx^\ep(T,\cdot)$ 
is then computed for all $\ep$ and $\dx$ considered, and they are presented as functions of $\ep$, on the left-hand side of
Fig. \ref{fig:UA} in logarithmic scale. Similarly, the right-hand side of Fig. \ref{fig:UA} displays the $\LL^1(0,T)$ norm of $I_\dx^\ep - I_{\dx_{\mathrm{ref}}}^\ep$, as functions of $\ep$ in logarithmic scale. Remark that, in both cases, these error curves are stratified, meaning that the approximation error in scheme \eqref{scheme:epsilon} is uniformly bounded with respect to $\ep$ when the discretization is fixed.

\begin{figure}[!ht]
\begin{center}
\begin{tabular}{@{}c@{}c@{}}
\includegraphics[width=0.5\textwidth]{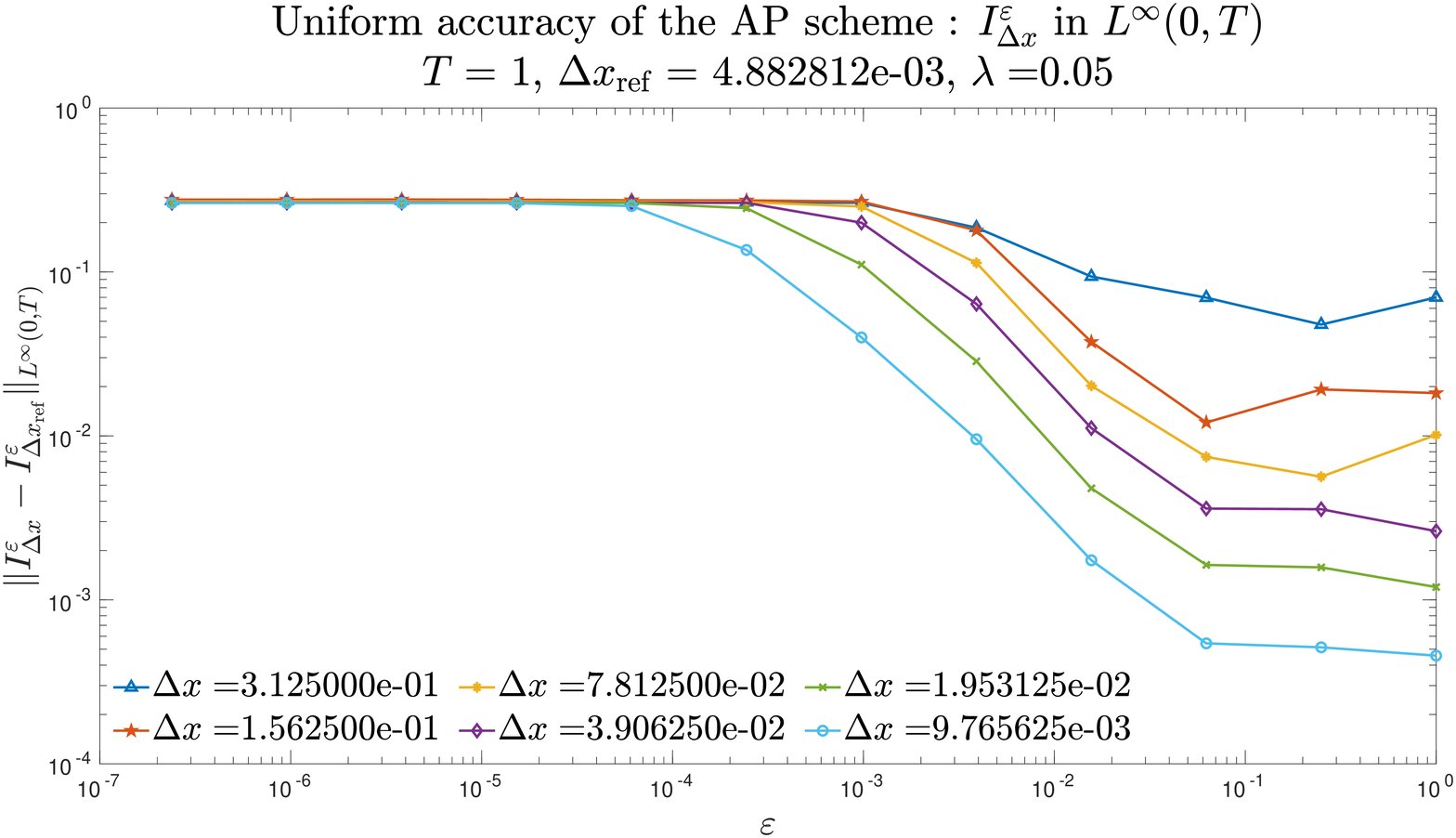}    &
\includegraphics[width=0.5\textwidth]{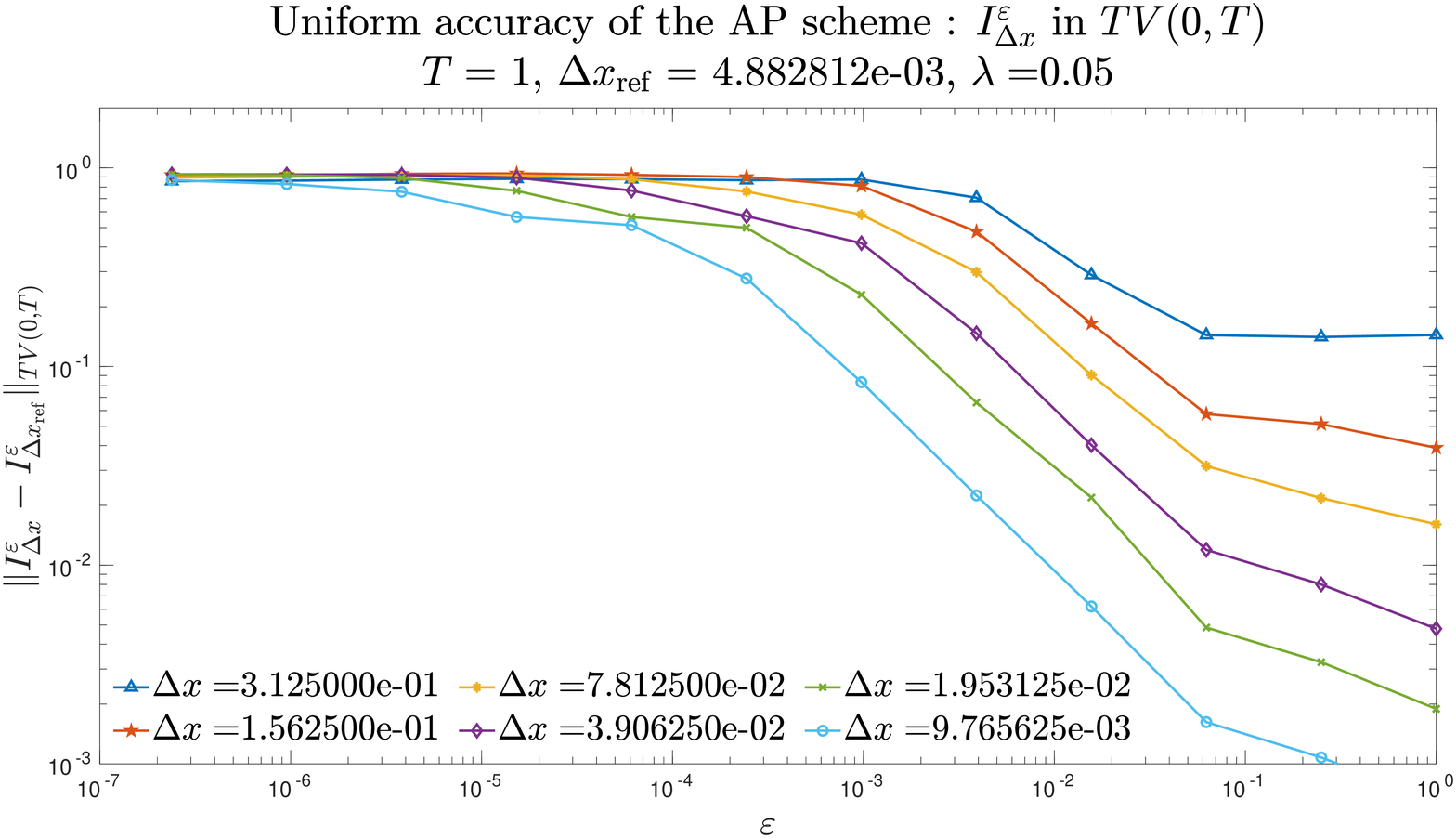}
\end{tabular}
\caption{Uniform accuracy test  for $I_\ep$ computed with \eqref{scheme:epsilon}: $\|I^\dx_\ep-I^{\dx_{\mathrm{ref}}}_\ep\|$ for a series of $\dx$ and as functions of $\ep$ (logarithmic scale). Left: $\LL^\infty(0,T)$ norm. Right: $TV(0,T)$ seminorm.
Parameters: $T=1$, $\lambda=5\cdot10^{-2}$, $u^\tin$ defined in \eqref{eq:u_in}, and $R$ in \eqref{eq:R}.}
\label{fig:Not_UA}
\end{center}
\end{figure}

The numerical tests above suggest that scheme \eqref{scheme:epsilon} enjoys uniform accuracy with respect to $\ep$ in $\LL^\infty$ norm for $u_\dx^\ep(T,\cdot)$ and in $\LL^1(0,T)$ norm for $I_\dx^\ep$. 
However, the lack of regularity of the solutions of \eqref{eq:u_epsilon} strongly influences the accuracy of the numerical resolution. To emphasize on this fact, 
remark that the uniform accuracy of \eqref{scheme:epsilon} is not true for $I_\dx^\ep$ in $\LL^\infty(0,T)$ nor in the total variation seminorm, denoted  $TV(0,T)$ in what follows. Indeed, the $\LL^\infty(0,T)$ norm and $TV(0,T)$ seminorm of $I_\dx^\ep - I_{\dx_{\mathrm{ref}}}^\ep$, as functions of $\ep$ in logarithmic scale are displayed in  Fig. \ref{fig:Not_UA}. Contrary to Fig. \ref{fig:UA}, the error curves are not stratified, and one can remark that 
 \[
  \sup\limits_\ep \left\|
  I_\dx^\ep - I_{\dx_{\mathrm{ref}}}^\ep
  \right\|_{\LL^\infty(0,T)} \underset{\dx\to 0} \nrightarrow 0,\;\;\;\;\text{and}\;\;\;\; 
  \sup\limits_\ep \left\|
  I_\dx^\ep - I_{\dx_{\mathrm{ref}}}^\ep
  \right\|_{TV(0,T)} \underset{\dx\to 0} \nrightarrow 0,
 \]
 meaning that \eqref{scheme:epsilon} does not enjoy uniform accuracy for $I_\dx^\ep$ in $\LL^\infty$ norm and $TV$ seminorm. The fact that these norms are poorly adapted to the study of the convergence of $I^\dx_\ep$ can be understood considering the jumps. For small $\ep$, $I_\ep$ is close to the discontinuous function $J$, so that the jumps are visually well approximated. The comparison between $I_\ep^\dx$ and $I_\ep^{\dx_{\mathrm{ref}}}$ is also good at first sight, but jumps may not be exactly simultaneous, making the difference $I_\ep^\dx-I_\ep^{\dx_{\mathrm{ref}}}$ have a thin peak around the jump. Such a peak is small in $L^1$ norm, but not in $L^\infty$ or $TV$.

\subsection{Extension to higher dimensions}
\label{sec:dim2}

Problems \eqref{eq:u_epsilon} and \eqref{eq:limit} are well-posed in any finite dimension $d$, but dimension $1$ was chosen for the presentation and the study of schemes \eqref{scheme:epsilon} and \eqref{scheme:limit} in this paper. 
However, schemes \eqref{scheme:epsilon}-\eqref{scheme:limit}  can be generalized to any dimension, and all the results of this paper still hold when $d\in\N^*$, the proofs being done exactly the same way but with heavier notations due to multi-indices.

\begin{figure}[!ht]
\begin{center}
\begin{tabular}{@{}c@{}c@{}}
\includegraphics[width=0.5\textwidth]{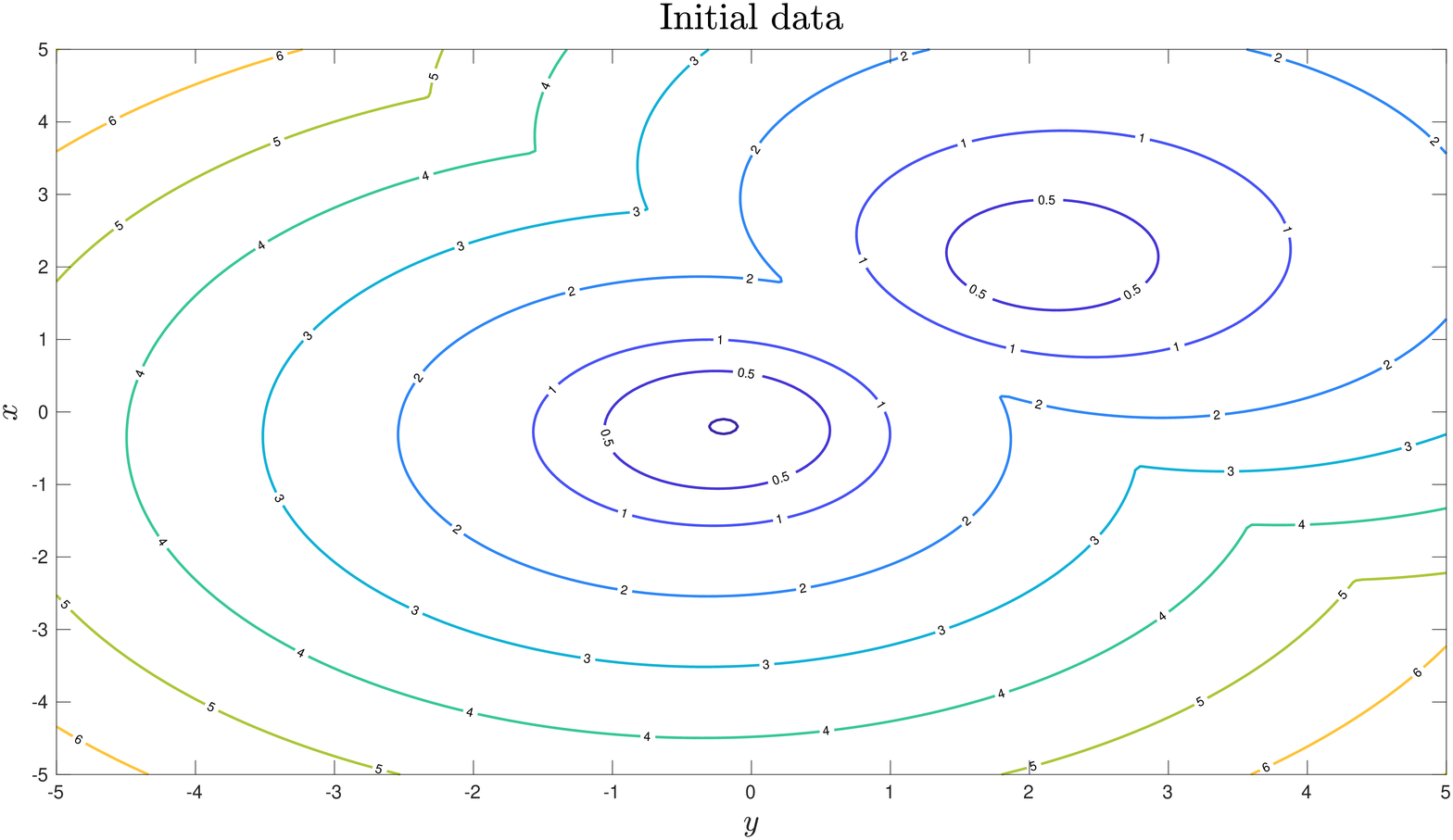}    &
\includegraphics[width=0.5\textwidth]{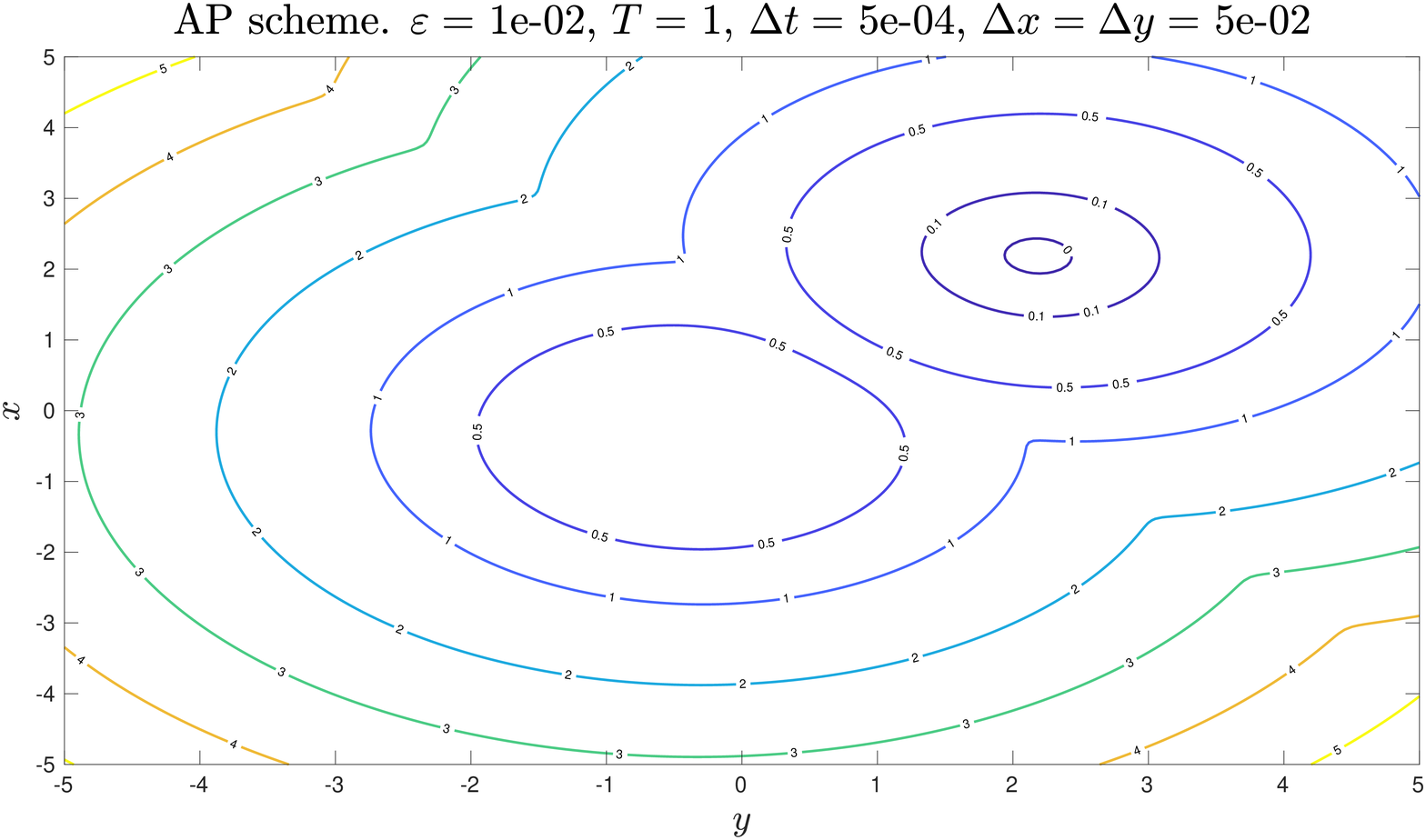} \\
\includegraphics[width=0.5\textwidth]{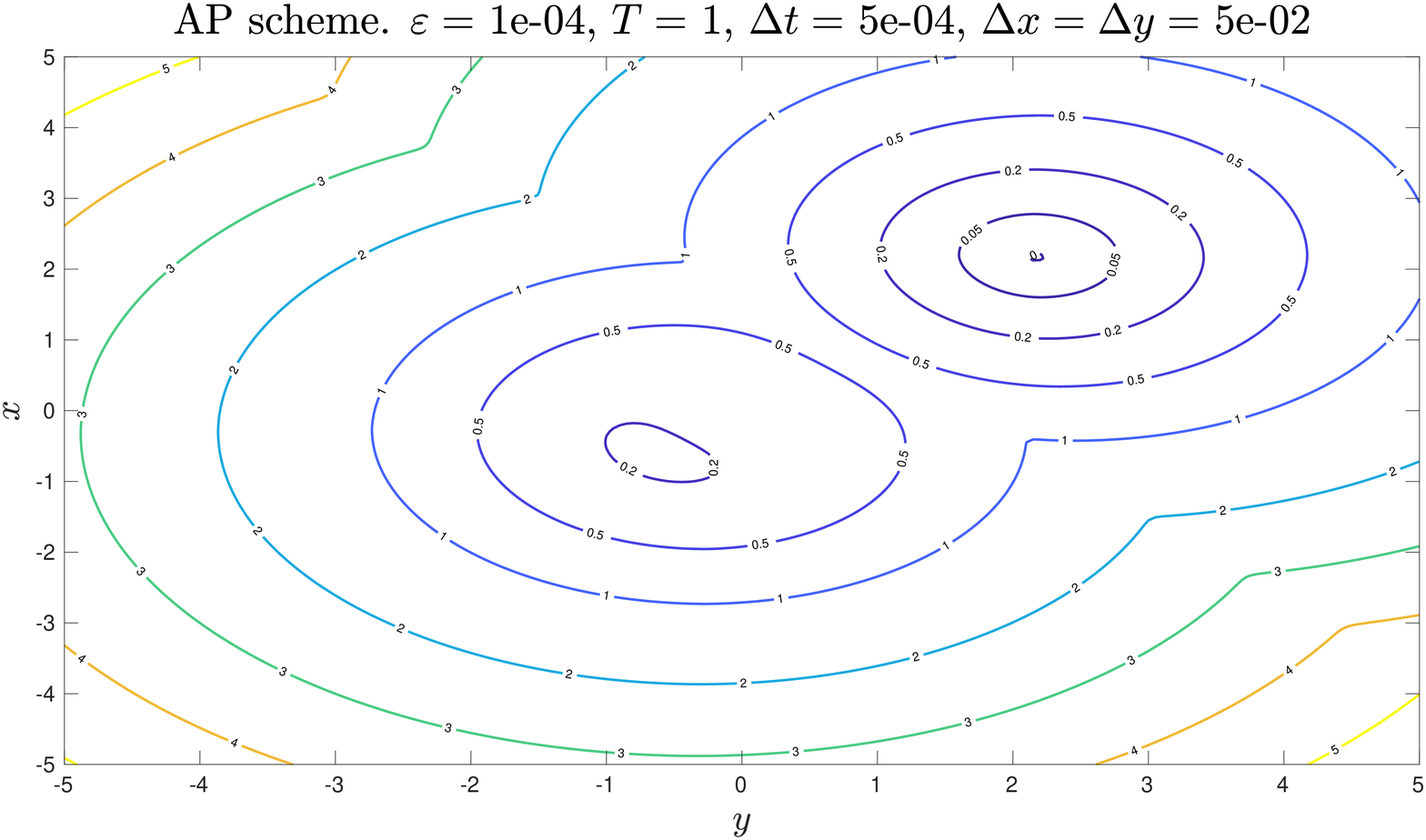}    &
\includegraphics[width=0.5\textwidth]{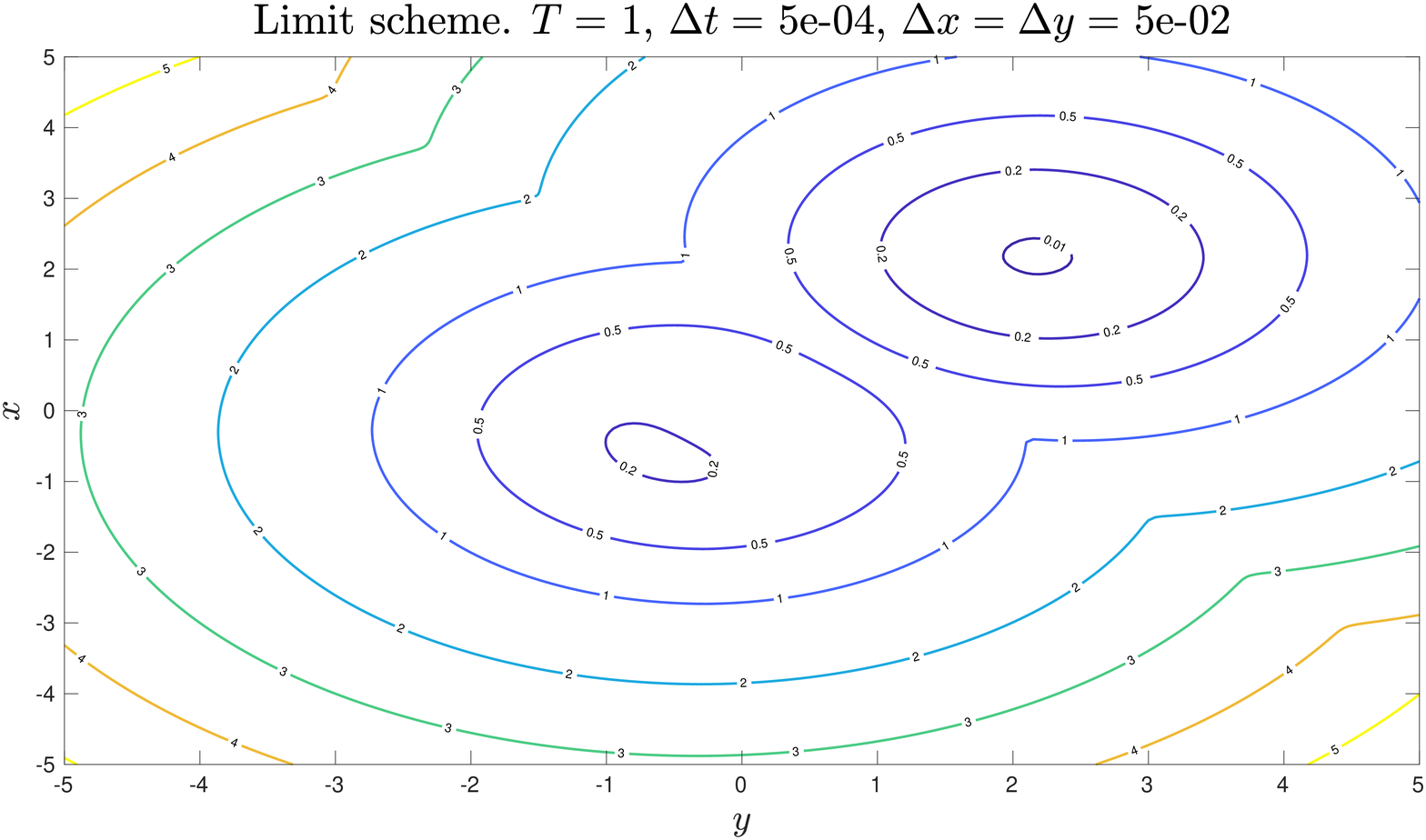}
\end{tabular}
\caption{Test with $d=2$. Top left: $u_\ep^\tin$. Top right: $u^\ep_\dt$ computed with \eqref{scheme:epsilon_dim2} and $\ep= 10^{-2}$. Bottom left:  $u^\ep_\dt$ computed with \eqref{scheme:epsilon_dim2} and $\ep= 10^{-4}$. Bottom right: $v_\dt$ computed with \eqref{scheme:limit_dim2}.
Parameters: $T=1$, $\dt=5\cdot10^{-4}$, $\dx=\dy=5\cdot10^{-2}$, $u_\ep^\tin$ and $v^\tin$ defined in \eqref{eq:u_in_dim2}, and $R$ in \eqref{eq:R_dim2}.
}
\label{fig:Dim2_u}
\end{center}
\end{figure}

\begin{figure}[!ht]
 \centering
\includegraphics[width=0.6\textwidth]{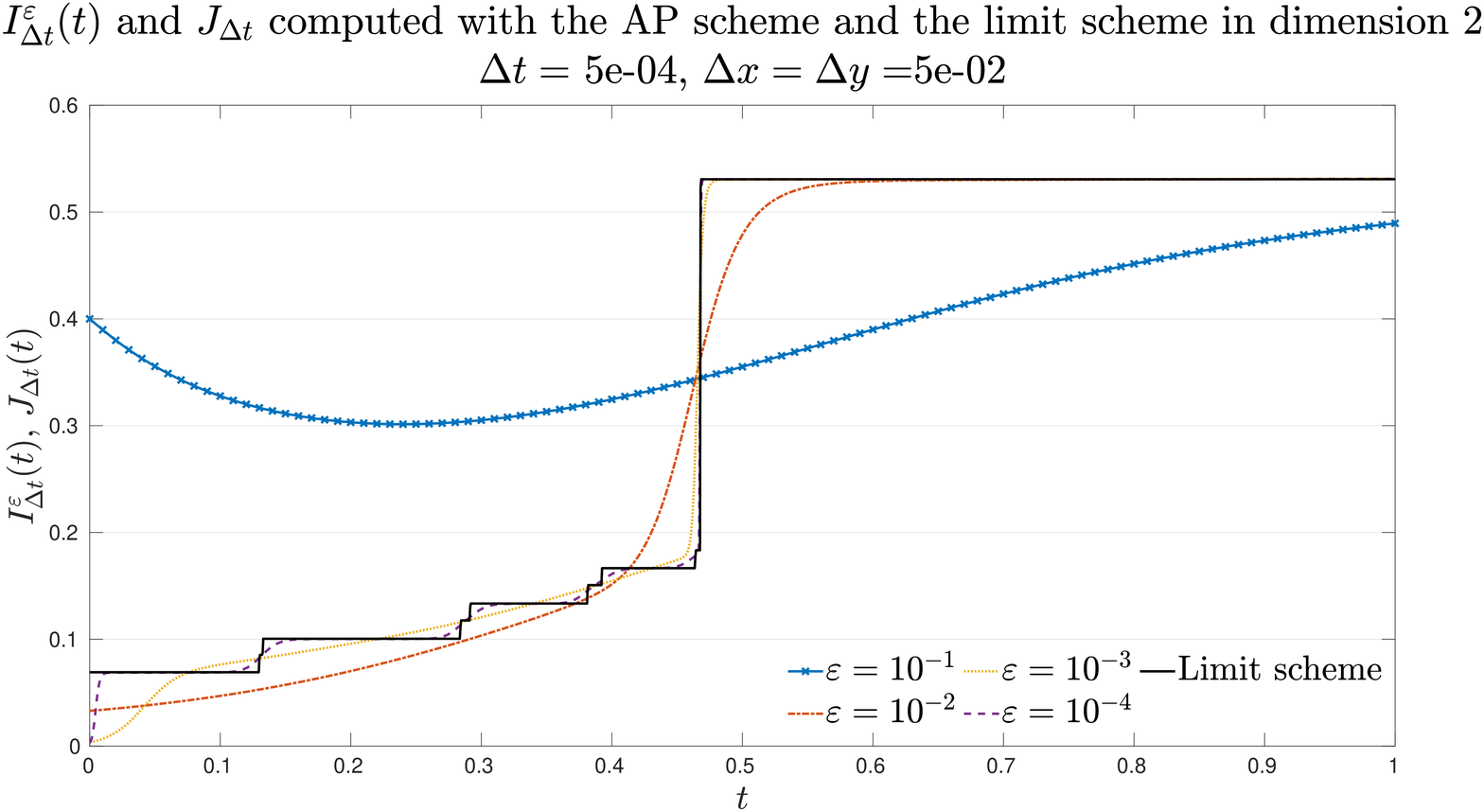}
\caption{Test with $d=2$: $I^\ep_\dt$ computed with \eqref{scheme:epsilon_dim2} for a series of $\ep$, and $J_\dt$ computed with \eqref{scheme:limit_dim2}.
Parameters: $T=1$, $\dt=5\cdot10^{-4}$, $\dx=\dy=5\cdot10^{-2}$, $u_\ep^\tin$ and $v^\tin$ defined in \eqref{eq:u_in_dim2}, and $R$ in \eqref{eq:R_dim2}.
}
\label{fig:dim2_I}
\end{figure}

We detail here the adaptation of schemes \eqref{scheme:epsilon} and \eqref{scheme:limit} in dimension $d=2$, and we provide some numerical tests to highlight the asymptotic-preserving property. The generalization to any dimension is straightforward. As in Section \ref{sec:results}, define $T$, $N_t$ and $\dt$ for the time discretization. Two trait steps are now needed, denoted $\dx$ and $\dy$ in what follows, and two trait grids are defined, namely
$ x_i=x_0+i\dx$ ($i\in\Z$), and $y_j=y_0+y_j\dy$  ($j\in\Z$), 
where $x_0,y_0\in\R$ are given. Let $n\in\ccl 0,N_t-1\ccr$, and $i,j\in\Z$. The schemes are given by
\begin{equation}
 \label{scheme:epsilon_dim2}
 \tag{$S_\ep^{d=2}$}
 \left\{ 
\begin{array}{r l}
 \ds\frac{u^{n+1}_{i,j}-u^n_{i,j}}{\dt} & \ds +
 H\left( \frac{u^n_{i,j}-u^n_{i-1,j}}{\dx}, \frac{u^n_{i+1,j}-u^n_{i,j}}{\dx} \right) 
 + H\left( \frac{u^n_{i,j}-u^n_{i,j-1}}{\dy}, \frac{u^n_{i,j+1}-u^n_{i,j}}{\dy} \right) \vspace{4pt} \\
 &\ds = \ep \frac{u^n_{i+1,j}-2u^n_{i,j} + u^n_{i-1,j}}{\dx^2} 
 + \ep \frac{u^n_{i,j+1}-2u^n_{i,j} + u^n_{i,j-1}}{\dy^2} - R\left(x_i,y_j,I^{n+1}\right) \vspace{4pt} \\
 \ds I^{n+1} &\ds= \dx\dy \sum\limits_{(i,j)\in\Z^2} \psi(x_i,y_j) \e^{-u^{n+1}_{i,j}/\ep},
\end{array}
 \right.
\end{equation}
and
\begin{equation}
 \label{scheme:limit_dim2}
 \tag{$S_0^{d=2}$}
 \left\{ 
\begin{array}{r l}
 \ds\frac{v^{n+1}_{i,j}-v^n_{i,j}}{\dt} & \ds +
 H\left( \frac{v^n_{i,j}-v^n_{i-1,j}}{\dx}, \frac{v^n_{i+1,j}-v^n_{i,j}}{\dx} \right) 
 + H\left( \frac{v^n_{i,j}-v^n_{i,j-1}}{\dy}, \frac{v^n_{i,j+1}-v^n_{i,j}}{\dy} \right) \vspace{4pt} \\
 &\ds = - R\left(x_i,y_j,I^{n+1}\right) \vspace{4pt} \\
 \ds \min\limits_{(i,j)\in\Z^2} v^{n+1}_{i,j} & \ds = 0,
\end{array}
 \right.
\end{equation}
where $H$ is defined in \eqref{eq:H}. They both can be implemented on a truncated domain, with or without approximations at boundaries, as presented above. The following tests use the version with approximation at boundaries and grids of constant size. Denoting $X=(x,y)\in\R^2$, we consider $u^\tin$, $v^\tin$ adapted from \eqref{eq:u_in},
\begin{equation}
\label{eq:u_in_dim2}
 u^\tin(X) = v^\tin(X) = \frac{\min\left(|X-\beta|^2; |X-\alpha|^2+ \delta\right)}{\sqrt{1+|X|^2}},
\end{equation}
with $\alpha=(2,2)$, $\beta=(-0.2,-0.2)$ and $\delta =1$. Similarly to \eqref{eq:R}, we define
\begin{equation}
\label{eq:R_dim2}
 R(X,I) = \e^{-I} \frac{|X|^2}{1+|X|^2}- I.
\end{equation}
Note that in both cases, $|\cdot|$ stands for the Euclidean norm on $\R^2$.
Fig. \ref{fig:Dim2_u} displays level lines of $u_\ep^\tin$ defined in \eqref{eq:u_in_dim2}, of $u^\ep_\dt$ computed with \eqref{scheme:epsilon_dim2} for $\ep=10^{-2}$ and $\ep=10^{-4}$, and of $v_\dt$ computed with \eqref{scheme:limit_dim2}. When $\ep$ is small, $u^\ep_\dt$ is similar to $v_\dt$. 
Moreover, one can notice that the  minimum of $u_\ep^\tin$ has jumped from the bottom left local minimum to the top right one. Fig. \ref{fig:dim2_I} highlights the stability of the component $I^\ep_\dt$ in \eqref{scheme:epsilon_dim2} when $\ep\to 0$. Indeed, it goes to the component $J_\dt$ of \eqref{scheme:limit_dim2}, and has discontinuities in the asymptotic regime. More generally, all the properties discussed in dimension $1$ are still statisfied.

\section*{Conclusion}

In this paper, we proposed and analyzed an asymptotic-preserving scheme for  parabolic Lotka-Volterra equations \eqref{eq:u_epsilon}, which model the evolution of a population density. The scheme \eqref{scheme:epsilon} we proposed is converging 
for fixed $\ep>0$,
and enjoys stability properties in the asymptotics. Moreover,  the limiting scheme \eqref{scheme:limit} is converging towards the unique viscosity solution of the constrained Hamilton-Jacobi equation \eqref{eq:limit}, which describes the asymptotic regime. 

The key ingredients for the construction of the asymptotic-preserving scheme are the monotonicity and the implicit treatment of the constraint. Thanks to these properties, the convergence of the scheme \eqref{scheme:epsilon} is proved, and so is its asymptotic behavior in the vanishing $\ep$ limit. The convergence of the limit scheme is based on compactness arguments, and once again  on the monotonicity of the scheme. It is indeed a usual property required for non-diffusive schemes for Hamilton-Jacobi equations. However, because of  the lack of regularity of the Lagrange multiplier associated with the non-negativity constraint, the scheme has to be regularized to prove its convergence. Eventually, the properties of the schemes have been discussed through numerical tests. 
Both \eqref{scheme:epsilon} and \eqref{scheme:limit} can also be generalized to any dimension. Moreover, numerical tests suggest that \eqref{scheme:epsilon} enjoys uniform accuracy in appropriate discrete function spaces, meaning that its precision is independent of $\ep$. 

A natural continuation of this work would be  the study of an asymptotic-preserving scheme for  integral Lotka-Volterra equations. This question will be adressed in a future work.

\bibliographystyle{plain}
\bibliography{Biblio}

\end{document}